\documentstyle[12pt,amssymb]{amsart}

\newcounter{lm}[section]
\newcounter{thm}[section]
\newcounter{prop}[section]
\newcounter{rem}[section]
\newcounter{cor}[section]

\newcounter{ex}[section]

\newcounter{defin}[section]
\newcounter{conv}[section]
\begin{document}
\newcommand{\thm}{\refstepcounter{thm} {\bf Theorem \arabic{section}.%
\arabic{thm}.} }

\renewcommand{\thethm}{\thesection.\arabic{thm}}

\newcommand{\nthm} [1] {\refstepcounter{thm} {\bf Theorem 
\arabic{section}.%
\arabic{thm}:} #1}

\newcommand{\nnthm} [1] {\refstepcounter{thm} {\bf Theorem 
\arabic{section}.%
\arabic{thm}} #1}

\newcommand{\nprop} [1] {\refstepcounter{prop} {\bf Proposition 
\arabic{section}.%
\arabic{prop}:} #1}

\newcommand{\prop} {\refstepcounter{prop}{\bf Proposition \arabic{section}.%
\arabic{prop}.} }

\renewcommand{\theprop}{\thesection.\arabic{prop}}

\newcommand{\cor} {\refstepcounter{cor}{\bf Corollary \arabic{section}.%
\arabic{cor}.} }

\newcommand{\ncor} [1] {\refstepcounter{cor} {\bf Corollary
\arabic{section}.%
\arabic{cor}:} #1} 

\renewcommand{\thecor}{\thesection.\arabic{cor}}

\newcommand{\lemma}{\refstepcounter{lm}{\bf Lemma \arabic{section}.%
\arabic{lm}.} }

\newcommand{\nlemma} [1] {\refstepcounter{lm} {\bf Lemma
\arabic{section}.%
\arabic{lm}:} #1}

\renewcommand{\thelm}{\thesection.\arabic{lm}}

\newcommand{\ex}{\refstepcounter{ex}{\bf Example \arabic{section}.%
\arabic{ex}.} }

\newcommand{\nex} [1] {\refstepcounter{ex} {\bf Example
\arabic{section}.%
\arabic{ex}:} #1} 

\renewcommand{\theex}{\thesection.\arabic{ex}}

\newcommand{\defin}{\refstepcounter{defin}{\bf Definition \arabic{section}.%
\arabic{defin}.} }

\renewcommand{\thedefin}{\thesection.\arabic{defin}}

\newcommand{\rem}{\refstepcounter{rem}{\bf Remark \arabic{section}.%
\arabic{rem}.} }

\newcommand{\nrem}{\refstepcounter{rem}{\bf \arabic{section}.%
\arabic{rem}.} }

\renewcommand{\therem}{\thesection.\arabic{rem}}

\newcommand{\rems}{{\bf Remarks }}

\newcommand{\conv}{\refstepcounter{conv}{\bf Convention \arabic{section}.%
\arabic{conv}.} }

\renewcommand{\theconv}{\thesection.\arabic{conv}}

\def\qed{\hfill $\Box$}
\newcommand{\vi}{{\varphi}}
\newcommand{\C}{{\bf C}}
\newcommand{\pr}{{\bf P}}
\newcommand{\R}{{\bf R}}
\newcommand{\Z}{{\bf Z}}
\newcommand{\Q}{{\bf Q}}
\newcommand{\N}{{\bf N}}
\newcommand{\tf}{{\tilde {f}}}
\newcommand{\tg}{{\tilde {g}}}
\newcommand{\th}{{\tilde {h}}}
\newcommand{\tp}{{\tilde {p}}}
\newcommand{\tx}{{\tilde {x}}}
\newcommand{\ty}{{\tilde {y}}}
\newcommand{\tz}{{\tilde {z}}}
\newcommand{\tu}{{\tilde {u}}}
\newcommand{\tq}{{\tilde {q}}}
\newcommand{\tv}{{\tilde {v}}}
\newcommand{\ts}{{\tilde {s}}}
\newcommand{\te}{{\tilde {e}}}
\newcommand{\tl}{{\tilde {l}}}
\newcommand{\tU}{{\tilde {U}}}
\newcommand{\tW}{{\tilde {W}}}
\newcommand{\tZ}{{\tilde {Z}}}
\newcommand{\tV}{{\tilde {V}}}
\newcommand{\tS}{{\tilde {S}}}
\newcommand{\tL}{{\tilde {L}}}
\newcommand{\tN}{{\tilde {N}}}
\newcommand{\tY}{{\tilde {Y}}}

\newcommand{\tvar}{{\tilde {\varphi}}}
\newcommand{\tsi}{{\tilde {\sigma}}}
\newcommand{\tpsi}{{\tilde {\psi}}}
\newcommand{\tchi}{{\tilde {\chi}}}
\newcommand{\tE}{{\tilde {E}}}
\newcommand{\tA}{{\tilde {A}}}
\newcommand{\tH}{{\tilde {H}}}
\newcommand{\tI}{{\tilde {I}}}
\newcommand{\tC}{{\tilde {C}}}
\newcommand{\G}{{\Gamma}}
\newcommand{\bS}{{\bar {S}}} 
\newcommand{\bC}{{\bar {C}}}
\newcommand{\bY}{{\bar {Y}}} 
\newcommand{\bH}{{\bar {H}}}
\newcommand{\bx}{{\bar {x}}}
\newcommand{\bb}{{\bar {b}}}
\newcommand{\bp}{{\bar {p}}}
\newcommand{\bq}{{\bar {q}}}
\newcommand{\bE}{{\bar {E}}}
\newcommand{\bF}{{\bar {F}}}
\newcommand{\bV}{{\bar {V}}}
\newcommand{\bU}{{\bar {U}}}
\newcommand{\bZ}{{\bar {Z}}}
\newcommand{\baf}{{\bar {f}}}

\newcommand{\cD}{{\cal {D}}}
\newcommand{\cE}{{\cal {E}}}
\newcommand{\cH}{{\cal {H}}}
\newcommand{\cG}{{\cal {G}}}
\newcommand{\cC}{{\cal {C}}}
\newcommand{\cS}{{\cal {S}}}
\newcommand{\cT}{{\cal {T}}}
\newcommand{\cN}{{\cal {N}}}
\newcommand{\cJ}{{\cal {J}}}
\newcommand{\cP}{{\cal {P}}}
\newcommand{\cQ}{{\cal {Q}}}
\newcommand{\cF}{{\cal {F}}}
\newcommand{\cO}{{\cal {O}}}
\newcommand{\cL}{{\cal {L}}}
\newcommand{\cM}{{\cal {M}}}
\newcommand{\cI}{{\cal {I}}}
\newcommand{\cV}{{\cal {V}}}

\newcommand{\hvar}{{\widehat {\varphi }}}
\newcommand{\hpsi}{{\widehat {\psi }}}
\newcommand{\hx}{{\widehat {x}}}
\newcommand{\hz}{{\widehat {z}}}
\newcommand{\hy}{{\widehat {y}}}
\newcommand{\hh}{{\widehat {h}}}
\newcommand{\hf}{{\widehat {f}}}
\newcommand{\hg}{{\widehat {g}}}
\newcommand{\hr}{{\widehat {r}}}
\newcommand{\hp}{{\widehat {p}}}
\newcommand{\ha}{{\widehat {a}}}
\newcommand{\hc}{{\widehat {c}}}
\newcommand{\hJ}{{\widehat {J}}}
\newcommand{\hZ}{{\widehat {Z}}}
\newcommand{\hE}{{\widehat {E}}}
\newcommand{\hM}{{\widehat {M}}}
\newcommand{\hA}{{\widehat {A}}}
\newcommand{\hI}{{\widehat {I}}}
\newcommand{\hH}{{\widehat {H}}}
\newcommand{\hsi}{{\widehat {\sigma}}}
\newcommand{\si}{{\sigma}}
\newcommand{\htt}{{\widehat {t}}}
\newcommand{\hmu}{{\widehat {\mu}}}
\newcommand{\hR}{{\widehat {R}}}
\newcommand{\hX}{{\widehat {X}}}
\newcommand{\hY}{{\widehat {Y}}}
\newcommand{\widehatpar}{{\widehat {\partial}}}
\newcommand{\de}{{\partial}}
\newcommand{\hde}{{\widehat {\partial}}}

\newcommand{\gi}{{G_{\infty }}}
\newcommand{\bk}{{\bar {k}}}
\newcommand{\br}{{\bar {r}}}
\newcommand{\bz}{{\bar {z}}}
\newcommand{\ba}{{\bar {a}}}
\newcommand{\bR}{{\bar {R}}}

\newcommand{\dX}{{\dot {X}}}
\newcommand{\dM}{{\dot {M}}}
\newcommand{\dD}{{\dot {D}}}
\newcommand{\dE}{{\dot {E}}}

\newcommand{\dimc}{{{\dim}_{\C}\,}}
\newcommand{\De}{{Deg \,}}
\newcommand{\gC}{{C_{{\rm geom}}}}

\newcommand{\is}{{\simeq}}
\newcommand{\var}{{\varphi}}
\newcommand{\df}{\mbox{def}_{\partial}}
\newcommand{\fr} [1] {{F^{#1}A}}
\newcommand{\fro} [1] {{F^{#1}_0A}}
\newcommand{\po}{{\pi_0}}
\newcommand{\id}{{\bf 1}}
\newcommand{\pe}{{\pi_1}}
\newcommand{\supp}{{\rm supp ~}}
\newcommand{\const}{{\rm const}}
\newcommand{\tors}{{\rm Tors ~}}
\newcommand{\Tor}{{\rm Tor ~}}
\newcommand{\Pic}{{\rm Pic ~}}
\newcommand{\Ker}{{\rm Ker ~}}
\newcommand{\Per}{{\rm Per ~}}
\newcommand{\ML}{{\rm ML}}
\newcommand{\AK}{{\rm AK}} 
\newcommand{\Dk}{{\rm Dk}}
\newcommand{\LND}{{\rm LND}}
\newcommand{\LNDG}{{\rm LND_{gr}}}
\newcommand{\GL}{{\rm GL}}
\newcommand{\Gr}{{\rm Gr ~}}
\newcommand{\gr}{{\rm gr }}
\newcommand{\Div}{{\rm Div ~}}
\newcommand{\proof} {{\bf Proof. ~}}
\newcommand{\loc}{{S^{-1}}}

\newcommand{\chr} {{\rm char ~}}
\newcommand{\grad} {{\rm grad ~}}
\newcommand{\spec} {{\rm spec ~}}
\newcommand{\reg} {{\rm reg ~}}
\newcommand{\dt} {{\rm det ~}}

\newcommand{\dg} {{\rm deg ~}}
\newcommand{\trdg} {{\rm tr.deg ~}}
\newcommand{\dgv} {{\rm deg}_{\varphi} ~}
\newcommand{\dgy} {{\rm deg}_y ~}
\newcommand{\dge} {{\rm deg}_1 ~}
\newcommand{\dgp} {{\rm deg}_\partial ~}

\newcommand{\ldm} {{\rm L-dim ~}}
\newcommand{\kdm} {{\rm K-dim ~}}
\newcommand{\dm} {{\rm dim ~}}
\newcommand{\cdm} {{\rm codim ~}}

\newcommand{\ol} {\overline}
\newcommand{\ok} {{\overline k}}
\newcommand{\ox} {{\overline x}}
\newcommand{\oy} {{\overline y}}
\newcommand{\om} {{\omega}}

\newcommand{\ff}{{\bf F}}
\newcommand{\T}{{\bf T}}
\newcommand{\sph}{{\bf S}}

\newcommand{\be}{\begin{equation}}
\newcommand{\ee}{\end{equation}}

\title {Polynomials with general $\C^2$--fibers are variables. I.}

\vspace{.5cm}

\author{Sh. Kaliman}
\thanks{Partially supported by NSA grant}

\date{}
\maketitle

ABSTRACT. 
Suppose that $X'$ is a smooth affine algebraic variety of dimension 3
with $H_3(X')=0$ which is a UFD and whose invertible functions
are constants. Suppose that $Z$ is a Zariski open subset of $X$
which has a morphism $p :Z \to U$ into a curve $U$ such that
all fibers of $p$ are isomorphic to $\C^2$. We prove that
$X'$ is isomorphic to $\C^3$ iff none of irreducible components
of $X' \setminus Z$ has non-isolated singularities. Furthermore,
if $X'$ is $\C^3$ then $p$ extends to a
polynomial on $\C^3$ which is linear in a suitable coordinate
system. As a consequence we obtain the fact formulated in the
title of the paper. \\[2ex]

%\tableofcontents

\bigskip

\section{Introduction}

We say that a nonconstant polynomial on $\C^n$ is a variable if it is
linear in a suitable polynomial coordinate system on $\C^n$.
Classification of such polynomials is a difficult and
important problem which is
solved only for $n=2$.
In 1961 Gutwirth [Gu] proved  the following fact
which was later reproved by Nagata [Na]:
every polynomial $p$ in two complex variables
whose general fibers are isomorphic to $\C$ 
(that is, there exists a finite subset $S$ of $\C$ such that
for every $c \in \C \setminus S$ the fiber $p^{-1} (c)$
is isomorphic to $\C$) is a variable.
In 1974--1975 Abhyankar, Moh, and Suzuki showed that a much
stronger fact holds: every irreducible polynomial 
in two complex variables, whose zero fiber is isomorphic
to $\C$, is variable [AbMo], [Su].
The Embedding conjecture formulated by Abhyankar and
Sathaye [Sa1] suggests that the similar fact holds in
higher dimensions: 

\bigskip

\noindent {\it Every irreducible polynomial $p$ in $n$ complex variables
with a fiber isomorphic to $\C^{n-1}$ is a variable.}

\bigskip

It seems that in the full generality the positive answer to the
Embedding
conjecture is not feasible in the near future but there is some progress
for $n=3$. In this dimension A. Sathaye, D. Wright, and P. Russell proved
some special cases of this conjecture
([Sa1], [Wr], [RuSa], see also [KaZa1]).
Then M. Koras and P. Russell
proved the Linearization conjecture for $n=3$
[KoRu2], [KaKoM-LRu] which
implies the following theorem: 

if $p$ is an irreducible polynomial on $\C^3$
such that it is quasi-invariant with respect to a regular
$\C^*$-action on $\C^3$ and its zero fiber is isomorphic to $\C^2$,
then $p$ is a variable.
\footnote {In fact, P. Russell indicated to the author that
the ``hard-case" of the Linearization conjecture 
is equivalent to this theorem. 
This equivalence can be extracted from [KoRu1].} 

This paper and the next joint paper of the 
author with M. Zaidenberg [KaZa2]
contain another step in the direction of the Embedding conjecture --
we prove the analogue of the Gutwirth theorem in dimension 3,
i.e. every polynomial with general $\C^2$--fibers is a variable.
It is worth mentioning that a special case of this theorem
follows from more general results of Miyanishi [Miy1] and Sathaye [Sa2]
(we are grateful to P. Russell who drew our attention to the
paper of Sathaye). They showed that if each fiber
$p^{-1}(c), \, c\in \C$ of a polynomial $p \in \C [x,y,z]$
is isomorphic to $\C^2$ and its generic fiber is also plane
\footnote{ If $K$ is the field of fractions of $\C [p]$
then this means that the ring $\C [x,y,z] \otimes_{\C [p]} K$
is isomorphic to the polynomial ring in two variables over K.}
then $p$ is a variable.

In fact, in our papers the analogue of the Gutwirth theorem in
dimension 3 is also a consequence of the following more general result.

\bigskip

{\bf Main Theorem.} {\it
Let $X'$ be an affine algebraic variety of dimension 3
such that $X'$ is a UFD 
\footnote{An affine algebraic variety
is called a UFD if its algebra of regular functions is a UFD},
all invertible functions
on $X'$ are constants, and

(1) the Euler characteristic of $X'$ is
$e(X')=1$; 

(2) there exists a Zariski open subset $Z$ of $X'$
and a morphism $p: Z \to U$ into a curve $U$
whose fibers
are isomorphic to $\C^2$;

(3) each irreducible component of $X' \setminus Z$ is a UFD.

Then $U$ is isomorphic to a Zariski open subset
of $\C$ and $p$ can be extended to a regular function on $X'$. Furthermore, 
$X'$ is isomorphic to $\C^3$ and $p$ is a variable.

The same conclusion remains true if we replace (1) and (3) by

($1'$) $X'$ is smooth and $H_3(X')=0$;
\footnote{All homology groups which we consider in this paper
have $\Z$-coefficients.}

($3'$) each irreducible component of $X' \setminus Z$
has at most isolated singularities. 

\medskip

In the case when conditions ($1'$) and (2) hold but (3) does not,
$X'$ is an exotic algebraic structure on $\C^3$ (that is,
$X'$ is diffeomorphic to $\R^6$ as a real manifold
but not isomorphic to $\C^3$) with a nontrivial
Makar-Limanov invariant.}

\medskip

The Makar-Limanov invariant was introduced in [M-L1], [KaM-L1]
(see also [KaM-L2], [Za], and [De]). For a reduced irreducible affine
algebraic variety $X'$ (and $X'$ from the Main Theorem
is reduced and irreducible since it is a UFD) 
this invariant is the subalgebra $\ML (X')$
\footnote{In papers of Makar-Limanov it is
denoted by $\AK (X')$.}
of the algebra of regular functions $\C [X']$ on $X'$
that consists of all functions which are invariant under any
regular $\C_+$-action on $X'$. 
If $\ML (X')$ coincides with the ring of constants
then we say that it is trivial. This is so, for instance, when $X'$
is isomorphic to $\C^n$. 

\medskip

The proof of the Main Theorem can be divided in three major steps.
The first step is the following strengthened version
of the theorem of Miyanishi [Miy1] (which is essentially
based on [Sa2])

\bigskip

{\bf Lemma I.} {\it
Let $X'$ be an affine algebraic variety of dimension 3
such that $X'$ is a UFD, all invertible functions
on $X'$ are constants, and

(1) the Euler characteristics of $X'$ is $e(X')=1$;

($2'$) there exists a Zariski open subset $Z$ of $X'$
which is a $\C^2$-cylinder over a curve $U$ (i.e. $Z$ is isomorphic
to the $\C^2 \times U$);

($3$) each irreducible component of $X' \setminus Z$ is a UFD.

Then $X'$ is isomorphic to $\C^3$.

(4) Furthermore, the curve $U$ is a Zariski open subset of $\C$,
the natural projection from $Z$ to $U$ can be extended to
a regular function on $X'$, and this function is a variable.
\footnote{In the original formulation of Miyanishi 
the statement (4) was absent
but it was, obviously, known to Miyanishi.
The last statement of this Lemma was also absent in the
paper of Miyanishi.}

The statement of this Lemma remains true if conditions (1) and (3)
are replaced by

($1'$) $X'$ is smooth and $H_3(X')=0$;

($3'$) each irreducible component of $X' \setminus Z$
has at most isolated singularities.}

\bigskip

We give a new proof of this theorem based
on the notion of affine modifications
which was studied in [KaZa1].
In brief the idea of the
proof is as follows.

An affine modification
is just a birational morphism of reduced irreducible
complex affine algebraic varieties $\sigma : X' \to X$.
The restriction of such a morphism to the complement
of the exceptional divisor $E$ of $X'$ is an isomorphism
between $X' \setminus E$ and $X\setminus D$ where
$D$ is a divisor of $X$. The image $C_0=\sigma (E)$ is
called the geometrical center of modification.
For some affine modifications (which are called below
affine cylindrical modifications) $E$ is isomorphic to the
direct product $\C^k  \times C_0$ 
(where $k= \cdm_X C_0 -1$) which enables us to compare
the topology of $X'$ and $X$. Cylindrical modifications may not
preserve some features of $X$ (for instance, normality)
and their geometrical centers are not always closed.
Therefore, we introduce a subset of so-called basic
modification (in the set of all cylindrical modification)
which enable us to control geometrical changes more
effectively. We show
that under the assumption of the Miyanishi theorem
$X'$ is an affine modification of $X = \C^3$ and the divisor
$D$ is the union of a finite number of parallel affine
planes in $\C^3$. Then the problem will be reduced to the case
when $D$ consists of one plane only. 
One of the central facts for the first step is Theorem 2.3
which says that $\sigma$ is the composition
$\sigma_1 \circ \cdots \circ \sigma_m$ where each
$\sigma_i : X_i \to X_{i-1}$ ($X'=X_m$ and $X=X_0$) is a
basic modification. If $m=1$ and $C$ is either a point
or a straight line in the plane $D$ then it is easy to check that $X'$ is
isomorphic to $\C^3$ and the other statements of the Miyanishi
theorem hold. In the general case of $m>1$, using the control
over topology, one can show that the center of $\sigma_1$
is either a point or a curve in $D$ which is isomorphic to $\C$.
If the center is such a curve then it can be viewed as a straight
line by the Abhyankar-Moh-Suzuki theorem whence $X_1$ is
isomorphic to $\C^3$. Now the induction by $m$ implies the
Miyanishi theorem.

\medskip

In the second step we prove

\bigskip

{\bf Lemma II.}
{\it Let $X'$ be an affine algebraic variety of dimension 3
such that $X'$ is a UFD, all invertible functions
on $X'$ are constants, and 
let the assumption ($1'$) and ($2'$) of Lemma I hold,
but (3) does not.
Then $X'$ is an exotic algebraic structure on $\C^3$ 
with a nontrivial Makar-Limanov invariant.}

\bigskip 

Our proof of the Miyanishi theorem is longer than the
original one but it
%(though we tried to make it shorter at the expense
%of generalities) but it
has some advantage besides a slightly stronger formulation. It helps us to
cope with the second step. 
Namely, under the assumption of Lemma II $X'$ is still an
affine modification of $X=\C^3$, $\sigma$ is still
a composition of basic modifications, and
one can reduce the problem to the case when $D$ is a coordinate plane.
It can be shown that the geometrical center of modification
is either a point or an irreducible contractible curve in $D$.
Besides the Abhyankar-Moh-Suzuki theorem we have another
remarkable fact in dimension 2 -- the Lin-Zaidenberg theorem [LiZa]
says that such a curve is given by $x^n=y^m$ in a suitable
coordinate system where $n$ and $m$ are relatively prime.
This allows us to present explicitly a system
of polynomial equations in some Euclidean space $\C^N$ whose
zero set is $X'$. Here we use the fact that basic
modifications of Cohen-Macaulay varieties 
are Davis modifications which were introduced
in [KaZa1] and which fit perfectly the aim of presenting
explicitly the result of a modification as a closed
affine subvariety of a Euclidean space.
This explicit presentation of $X'$ as a subvariety
of $\C^N$ enables us to
compute the Makar-Limanov invariant of $X'$, using
the technique from [KaM-L1], [KaM-L2].
If condition ($3b'$) holds and $X'$ is smooth
then the Makar-Limanov invariant 
of $X'$ is non-trivial whence $X'$ is not isomorphic to $\C^3$.
On the other hand we show that $X'$ is contractible
and, therefore, it is diffeomorphic
to $\R^{6}$ by the Dimca-Ramanujam theorem [Di], [ChDi]     
which concludes the second step.

Third step is

\bigskip

{\bf Lemma III.}
{\it Let $q: Z \to U$ be a morphism of an affine algebraic
variety $Z$ into a curve $U$
such that every fiber of $q$ is isomorphic to $\C^2$.
Then there exists a Zariski open subset $U^*$ of $U$ such that
for $Z^*=q^{-1}(U^*)$ and $r=q|_{Z^*}$ the morphism $r: Z^* \to U^*$ is a
$\C^2$-cylinder over $U^*$ (that is, there exists
an isomorphism $\varphi : Z^* \to \C^2 \times U^*$ for which the 
composition of the projection to
the second factor and $\varphi$ coincides with $r$).}

\bigskip

The proof of Lemma III will be the content of the next
joint paper of the author and M. Zaidenberg [KaZa2].

The combination of Lemmas I, II, and III implies  Main Theorem.

\bigskip

It is our pleasure to thank M. Zaidenberg for his suggestion to
check Lemma II and many fruitful 
discussions. Actually, the idea of this paper arose during
the joint work of the author and M. Zaidenberg on [KaZa1].
Later M. Zaidenberg decided not to participate in the project
due to other obligations and the author had to finish it alone.

It is also our pleasure to thank I. Dolgachev
whose consultations were very useful for the author.

\bigskip

\section{Affine Modifications}

\bigskip

\subsection{Notation and Terminology}

In this subsection we present central definitions and notation 
which will be used in the rest of the paper.
The ground field in this paper will always be the field of complex numbers $\C$.
But it should be noted that all facts of this section hold
for every field of characteristics zero with the exception
of the results where the homology or fundamental groups are
mentioned.
%\footnote{Using this fact and the ``Lefschetz principle" (e.g., see [BCW])
%one can extend some other results of this paper to the case of
%a field $k$ of characteristics 0. For instance,
%every polynomial from $k^{[3]}$ with general $k^2$-fibers
%is a variable.}

\medskip

\defin We remind that an affine domain $A$ over $\C$ is
just the algebra of regular functions $\C [X]$ on a reduced
irreducible complex affine algebraic variety $X$.
Let $I$ be an ideal in
$A$ and $f \in A \setminus \{ 0 \}$.
By the affine modification of $A$ with locus $(I,f)$ we mean
the algebra $A': = A[I/f]$ together with the natural 
embedding $A \hookrightarrow A'$. 
That is, if
$b_0, \dots , b_s$ are generators of $I$ then $A'$ is
the subalgebra of the field ${\rm Frac} (A)$ of fractions of $A$
which contains $A$ and which is generated over $A$ by
the elements $b_0/f, \dots ,b_s/f$. It can be easily checked 
[KaZa1] that $A'$ is
an affine domain provided $A$ is.
Hence the spectrum of $A'$ is an
affine algebraic variety $X'$
and the natural embedding $A \hookrightarrow A'$
generates a morphism
$\sigma : X' \to X$.
Sometimes we refer to $\sigma$ as an affine modification or we say
that $X'$ is an affine modification of $X$.
The reduction $D$ (resp. $E$) of the divisor $f^{-1}(0) \subset X$ (resp.
$(f\circ \sigma)^{-1}(0) \subset X'$) will be called
the divisor (resp. the exceptional divisor) of the modification.
The subvariety of $X$ defined by the ideal $I$ 
(and sometimes the ideal $I$ itself) will be 
called the center of the modification, its reduction (which coincides,
of course, with the zero set of $I$ in $X$) will be called the reduced
center of the modification, and $\sigma (E)$
will be called the geometrical center of modification.
\footnote{Affine modifications appeared under different names
in several papers 
to which the author
and M. Zaidenberg did not pay sufficient attention
while we were writing [KaZa1].
In particular, we missed at that time the notion
of N\'{e}ron's blowing-ups [N\'{e}] which appeared
as early as in 1964 (see also [Ar], [WaWe]).}

\medskip

\rem (1) If $A'$ is as above and $f \notin I$ then consider the ideal
$J= \{ I,f \}$ in $A$ generated by $I$ and $f$. Clearly,
$A' =A[J/f]$ whence we shall suppose that $f \in I$ in the rest
of the paper. 

(2) It is easy to produce examples (and some of them appear below)
which show that the center, the reduced center, and the 
closure of the geometrical
center of a modification may be different but in all cases the
geometrical center is contained in the reduced center. Indeed, otherwise
one can choose an element $a \in A$ which vanishes 
identically on the latter but not
on the former. By Nullstellensatz $a^n \in I$ for some natural $n$.
On the other hand $a^n/f$ takes on the $\infty$-value on $X'$ whence
this function cannot be regular on $X'$. Contradiction.  

\medskip

\defin Let $p : Y \to Z$ be a morphism of 
algebraic varieties
%and let $Y$ be another algebraic variety. 
We say that $p$ is
a $\C^s$-cylinder over $Z$ if there exists an isomorphism
$\varphi : Y \to \C^s \times Z$ such that $p \circ \varphi^{-1}$
is the projection to the second factor.

\medskip

\defin 
Let $\sigma (E)$ be 
the geometrical center of the affine modification $\sigma$
from Definition 2.1.
Suppose that
$\sigma (E)$ is not just a constructive set but an algebraic variety
of pure dimension.
We say that $\sigma$ is a cylindrical modification of rank $s$
if $\sigma |_E : E \to \sigma (E)$ is a $\C^s$-cylinder 
where $s+1$ is the codimension of the geometrical center in $X$.

It is useful to know when
the geometrical center of a cylindrical modification coincides
with the reduced center and, in particular, is closed
(then we can better control the change of
topology under the modification). Here is the definition of
semi-basic modifications which are cylindrical and
whose geometrical centers coincide with the reduced
ones (see Proposition 2.7 below).

\medskip

\defin
Let $b_0, \dots , b_s $ be a sequence in an affine domain
$A= \C [X]$ which generates an ideal $I$.
We say that this sequence is semi-regular if
the height of $I$ is $s+1$
(or, equivalently, the zero set of $I$ in $X$ is
a subvariety of pure codimension $s+1$).
If in addition $b_0=f$ then the affine modification
$A\hookrightarrow A'$ with locus
$(I,f)$ will be called semi-basic of rank $s$,
and $b_0, \dots , b_s $ will be called a representative
system of generators for this modification.

\medskip  

We shall see that semi-basic modifications preserve
Cohen-Macaulay rings but they do not preserve normality.
For this purpose we need to consider a more narrow class
of modifications.

\medskip

\defin
Let the notation of Definition 2.4 hold. A semi-regular
sequence $b_0, b_1, \ldots , b_s $ is called
an almost complete intersection
provided the following two
conditions hold

(i) none of the irreducible components of the 
set $C$ of the common zeros
of $I$ in $X$ is contained in the
singularities of $X$;

(ii) for every irreducible component $C_i$ of $C$
there exists its Zariski open subset
$C_i^0 \subset {\rm reg} \, X$ which is a complete intersection given by
$b_0= \dots =b_s=0$ (in a neighborhood of $C_i^0$).
That is,  the gradients
of $b_0, \dots ,b_s$ are linearly independent at
generic points of each irreducible component of $C$.

If in addition $b_0=f$ then the affine modification
$A\hookrightarrow A'$ will be called basic of rank $s$,
and again $b_0, \dots , b_s $ will be called a representative
system of generators.

\medskip

We shall need also affine modifications which are
cylindrical (resp. semi-basic,
basic) only locally. These notions will appear in the next
subsection.

\medskip

\conv 
Further in this paper $X$ and $X'$ 
will {\bf always} be reduced irreducible affine algebraic
varieties. The algebra of regular functions of an affine
algebraic variety $Y$ will be $\C [Y]$. We put
$A = \C [X]$ and $A' = \C [X']$, that is, $A$ and $A'$
will {\bf always} be the affine domains that correspond to the affine
varieties $X$ and $X'$ respectively. Furthermore,
we suppose that the notation
$A\hookrightarrow A'$ is fixed throughout the paper.
It will {\bf always} mean an affine modification with locus
$(I,f)$. The corresponding morphism of the algebraic
varieties will {\bf always} be denoted by
$\sigma : X' \to X$. The divisor, the exceptional divisor,
and the reduced center of the modification will be
{\bf always} denoted by $D,E$, and $C$ respectively.

\medskip   

We shall also use the following notation in
the rest of this section: if $Y$ is
an affine algebraic variety and $B=\C [Y]$ then
for every closed algebraic subvariety $Z$ of $Y$ the
defining ideal of $Z$ in $B$ will be  denoted by $\cI_B (Z)$.
For every ideal $J$ in $\C [Y]$ we denote by $\cV_Y (J)$
the zero set of this ideal in $Y$. 

\bigskip

\subsection{General Facts about Affine Modifications}

\bigskip

We shall list first several useful facts from [KaZa1].

\medskip 

\thm {\it  {\rm (1) [KaZa1, Lemmas 1.1 and 1.2]}
Let $ A \hookrightarrow A'$ be an affine modification.
Then the fields of fractions
${\rm Frac} (A)$ and ${\rm Frac} (A')$ coincide, i.e. $\sigma$
is a birational morphism. The restriction of $\sigma$ to 
$X' \setminus E$ is an isomorphism between $X' \setminus E$
and $X \setminus D$.

{\rm (2) [KaZa1, Th. 1.1]}
Every birational morphism $X' \to X$ of affine
algebraic varieties is a modification. That is, there exist
an ideal $I \subset A$ and $f \in A$ such that $A' =A[I/f]$.

{\rm (3) [KaZa1, Prop. 1.2]}
Let $f=f_1f_2$ and $A'=A[I/f]$. Then $A'=A^1[I^2/f_2]$
where $A^1=A[I_1/f_1]$, $I_1$ is the ideal in $A$ generated
by $I$ and $f_1$, and $I^2$ is the ideal in $A^1$ generated
by $I/f_1$.

{\rm (4) [KaZa1, Prop. 3.1 and Th. 3.1]}
Let $ A \hookrightarrow A'$ be an affine modification.
Suppose that $E$ and $D$ are topological manifolds 
and they have the same number
of irreducible components. Furthermore, for every such component
$E_0$ of $E$ there exists a unique component $D_0$ of $D$
for which $E_0= \sigma^* (D_0)$ and
$\sigma (E_0) \cap {\rm reg} \, D_0 \ne \emptyset$.
Suppose also that $\sigma |_E : E \to D$ generates an isomorphism of
the homology groups. Then $\sigma : X' \to X$
generates isomorphisms of the homology and fundamental groups.

{\rm (5) [KaZa1, Cor. 2.1]} 
Let $X_1 = $spec$\,A_1$ be an
irreducible
closed subvariety of $X.$ Let the ideal
$I_1 \subset A_1$ consist of the restrictions to $X_1$ of the elements of $I.$
Suppose that $f_1 := f\,|\,X_1 \neq 0.$  Consider the
modification $A_1 \hookrightarrow A_1'$ with locus $(I_1,f_1)$
and the corresponding morphism of algebraic varieties
$\sigma_1 : X_1' \to X_1$ where  $X_1' = $spec$\,A_1'$.
Then there is a unique closed embedding $i'\,:\,X_1' \hookrightarrow X'$ 
making the following diagram commutative:

\medskip

\begin{picture}(150,75)
\unitlength0.2em
\put(45,25){$X_1'$}
\put(60,25){$\hookrightarrow$}
\put(62,29){{$ i'$}}
\put(75,25){$X^\prime$}
%\put(88,25){$\hookrightarrow$}
%\put(101,25){$\C^r \times \C^s$}
\put(45,5){$X_1$}
\put(60,5){$\hookrightarrow$}
\put(62,9){{$ i$}}
\put(75,5){$X$}
%\put(82,5){.}
%\put(92,5){$\hookrightarrow$}
%\put(108,5){$\C^r\,\,\,$}
\put(47,22){$\vector(0,-1){11}$}
\put(50,16){{$\rm \sigma_1$}}
\put(77,22){$\vector(0,-1){11}$}
\put(80,16){{$\rm \sigma$}}
%\put(110,22){$\vector(0,-1){11}$}
%\put(113,16){{$\rm pr_1$}}
\end{picture}

\medskip

\noindent where $i : X_1 \hookrightarrow X$ is the natural embedding.
In particular, affine modifications commute 
with direct products.}

\medskip

We need to discuss the behavior of affine modifications under
localizations (this should have been done in [KaZa1]).
Let $S$ be a multiplicative subset of $A$
and $S^{-1}A$
(resp. $S^{-1}A'$) be the localization of $A$
(resp $A'$) with respect to $S$. Every ideal
$I$ in $A$ generates an ideal $S^{-1}I$
in $\loc A$. The intersection of $\loc I$ with
$A$ is an ideal $S(I)$ which contains $I$.
The following fact is an immediate consequence of
the definitions of affine modifications and localizations.

\bigskip

\prop {\it In the notation above we have
$\loc A' = (\loc A) [\loc I/f]= (\loc A) [\loc (S(I))/f]$.
That is, localizations and affine modifications commute.}

%\bigskip
%\rem Let $T$ be the union of the zero sets of the elements of $S$.
%Consider the case when $T$ is an algebraic variety.
%Recall that geometrically this localization
%means switching from $X$ to $X\setminus T$,
%from the ideal $I$ to the ideal $\loc I=\loc (S(I))$.
%[AtMa]).
%Therefore,
%Proposition 2.2 says that $X' \setminus \sigma^{-1} (T)$
%is an affine modification of $X \setminus T$ with center
%$\loc I$.
%It is interesting to discuss the zero set
%of $\loc I=\loc (S(I))$.
%Consider a minimal primary decomposition
%$I= \bigcap_{i=1}^m \cQ_i$ where each $\cQ_i$ is primary
%and its associate prime is $\cP_i$.
%(This implies that the reduced center $C$ of the modification
%is the union of $m$ irreducible
%subvarieties $C_i = \cV_X (\cP_i )$.)
%Suppose that $S$ meets $\cP_i$ iff $i > k$.
%Then $\loc I = \bigcap_{i=1}^k \loc \cQ_i$,
%each $\loc \cQ_i$ is primary in $\loc A$ for $i\leq k$,
%and $S(I)= \bigcap_{i=1}^k \cQ_i$ [Ei, Th. 3.10 (d)].
%[AtMa, Prop. 4.9].
%A component $C_i$ is contained in $T$ iff $i>k$,
%and $\cV_{X\setminus T} (\loc I)=
%\bigcup_{i=1}^k C_i'$ (see the end of subsection 2.1 for notation) where
%$C_i'=C_i\setminus T$.

\medskip

\defin Suppose that $B$ is a localization of an affine domain,
$J$ is an ideal in
$B$ and $g \in B \setminus \{ 0 \}$.
By the local affine modification of $B$ with locus $(J,g)$ we mean
the algebra $B': = B[J/g]$ together with the natural
embedding $B \hookrightarrow B'$.
By Proposition 2.1 $B'$ is also a localization of
an affine domain. 
Hence the spectrum of $B$ (resp. $B'$) is a
(germ of an) affine algebraic variety $Y$ (resp. $Y'$)
and the natural embedding $B \hookrightarrow B'$
generates a morphism
$\delta : Y' \to Y$. We can define now the divisor,
the exceptional divisor, the (reduced, geometrical) center
of this local modification exactly in the same manner we did
for affine modifications.

\medskip

\rem By Proposition 2.1 each local affine modification
$B \hookrightarrow B'$ is just a localization of
an affine modification $A \hookrightarrow A'$
which respect to a multiplicative system $S \subset A$.
This implies that Theorem 2.1 (1)-(3) and (5) hold for
local affine modifications as well. Similarly,
some facts below (including Theorem 2.3) can be easily
reformulated for the local case.

\medskip

\defin  (1) A local affine modification
$B \hookrightarrow B'$ is called cylindrical (resp.
semi-basic, basic) if an affine modification
$A \hookrightarrow A'$ from Remark 2.2 can be chosen
cylindrical (resp. semi-basic, basic).

(2) Let $A \hookrightarrow A'$ be an affine modification.
Suppose that $M$
is a maximal ideal of $A$. Recall that the localization $A_M$ of $A$
near $M$ is the localization of $A$ with respect to the
multiplicative system $S = A\setminus M$. Let $I_M$ denote the
ideal generated in $A_M$ by $I$. We say that
this affine modification is locally
cylindrical (resp. semi-basic, basic) if for every point
of the geometrical center $\sigma (E)$ and the maximal ideal $M$, that vanishes
at this point, the local affine modification
$A_M \hookrightarrow A_M[I_M/f]= S^{-1}A'$ is
cylindrical (resp. semi-basic, basic).

\medskip

As we mentioned before there exist affine modifications
whose reduced center is different from the closure
of the geometrical one.
That is, for such a modification 
the natural projection
$E \to C$ is not dominant where $E$ is the exceptional divisor of
the modification and $C$ is its reduced center.
(Consider, for instance, $A' =A[\{ f \}/f^2]$. Here $E$ is empty
and $C$ is not.)
In order to have control over the change of topology
under an affine modification we need to understand when
the reduced center
coincides with the closure of the geometrical one. 
This requires the notion of the largest ideal of
a modification $A'=A[I/f]$ which was introduced in [KaZa1].

\medskip

\defin
The ideal $K= \{ a\in A | \, a/f \in A' \}$ in $A$ is called the
$f$-largest ideal of the modification $A\hookrightarrow A'$.
Clearly, $I \subset K$ and $A'=A[K/f]$.
%depends only on
%$A, A'$, and $f$.
When $A$ and $A'$ are fixed we
denote this largest ideal $K$ by $I_f$.

\medskip

\prop
{\sl Let $A\hookrightarrow A'$ be an affine modification such that $I=I_f$.
Then the reduced center of the modification coincides
with the closure of the geometrical one.
(This means that for every component $C_0$ of the reduced center,
$\sigma^{-1}(C_0)$ is a hypersurface in $X'$
whose image is dense in $C_0$.)}

\proof
Assume that $\sigma^{-1}(U)$ is empty for
some Zariski open subset $U$ of $C_0$. Choose a
regular function $a\in A$ so that $a$
vanishes on each component of the reduced center except for $C_0$
but we require that $a$ vanishes also on $C_0 \setminus U$.
Then $a'=a\circ \sigma$ vanishes on the
exceptional divisor $E$.
Note that the zeros of $f':= f\circ \sigma$ on $X'$
coincide with $E$.
By Nullstellensatz for some $n>0$ the element $(a')^n$ is
divisible by $f'$ in $A'$. Since
$(a')^n/f'=a^n/f \in A'$ we have $a^n\in I_f$
which shows that the reduced center does not contain $C_0$. Contradiction.
\qed

\medskip

\rem
As a consequence of Proposition 2.2 we see that for $I=I_f$ the number
of irreducible components of the exceptional divisor $E$ is at least
the same as the number of irreducible components of the reduced
center of the modification. There is a better estimate of
the number of irreducible components of $E$ (which will 
not be used further). Namely, it can be shown that
in the case of a normal affine domain $A$ this number
is greater than or equal to the number of ideals in a
minimal primary decomposition of $I$.

\medskip

Here is another property of the largest ideals of affine modifications
which will be used again and again in this section.

\medskip

\prop {\it 
Let $g \in A\setminus \{ 0 \}$ and
$f=g^n$ for a natural $n$.
Suppose that the defining ideal $\cI_{A'} (E)$ (see the end of
subsection 2.1 for notation) of the exceptional divisor
$E$ of the modification
$A \hookrightarrow A'$ coincides with the principal
ideal in $A'$ generated by $g$.
Then $(\cI_{A} (C))^n\subset I_f$ (in particular, for $n=1$
we have $\cI_{A} (C) =I_f$).

Furthermore, for every ideal $J$ in $A$ which is contained
in $\cI_{A} (C)$ the algebra $A_1: = A[J/g]$ is
contained in $A'$.}

\proof Note that for every $a \in (\cI_{A} (C))^n$ we have
$a \in (\cI_{A'} (E))^n$ whence $a/f \in A'$.
By Definition 2.8 $a \in I_f$ which is the first statement.
This implies that $g^{n-1}J \subset I_f$.
Hence $A_1 =A[g^{n-1}J/f] \subset A[I_f/f]=A'$.        
\qed

\medskip

\rem Suppose that $A_1$ be an affine domain such that
$A \subset A_1 \subset A'$ (for instance, $A_1$ is
from Proposition 2.3).  

(1) Since every generator of $A_1$ over $A$ is of form
$b/f^m$ where $b \in I^m$ we see that there exists
an ideal $L_1$ in $A$ such that $A_1 = A[L_1/f^m]$ for some
$m \geq 0$. 

(2) Furthermore there exists an ideal $K_1$ in
$A_1$ such that $A'=A_1[K_1/f]$ (it is enough to consider
the ideal generated by $I$ in $A_1$ as $K_1$).
Hence the modification $A_1 \hookrightarrow A'$ with
locus $(K_1,f)$ has the same
exceptional divisor $E$ as
the modification $A \hookrightarrow A'$ and the divisor of
$A_1 \hookrightarrow A'$ coincides with the exceptional
divisor of the affine modification $A \hookrightarrow A_1$
with locus $(L_1,f^m)$.  In fact,
under the assumption of Proposition 2.3 we can make a stronger
claim which will help us later to show that the number of
factors in the desired decomposition of an affine modification
into basic modifications is finite.

\medskip

\cor (cf. [WaWe, Prop. 1.2])
{\it Suppose that $J=\cI_A(C)$ in Proposition 2.3.
Then there exists an ideal $K_1$ in $A_1$ such that
$A'=A_1[K_1/g^{n-1}]$. That is, $A_1 \hookrightarrow A'$
may be viewed as an affine modification with locus $(K_1, g^{n-1})$.}

\proof Let $b_0=g^n,b_1, \dots , b_s$ be generators of $I$.
Note that $b_i/g \in A_1$ for every $i$. The ideal $K_1$ in $A_1$
generated by $g^{n-1},b_1/g, \dots , b_s/g$ is the desired ideal.
\qed

\medskip

We remind that we are planing to show that $X'$ from the Main
theorem is a modification of $\C^3$. Since both of these threefolds
are UFDs we need to have a closer look at affine modifications
$A \hookrightarrow A'$ of UFDs. We shall see that the assumption
of Proposition 2.3 on $\cI_{A'}(E)$ is automatically true in this case.

\medskip

\prop {\sl
Let $A \hookrightarrow A'$ be an affine modification, let $A$ be 
a UFD, and let the following conditions hold

(i) $f=g^n$ where $g \in A$ is irreducible;

(ii) the closure of the geometrical center of
the modification coincides with the reduced center
(by Proposition 2.2 this is so, for instance, when
$I$ is the $f$-largest ideal $I_{f}$ of the modification );

(iii) $A' \ne A[1/f]$ or, equivalently, the exceptional divisor
$E$ of the modification is not empty.

Then

(1) $g$ is irreducible as an element of $A'$;

(2) if $A'$ is also a UFD then $E$ and, 
therefore, $C$ are irreducible, and $\cI_{A'} (E)$ coincides with
the principal ideal generated by g.}

\proof
Let $g^k=a'b'$ where $a'=a/f^{l},b'=b/f^{m},a \in I^l$, and $b\in I^m$. Hence 
$g^{k+nl+nm}=ab$ in $A$. Since $A$ is a UFD we
have $a=ug^s$ and $b=vg^r$
where $s+r=k+nl+nm$ and $u, v$ are units. 
If $s<nl$ then $a'=u/f^{nl-s}$ whence $1/f \in A'$.
This contradicts (iii). 
Thus $s\geq nl $ and, similarly, $r \geq nm$. Hence
$a'=ug^{s-nl}$ and $b'=vg^{r-nm}$ are elements of $A$
which implies (1).

Assume that $E=E_1 \cup E_2$ where $E_1,E_2$ are effective divisors
of $X'$ (where $X'$ is as in Definition 2.1) 
without common irreducible components. If $A'$
is a UFD $E_k$ is the zero set of some $a_k' \in A'$.
Hence $g=u(a_1')^{n_1}(a_2')^{n_2}$ where
$u$ is a unit. Since $a_k'$ is not
a unit this contradicts to the fact that $g$ is irreducible
in $A'$ whence we have (2).
\qed

\rem (1) Remark 2.3 implies that if $A$ and $A'$ in Proposition 2.4
are UFDs then 
the ideal  $I_f$ is primary.   

(2) The following generalization of Proposition 2.4 will appear
in a coming paper [KaVeZa] of M. Zaidenberg, S. Venereau, and the author.
Let $A \hookrightarrow A'$ be an affine modification such that
$A$ and $A'$ are UFDs which have the same units. Then the numbers
of irreducible components in $D$ and $E$ are the same.

\medskip

\ex  (1) Lemma 2.3 implies that if $A$ is a UFD and $E$ is not
irreducible then $A'$ is not a UFD.  
Consider, for instance, $A=\C [x,y]$ and $f=x$. Let $I$
be generated by $x$ and $y^2-y$. Then $X'$ is the surface
in $\C^3$ (with coordinates $x,y,z$) given by
$xz=y^2-y$. This is a so-called Danielevski surface. These surfaces
are not UFDs. The exceptional divisor in this case consists of
two components $x=y=0$ and $x=y-1=0$.

(2) In the case when $E$ is irreducible $A'$ may be not a UFD
anyway, if, say,  $\cI_{A'} (E)$ is not
a principal ideal. Let $A$ and $f$ be as in the first example,
and let $I$ be generated by $x$ and $y^2$. Then $X'$ is the surface
in $\C^3$ given by $xz=y^2$. It is not a UFD.

\medskip

In general the units of $A$ and $A'$ differ 
(consider, for instance, $A'=A[1/f]$) but under
some mild assumption this is not the case
(in particular, in Proposition 2.4 $A$ and $A'$ 
have the same units).

\medskip

\prop {\sl
Let $A \hookrightarrow A'$ be an affine modification.
Suppose that for every 
natural $k$ each irreducible
divisor $g$ of $f^k$ in $A$  is not a unit in $A'$
or, equivalently, the set
$(g\circ \sigma)^{-1}(0)$ is not empty.
(When $A$ is a UFD it is enough, of course, to consider the
irreducible divisors of $f$ only.)
Then the units of $A'$ and $A$
are the same.}

\proof
Note that $A'$ is a subalgebra of $A[1/f]$. Thus its units
are also units of $A[1/f]$. The units of the last algebra
are the products of irreducible divisors of $f^k$ and
the units of $A$. By the assumption
these divisors are not invertible functions on $X'$
whence the units of $A'$ coincide with the units of $A$.
\qed

\medskip

\prop {\it Let $I_j$ be an ideal in $A$ for
$j=1, \dots , k$, and let $f_j \in I_j \setminus \{0 \}$.
Suppose that $f=f_1 \cdots f_k$ and 
$I =(f/f_1)I_1 + \ldots + (f/f_k)I_k$.
Let $A_j =A[I_j/f_j]$ and let
$\sigma_j : X_j \to X$ be the morphism of affine algebraic
varieties associated with the affine modification
$A \hookrightarrow A_j$ with locus $(I_j,f_j)$. 
Suppose that $E_j$ is the exceptional
divisor of this modification. These morphisms
define the affine variety $Y=X_1 \times_{X} X_2 \times_{X} \cdots \times_X X_k$
and its subvariety $Y^*=(X_1\setminus E_1)\times_{X} \cdots \times_X 
(X_k\setminus E_k)$. 

(1) The variety $X'$ is isomorphic
to the closure $\bY^*$ of $Y^*$ in $Y$ and under this isomorphism
$\sigma$ coincides with the restriction of the natural projection
$\tau : Y \to X$ to $\bY^*$.

(2) If $f_j$ and $f_l$ have no common zeros on $X$ for
every pair $j \ne l$ then $X'$ is isomorphic to $Y$.}

\proof
Let $D_j$ be the zero locus of $f_j$ on $X$. 
Then $D = \bigcup_{j=1}^k D_j$. Since 
the restriction of $\sigma_j$ to $X_j \setminus E_j$ is 
an isomorphism between $X_j \setminus E_j$ and
$X\setminus D_j$ we see that $Y^*$ is isomorphic to
$X\setminus D$. In particular, the natural projection
$\bY^* \to X$ generates an isomorphism of the fraction fields
of the algebras $B: =\C [\bY^* ]$ and $A$. 
That is, $B \subset {\rm Frac} \, (A)$.
The natural projection
$\bY^* \to X_j$ enables us to treat $A_j$ as a subalgebra of $B$.
Furthermore, since $\bY^*$ is the subvariety of $Y$
we see that $B$ is generated by these subalgebras
$A_1, \dots , A_k$. It remains to note that $A'$
is also generated by $A_1, \dots , A_k$ 
since $I =(f/f_1)I_1 + \ldots + (f/f_k)I_k$.
Hence $A'=B$ and $X'$ is naturally isomorphic to $\bY^*$
which yields (1).

For the second statement it suffices to prove that
$Y$ is irreducible. Assume that $Y$ has an irreducible component $Y_1$
different from $\bY^*$. Then the image of this component
under the projection $\tau : Y \to X$ must be contained in $D$
since the restriction of $\tau^{-1}$ to $X\setminus D$ is an isomorphism
between $X\setminus D$ and $Y^*$. We can suppose that
this image is contained in $D_1$. 
Put $T= \bigcup_{j=2}^k D_j$ and
consider $\theta : Y\setminus \tau^{-1} (T)
\to X\setminus T$ where $\theta$ is the restriction of $\tau$.
Since for $j \geq 2$ the restriction of $\sigma_j$ to $X_j \setminus 
\sigma_j^{-1} (T)$ is an isomorphism between this variety and
$X\setminus T$ we see that $Y\setminus \tau^{-1} (T)$ is
isomorphic to $X_1\setminus \sigma_1^{-1} (T)$ and $\theta$
coincides with the restriction of $\sigma_1$ to
$X_1\setminus \sigma_1^{-1} (T)$ under this isomorphism.
Thus $\sigma_1^{-1} (X \setminus T) \simeq 
\theta^{-1} (X\setminus T) \simeq \tau^{-1} (X\setminus T)$.
Note that $T$ does not meet $D_1$ by the assumption of this Proposition.
Hence $\tau^{-1} (X\setminus T)$ contains $Y_1$ and, therefore,
it is not irreducible. On the hand $\sigma_1^{-1} (X \setminus T)
\subset X_1$ is irreducible. This  contradiction concludes (2).
\qed

\medskip

\rem
Let us discuss the coordinate meaning of Proposition 2.6 (2).
Suppose for simplicity that $k=2$.
Let $X$ be a closed affine subvariety of $\C^n$ with
a coordinate system $\bx$ and let
$X_j$ be a closed affine subvariety of $\C^{n_j}$
with a coordinate system $(\bx , \bz_j)$.
That is, $\C^{n_j}$ contains the above sample of $\C^n$
as a coordinate $n$-plane.
Suppose that $X_j$ coincides with the zeros
of a polynomial system of equations $P_j(\bx , \bz_j )=0$
and $\sigma_j$ can be identified
with the restriction of the natural projection $\C^{n_j} \to \C^n$.
Consider the space $\C^{n_1+n_2-n}$
with coordinates $(\bx , \bz_1 , \bz_2)$.
Then Proposition 2.6 (2) implies that
the set of zeros of the system $P_1(\bx ,\bz_1)=P_2(\bx ,\bz_2)=0$
in this space is isomorphic to $X'$.

\bigskip

\subsection{Semi-basic Modifications}

In general the topologies of the exceptional divisor $E$
and the reduced center $C$ of an affine modification 
are not related very well even in the case when
the reduced center coincides with the closure of the geometrical one.
The natural projection $E \to C$ 
may not be surjective and, furthermore,
its generic fibers may be not connected.  

\medskip

\ex
Consider $A$ equal to the ring $\C [x,y,z]$ of polynomials
in three variables. Let $f=x$ and the ideal $I$ be
generated by $x$ and $y$.
Then $A'=A[I/f]=\C [x,u,v]$ where $y=xu,z=v$. In this
case $C$ is the $z$-axis and $E$ is the $uv$-plane.
Let $\Gamma$ be a closed curve in the $uv$-plane whose
projection to the $v$-axis is dominant but neither surjective
nor  injective. Suppose that $g(u,v)=0$ is the defining
equation of this curve in the $uv$-plane. Consider the ideal
$J$ in $A'$ generated by $f$ and $g$,
and put $A''=A'[J/x]$. The sequence $f,g$ is regular in $A'$, and
we shall see later in this subsection
(Proposition 2.7) that its exceptional
divisor $F$ is isomorphic to $\Gamma \times \C$.
Note that the natural embedding of $A$ into $A''$ is also
a modification by Theorem 2.1 (2) which has the same exceptional divisor
$F$.
The projection $F \to C$ is the composition of
the projections $F \to \Gamma $ and $\Gamma \to C$.
This yields the desired example.

\medskip
           
This example suggests also an approach to what should
be done in order to track the change of topology. We shall try
to present an affine modification
$A \hookrightarrow A'$ as a composition of
basic modifications.
If $A_1$ is an affine domain such that $A\subset A_1 \subset A'$
then this modification is the composition of
$A \hookrightarrow A_1$ and $A_1 \hookrightarrow A'$,
where the last two embeddings are affine modifications
by Theorem 2.1 (2) (see also Remark 2.4).
When $f=g^n$ then Proposition 2.3 suggests to
look for $A_1$ in the form $A_1 =A[J/g]$ where
$J \subset \cI_A (C)$. Our first aim in this subsection is
to show that $J$ can be chosen so that the affine modification
$A \hookrightarrow A_1$ with locus $(J,g)$ is semi-basic, and, under
some additional assumption, even basic.

\medskip

\lemma {\it
Let  $C$ be a closed reduced subvariety of $X$ of codimension $s+1$
and let $I=\cI_A (C)$. Suppose that $f \in I \setminus \{ 0 \}$. 

(1) Then one can choose a
semi-regular sequence $f=b_0, \ldots , b_s$
whose elements are contained in $I$.

(2) Let this sequence generate an ideal $J$.
If none of the irreducible components of $C$
and none of the irreducible components
of the zero divisor of $f$ is
contained in the singularities ${\rm sing} \,X$ of $X$
then the sequence above can be chosen so that
none of the irreducible components of $\cV_X (J)$ is
contained in ${\rm sing} \,X$.

(3) If the assumption of (2) holds and the zero
multiplicity of $f$ at generic points of each
irreducible component of $C$ is 1, then the sequence
$b_0=f, b_1, \dots , b_s$ can be chosen so that
it is an almost complete intersection.

(4) There exists a finite-dimensional subspace $S$ of $I$ such that
(1)-(3) are true when $b_1, \dots , b_s$ are generic points of any 
finite-dimensional subspace of $I$ which contains $S$.}

%(4) Let the assumption of (2) hold and let $C^i 
%\, i=1, \dots ,k$ be closed irreducible
%subvarieties of $C$. Suppose that
%none of these subvarieties is contained
%in ${\rm sing} \,X$ of $X$
%and the zero
%multiplicity of $f$ at generic points of each $C^i$ is 1.
%Then $b_0, \dots , b_s$ from (3) can be chosen so 
%that, in addition, the gradients
%of these regular functions on $X$ are linearly independent
%at generic points of $C^i$ for every $i=1, \dots ,k$.}
%(5) Let $b_1^0, \dots , b_s^0 \in I$ and let $k\in \N$.
%Then in (1) and (2) each $b_j$ $(j=1, \dots ,s)$ can be
%chosen in the form $b_j=b_j^0 +b_j^1$ where $b_j^1 \in I^k$.}
\proof 
Suppose that $X$ is a closed subvariety of $\C^n$ and
$\bx = (x_1, \dots , x_n)$ is a coordinate system on $\C^n$.
Let $g_0, g_1, \dots , g_r$ be generators of $I$ and let
$f$ be one of them (say, $f=g_0$).
Consider the set $S_m$ of $(m+1)$-tuples $\bb_m =(b_0=f, b_1, \dots , b_m)$
such that each $b_i \, (i \geq 1)$ is of 
the form $\sum_{j=0}^r l_j(\bx )g_j$
where each $l_j(\bx )$ is a linear polynomial on $\C^n$.
Let $W(\bb_m)$  be the set of common zeros of $b_0=f, b_1, \dots , b_m$.
We want to show that for $m\leq s$ the statements (1)-(3)
are true with $s$ replaced by $m$. By the assumption this is so
when $m=0$. Suppose that this is true for $m<s$ and show that
(1)-(3) holds for $m+1$. For every point $o \in X \setminus C$
we can find $g_j$ which does not vanish at $o$.
Thus, taking, a linear combination of $g_0, g_1, \dots , g_r$
as $b_{m+1}$ we can suppose that $b_{m+1}$ does not vanish
at generic points of every irreducible component of $W(\bb_m)$.
Hence the sequence $b_0, \dots , b_{m+1}$ generates an 
ideal $J_{m+1}$ such that $\cV_X (J_{m+1})$ has codimension $m+2$ in $X$.

If the assumption of (2) holds then
$T = {\rm sing } \, X \cap W(\bb_m )$ does not contain
any irreducible component of $W(\bb_m)$. We can assume
that $b_{m+1}$ does not vanish at generic points of 
every irreducible component of $T\setminus C$. Hence
${\rm sing } \, X \cap W(\bb_{m+1} )$ does not contain
any irreducible component of $W(\bb_{m+1})$ which is (2).

Since $S_{m+1}$ can be viewed as an algebraic variety and
$W(\bb_{m+1})$ depends algebraically on $\bb_{m+1}$ we see that
the number of irreducible components of $W(\bb_{m+1})$ at generic points
of $S_{m+1}$ is the same. 
Thus if we perturb $b_{m+1}$ by an element of form
$\sum_{j=0}^r l_j(\bx )g_j$ we shall still have (1) and (2),
i.e. (1) and (2) hold
for generic $\bb_{m+1}\in S_{m+1}$. Since the codimension of $C$ is $s+1$
and $I= \cI_A (C)$ for every generic point of $C$ (which is 
a smooth point and which belongs to ${\rm reg} \, X$)
we can always find $m+1$ elements among $g_1, \dots , g_r$
such that the gradients of these elements and
the gradient of $g_0=f$ are linearly independent at this point.
Thus one can suppose that the gradients of $b_0, \dots , b_{m+1}$
are linearly independent at generic points of $C$.
%(and, similarly, $C^i$ for every $i$ which yields the desired
%induction for (4)).
Consider a generic point $o \in W (\bb_{m+1}) \setminus C$.
One can suppose that $g_1(o) =1$. Choose a linear polynomial
$l_1 (\bx )$ such that $l_1 (o)=0$. Then
the gradients of $g_1l_1$ and $l_1$ at $o$ are the same. Therefore,
perturbing $b_{m+1}$ by a function of form $g_1l_1$
one can suppose that the gradients of
$b_0, \dots , b_{m+1}$ are linearly independent at $o$
whence they are linearly independent at generic points of
the irreducible component of $W(\bb_{m+1})$ that contains $o$.
Since the number of irreducible components
of $W(\bb_{m+1})$ does not change in a neighborhood 
of a generic point $\bb_{m+1}$ of $S_{m+1}$
we can perturb $b_{m+1}$ so that
we have the linear independence of the gradients at generic
points of each irreducible component of $W(\bb_{m+1})$
which is (3). 

For (4) it suffices to put $S=S_{s}$.
%For (5) we choose the sequence $g_0, \dots , g_r$ so that it
%contains also generators of $I^k$. Note that if
%$o \in X\setminus C$ then
%one of these generators does not vanish at $o$.
%Thus we can accomplish (1) and (2) using only
%perturbation of $b_1^0, \dots , b_s^0$ by linear
%combinations of these generators of $I^k$.
\qed

\medskip

%\rem Let the assumption of Lemma 2.1 (3) hold and let
%$C_0$ be a reduced subvariety of $C$ such that 
%none of the irreducible components of $C_0$
%is contained in the singularities of $X$.
%Suppose also the zero
%multiplicity of $f$ at generic points of each
%irreducible component of $C_0$ is 1. Note that Lemma 2.1 (4) implies
%that $b_1, \dots , b_s$ can be chosen so that the gradients
%of $b_0, \dots , b_s$ are linearly independent not only
%at generic points of $\cV_X (J)$ but also
%at generic points of each irreducible component of $C_0$.
%\medskip
\prop {\it Suppose that
$A \hookrightarrow A'$ is a semi-basic modification of rank $s>0$.
Then it is a cylindrical modification of rank $s$. Furthermore,
the reduced and geometrical centers of this modification
coincide.}

\proof We shall begin with an example which was also described
in [KaZa1]. Let $J_0$ be the maximal ideal in $\C^{[s+1]}=
\C [x_0, x_1, \dots ,x_s]$ generated by all coordinates
(i.e. it corresponds to the origin $o$ in $\C^{s+1}$).
Put $B_0=\C^{[s+1]} [J_0/x_0]$ and consider the
modification $\C^{[s+1]} \hookrightarrow B_0$ with locus
$(J_0,x_0)$. Then $B_0$ is isomorphic to
$\C [x_0, y_1, \dots ,y_s]$ and $x_i =x_0y_i$ for $i=1, \dots ,s$.
That is, $Z_0: = {\rm spec} \, B_0$ may be viewed as the
subvariety of $\C^{2s+1}$ (whose coordinates are
$x_0, x_1, \dots ,x_s, y_1, \dots ,y_s$) given by the system
of equations $x_i-x_0y_i=0, \, i \geq 1$. Let $\rho : \C^{2s+1} \to \C^{s+1}$
be the natural projection to the first $s+1$ coordinates.
Then our modification is nothing but the restriction of $\rho$
to $Z_0$. Its reduce and geometrical
centers are $o$ and the exceptional divisor
$E_0 =\rho^{-1} (o) \simeq \C^s$.

Put $Z= \C^{s+1} \times X$ and $B =\C [Z]$.
That is, $B=\C^{[s+1]} \otimes A = A^{[s+1]}=A[x_0, x_1, \dots ,x_s]$.
Put $J=J_0B$ and consider the modification
$B \hookrightarrow B'$ with locus $(J,x_0)$.
Since modifications commute with direct products
(see Theorem 2.1 (5)) we see that $Z' := {\rm spec} \, B' =
Z_0 \times X$ and the above modification is the restriction $\delta$
to $Z' \subset \C^{2s+1} \times X$ of the natural
projection $(\rho , {\rm id} ) : \C^{2s+1} \times X
\to \C^{s+1} \times X =Z$. In particular, its 
reduced and geometrical centers are $C^0 = o \times X$
and the exceptional divisor $E^0 = E_0 \times X$.

Let $b_0=f, b_1, \dots , b_s$ be a representative system of
generators for $A \hookrightarrow A'$. In particular,
this sequence generates $I$. Consider the embedding
$i : X \hookrightarrow Z$
given by the system of equations
$x_i-b_i=0, \, i=0, \dots ,s$. 
Then the restriction of
$J$ to $X$ coincides with $I$. By Theorem 2.1 (5) we have 
the commutative diagram

\medskip

\begin{picture}(150,75)
\unitlength0.2em
\put(45,25){$X'$}
\put(60,25){$\hookrightarrow$}
\put(62,29){{$ i'$}}
\put(75,25){$Z^\prime$}
%\put(88,25){$\hookrightarrow$}
%\put(101,25){$\C^r \times \C^s$}
\put(45,5){$X$}
\put(60,5){$\hookrightarrow$}
\put(62,9){{$i$}}
\put(75,5){$Z$}
%\put(82,5){.}
%\put(92,5){$\hookrightarrow$}
%\put(108,5){$\C^r\,\,\,$}
\put(47,22){$\vector(0,-1){11}$}
\put(50,16){{$\rm \sigma$}}
\put(77,22){$\vector(0,-1){11}$}
\put(80,16){{$\rm \delta$}}
%\put(110,22){$\vector(0,-1){11}$}
%\put(113,16){{$\rm pr_1$}}
\end{picture} 

\medskip

\noindent where $i' : X' \hookrightarrow Z'$ is a closed embedding.
The reduced center of $\sigma : X' \to X$ coincides 
with $C =C^0\cap i(X)$,
and it is of codimension $s+1$ in $X$ since $\sigma$
is semi-basic.
Since $E$ is of codimension 1 in $X'$ we see that
each generic fiber $F$ of $\sigma |_E : E \to \sigma (E) \subset C$
must be at least of dimension $s$. But $F$
is contained in a fiber $F^0\simeq \C^s$ of 
$\delta |_{E^0} : E^0 \to C^0$. Hence $\dim F =s$ and
$\sigma (E)$ is dense in $C$. Furthermore, since
$i'$ is a closed embedding $F=F^0$ and
$\sigma (E)=C$ whence $E$ is a $\C^s$-cylinder
over $C$ and the reduced and geometrical centers
of this modification coincide.
\qed

\medskip

\ex Not every cylindrical modification is semi-basic
(even locally). Consider the algebra $A$ of regular functions
on $X = \{ xy=zt \} \subset \C^4_{x,y,z,t}$. Let $f=x^2$
and the ideal $I$ be generated by $x^2, yx,y^2,$ and $z$.
The the reduced (and the geometrical) center $C$ of the modification
$A \hookrightarrow A'$ is the line $x=y=z=0$. Clearly
this modification is not semi-basic since this line
cannot be given by the zeros of $x^2$ and one more regular
function. One can present $X'$ as the hypersurface
$\{ u=vt \} \subset \C^4_{x,u,v,t}$ where $xu=y$ and $x^2v=z$
and the exceptional divisor $E$ in $X'$ coincides with the
zeros of $x$. It is easy to see that the projection
$E \to C$ is the cylinder.

\bigskip

\subsection{Davis modifications}

It is useful to compare semi-basic modifications with
Davis modifications which were introduced in [KaZa1]

\medskip

\thm ([Da], see also [Ei, Ex. 17.14])
{\it Let $f=b_0, b_1, \dots , b_s$ be generators of an 
ideal $J$ in a Noetherian domain $B$.
Consider the
surjective homomorphism 
$$\beta\,:\,B^{[s]}=B[y_1,\dots,y_s]
\longrightarrow\!\!\!\to B[J/f] = B[b_1/f,\dots,b_s/f] \simeq B'$$
where $y_1, \dots , y_s$ are independent variables and $\beta(y_i)=
b_i/f,\,\,\,i=1,\dots,s.$
Denote by $J'$ the ideal of the polynomial algebra $B^{[s]}$
generated by the elements $L_1,\dots,L_s \in {\rm ker}\,\beta$
where $L_i = fy_i - b_i.$
Then ${\rm ker}\,\beta$ coincides with $J'$ 
iff $J'$ is a prime ideal.
The latter is true, for instance, if
the system of generators $b_0 = f,\,b_1,\dots,b_s$ of the ideal $J$
is regular \footnote{I.e. the ideal $(b_0,\dots,b_{s})$ is
proper and for each $i=1,\dots,s$ the image of $b_i$
is not a zero divisor in $B/(b_0,\dots,b_{i-1}).$}.}

\medskip

\defin Let $B$ be (a localization of) an affine domain.
When $J'$ from Theorem 2.2 is prime the (local) affine modification
$B \hookrightarrow B'$ with locus $(J,f)$ is called Davis,
and the sequence $b_0, b_1, \dots , b_s$ is called a representative
system of generators for this modification.

\medskip

\rem
It is easy to see that
every (local) affine Davis modification $B \hookrightarrow B'$
as above is cylindrical of rank $s$ [KaZa1, Prop. 1.1 (c)].
The reduced and the geometrical centers
coincide (but the center of the modification can still be different
from the reduced center, see Examples 1.2 and 1.5 in [KaZa1]).
In particular, the codimension of every irreducible component of
the reduced center 
is $s+1$ unless this center is empty. 
Hence this Davis modification is not only cylindrical
but automatically semi-basic in the case of a non-empty
reduced center.

\medskip 

Recall that if $J$ is an ideal of an affine domain
$A$ and $M$ is a maximal ideal in $A$ then we
denote the localization
of $A$ near $M$ (i.e. the localization 
with respect to the multiplicative system $S =A\setminus M$)
by $A_M$ and the ideal
generated by $J$ in $A_M$ by $A_M$.

\medskip

\prop {\it 
Let $A \hookrightarrow A'$ be an affine modification, and
let $b_0=f, b_1, \dots , b_s$ be a system of
generators of $I$. If $M$ is a maximal ideal in $A$ then
by Proposition 2.1
$A_M \hookrightarrow S^{-1} A'$ is the local affine modification
with locus $(I_M, f)$. Suppose that for every maximal ideal $M$ this
local modification is
Davis and $b_0, \dots , b_s$ is a representative
system of generators in the ideal $I_M$.
Then $A \hookrightarrow A'$ is a Davis modification.}

\proof Let $A^{[s]}=A[y_1, \dots ,y_s]$ and let
$I'$ be the ideal in $A^{[s]}$ generated by
$L_i=y_if -b_i, \, i=1, \dots ,s$.
Put $Y= \C^s \times X$ (in particular, $A^{[s]}=\C [Y]$).
Let $Y_1$ be the subvariety of $Y$ defined by the ideal $I'$
in $A^{[s]}$.
We need to show that $I'$ is prime,
i.e. $Y_1$ is reduced irreducible.
Choose a maximal ideal $M'$ in $A^{[s]}$ which 
vanishes at
a point $x' \in X'$. Let $x= \sigma (x')$
where $\sigma : X' \to X$ is the natural projection and
let $M$ be the maximal ideal of $A$ that vanishes
at $x$. Then $A \setminus M \subset A^{[s]}\setminus M'$ and
$A^{[s]}_{M'}$ is a further localization of $S^{-1}A^{[s]}$.
Since $A_M \hookrightarrow S^{-1} A'$ is a Davis
modification and $b_0, \dots , b_s$ is a representative
system of generators  of this modification by assumption, the ideal
$S^{-1}I'$ is prime in $S^{-1}A^{[s]}$.
But the localization $I_{M'}'$ of this ideal must be also prime. 
Hence the germ of $Y_1$
at $x'$ is reduced irreducible. 

If we want to claim the same 
about $Y_1$ we need to show that it is connected (for instance,
irreducible). Note that $E_1=Y_1 \cap f^{-1}(0)
\simeq \C^s \times C$ where $C=\{ b_0 = \dots = b_s=0 \}$
is the reduced center of the modification. Since
the localizations of our modification are Davis the codimension
of irreducible component of $C$ in $X$ must be $s+1$ by Remark 2.7. Hence
$\dim E_1 = \dim X -1$ unless $E_1$ is empty. By construction
$Y_1 \setminus E_1$ is isomorphic to $X \setminus D$ and,
therefore, irreducible. Furthermore, since the codimension
of each irreducible component of $Y_1$ in $Y$ is at most $s$
(i.e. the dimension of such a component is at least $\dim X$),
we see that the numbers of irreducible components of
$Y_1$ and $Y_1 \setminus E_1$ are the same. 
Thus $Y_1$ is irreducible.
\qed

\medskip

Another notion we have to use is Cohen-Macaulay rings.
We send the readers for the definition and properties of
these objects to [Ei] or [Ma]. We reverse first Remark 2.7
and show that every semi-basic modification of
a Cohen-Macaulay affine domain is Davis.

\medskip

\prop {\it Let $A$ be Cohen-Macaulay and
$I$ the ideal generated by a semi-regular
sequence $b_0=f, b_1, \dots ,b_s$.
Then the affine modification $A \hookrightarrow A'$
is Davis.}

\proof Let $M$ be a maximal ideal in $A$. Then $A_M$ is
also Cohen-Macaulay [Ma, Th. 30]. In the local ring $A_M$
every semi-regular sequence is regular [Ma, Th. 31].
Thus the modification $A_M \hookrightarrow A_M[I_M/f]=
S^{-1}A'$ (where $S=A\setminus M$) is Davis by Theorem 2.2.
Hence $A \hookrightarrow A'$ is Davis by Proposition 2.8.
\qed    

\medskip

\prop {\it Suppose that
$A \hookrightarrow A'$ is a Davis modification.
Let $A$ be Cohen-Macaulay.
Then $A'$ is also Cohen-Macaulay.}

\proof 
Let $L_1, \dots , L_s \in A^{[s]}$ be as in the proof of Proposition 2.8.
Since $A$ is Cohen-Macaulay
$A^{[s]}$ is Cohen-Macaulay as well [Ei, Prop. 18.9].
The ideal $I'$ generated by $L_1, \dots , L_s$ has height $s$, i.e. its
zero set has codimension $s$ in ${\rm spec} \, A^{[s]}$.
Hence $A' \simeq A^{[s]}/I'$ is Cohen-Macaulay by [Ei, Prop. 18.13].
\qed

\bigskip

\subsection{Basic Modifications and Preservation of Normality
and UFDs}

We saw (Example 2.1) that semi-basic modifications do not
preserve UFDs. Furthermore, they do not preserve normality in general,
and we shall need normality.

\medskip 

\ex
Let $A=\C [x,y]$, $f=x^2$ and $I$ is generated by $f$ and
$y^2$. Consider $A \hookrightarrow A'$. Then $A'$ is not
normal in the worst possible scenario: $X'$ is given
in $\C^3$ by $x^2z=y^2$ and it has selfintersection points
in codimension 1.

\medskip

We shall show that in the case of Cohen-Macaulay varieties
normality survives under basic modifications, but let us
emphasize first some nice properties of these modifications.

\medskip

\rem Let $b_0=f, \dots , b_s$ be a representative system of
generators of a basic modification $A \hookrightarrow A'$.
Note $b_0, \dots , b_s$ may be viewed as elements
of a local holomorphic coordinate system at a generic
point $x$ of the reduced center of the modification.
This implies that
every point $y \in \sigma^{-1}(x)$ is a smooth point
of $X'$ and the zero multiplicity of $f\circ \sigma$
at $y$ is 1. Actually one can see that locally
this modification at $x$ is nothing but a usual
(affine) monoidal transformation.

Remark 2.8 and Theorem 2.1 (4) imply

\medskip

\prop {\it Let $A \hookrightarrow A'$ be
a basic modification.
Suppose that $C$ (and, therefore,
$E$) is irreducible and a topological manifold.
Suppose also that the natural embedding of $C$ into $D$ generates
an isomorphism of the homology of $C$ and $D$. Then $\sigma$ generates
isomorphisms of the fundamental groups and the homology
groups of $X$ and $X'$.}

\medskip

\prop {\it Let $A \hookrightarrow A$ be a basic modification.
Suppose that $A$ is normal and Cohen-Macaulay.
Then $A'$ is normal and Cohen-Macaulay.}

\proof By Proposition 2.9 this modification is Davis.
Thus by Proposition 2.10 $A'$ is Cohen-Macaulay.
Note that if the singularities of $X'$ is at least of codimension 2
then $X'$ is normal by [Ha, Ch. 2, Prop. 8.23].
Since $X$ is normal
the codimension of $\sigma^{-1} ({\rm sing} \, X \setminus D)
= \sigma^{-1} ({\rm sing} \, X \setminus C)$ in 
$X'$ is at least 2 whence we can ignore this subvariety.
Let $C^0$ be the subset of the reduced center $C$, at the points
of which the gradients of a representative system of generators
are linearly independent. By the definition of a basic
modification the codimension of $C\setminus C^0$ in $C$ is
at least 1. Since $\sigma$ is cylindrical
the codimension of $\sigma^{-1} (C\setminus C^0)$ in
$E$ is at least 1 and in 
$X'$ is at least 2, and we can ignore these points again.
The other points of $X'$ are smooth by Remark 2.8.
\qed

\medskip

As soon as we control normality we can take care of
preservation of UFDs under affine modifications.

\medskip

\prop {\it Let $A \hookrightarrow A'$ be an affine modification
of normal affine domains such that $E$ and $D$
are irreducible, and let $f=g^{n}$ where $g \in A$.

(1) Suppose that $A'$ is a UFD and 
the defining ideal of $D$
is generated by $g$.
Then $A$ is a UFD.

(2) Suppose that $A$ is a UFD and 
the defining ideal
of $E$ in $A'$ is generated by
the regular function $g':=g\circ \sigma$ on $X'$. 
Then $A'$ is a UFD.}

\proof (1) Let $S$ be a closed irreducible hypersurface in $X$
which is different from $D$
and let $S'$ be its strict transform (i.e. the closure
of $\sigma^{-1}(S\setminus D)$ in $X'$). Since $A'$ is a
UFD the defining ideal of $S'$ in $A'$ coincides with
the principal ideal generated by a regular function
$h'$ on $X'$. Note that $h'=h/g^{k}$ 
where $h\in A$ is not
divisible by $g$. 
Hence $S\setminus D$ coincides with the
zeros of $h$ on $X \setminus D$ and
the zeros of $h$ in $X$ does not contain $D$.
The zero multiplicity of $h$ at generic points of $S$ is the same
as the zero multiplicity of $h'$ at generic points of $S'$,
i.e. it is 1. If $e$ is another function which vanishes on $S$
then $e/h$ is regular at these generic points and on
$X\setminus D$ whence it is regular 
on $X$ except a subvariety of codimension 2.
Hence $e/h$ is holomorphic on $X'$ [Rem, Lemma 13.10]
and, therefore, regular (e.g., see [Ka2]).
Thus the defining ideal of $S$
is principal and $A$ is a UFD.

(2) Let $S'$ be a closed irreducible
algebraic hypersurface in $X'$. 
We disregard the case when
$S'$ coincides with $E$ since the defining ideal of $E$
is generated by $g'$.
Then $\sigma (S')$ is a constructive set and its closure $S$
is an irreducible hypersurface in $X$. By the assumption 
the defining ideal of $S$ in $A$ is the principal ideal generated by
a regular function $h$ on $X$.
Suppose that $h':=h \circ \sigma$ has zero
multiplicity $r$ at generic points of $E$. Then
$e= h'/(g')^{r} $ is regular 
at these generic points
(and on $X' \setminus E$, of course). By the same argument about
deleting singularities in codimension 2
we conclude that $e$ is regular on $X'$.
By construction its zeros
on $X' \setminus E$ coincide with $S' \setminus E$ and these
zeros do not contain $E$. 
Hence $S' = e^{-1} (0)$ and furthermore
the zero multiplicity of $e$ at generic points of $S'$ is 1
(since it is the same of the zero multiplicity of $h$ at generic
points of $S$). Using again the argument about deleting singularities
in codimension 2 we see that every regular function which vanishes
on $S'$ is divisible by $e$. 
Thus $X'$ is a UFD.
\qed

\medskip

\lemma {\it Let $A \hookrightarrow A'$ be an affine modification
and $\dim X =3$.
Let $\cS'$ be
the germ of an analytic
surface in $X'$ at $x' \in E$ such that $\cG'=\cS' \cap E$ is
the germ of a curve which meets 
$\sigma^{-1} (x)$ at $x'$ only, where $x=\sigma (x')$. 
Then $\sigma (\cS' )$ is
the germ of an analytic surface at $x \in X$.}

\proof By the assumption $x \notin \sigma (\de \cS' )$ where
$\de \cS'$ is the boundary of $\cS'$. Take a small neighborhood
$V$ of $x$ in $X$ which does not meet $\sigma (\de \cS' )$
and consider the intersection $\cS$ of $\sigma (\cS' )$ and $V$.
Show that it is closed.
Consider a sequence
of points $x_i \in \cS, i=1,2, \ldots$ which converges
to $x_0 \in V$. If this sequence is contained in $\cG =\sigma (\cG')$
then $x_0$ belongs to $\cG$ since it follows from the assumption
that the morphism $\sigma |_{\cG'} : \cG' \to \cG$
is finite. Suppose that none of the points from the
sequence is in $\cG$ (and, therefore, $D$). Since
$\sigma |_{X'\setminus E} : X'\setminus E \to X\setminus D$
is an isomorphism $x_i'=\sigma^{-1}(x_i)$ is a point in $\cS'$.
Let $\{ x_i' \}$ converges to a point $x_0'$. Note that
$x_0' \notin \de (\cS' )$ since otherwise $\sigma (\de \cS' )$
meets $V$. Thus $x_0' \in \cS'$ whence $x_0=\sigma (x_0')$
belongs to $\cS$. This implies that $\cS \setminus D$
is a closed analytic subset of
$V \setminus D$ and none of the irreducible components
of the germ of $D$ at $x$ is
contained in $\cS$.
Then by Thullen's theorem (e.g., see
[GrRem], th. 2.1) or by Remmert's theorem 
(e.g., see [BeNa], thm. 1.2) $\cS$ is the germ of an
analytic hypersurface at $x\in X$.
\qed

\medskip

\defin We say that $X$ is a local holomorphic UFD if for 
every $x \in X$ the ring of germs of holomorphic functions
at $x$ is a UFD (note that when $X$ is smooth it is
a local holomorphic UFD by the theorem of Auslander and
Buchsbaum, e.g.
see [Ei, Ch. 19, th. 19.19]).  Let $\dim X =3$.
We say that $X$ is a local holomorphic UFD with respect
to the modification $A \hookrightarrow A'$ if the defining
ideal of every germ of analytic surface $\cS$ as in Lemma 2.2
is principal in the ring of germs of holomorphic functions.
Of course, if $X$ is a local holomorphic UFD then it is
a local holomorphic UFD with respect to $A \hookrightarrow A'$.

\medskip

\prop {\it Let $A \hookrightarrow A'$ be an affine modification
of normal affine domains, and let $f=g^n, g \in A$.
Suppose that $A \hookrightarrow A_1$ is another affine
modification with locus $(J,f)$ such that $A_1 \subset A'$.
Let the defining ideal of the exceptional divisor $E_1$
in $A_1$ is generated by 
$g_1=g\circ \delta_1$ where $\delta_1 : X_1 \to X$ is 
the associate morphism of algebraic varieties.
Suppose that $X$ is a local holomorphic UFD
with respect to $A \hookrightarrow A'$. 
Then $X_1$ is a local holomorphic UFD
with respect to $A_1 \hookrightarrow A'$.}

The proof of this Proposition is the exact repetition of
the proof of Proposition 2.13 (2) with $S$ and $S'$ replaced
by $\cS$ and $\cS_1 = \delta_1^{-1}(S)$.

\rem The author believes that if $X$ is a UFD then it is
automatically a local holomorphic UFD.
This fact is likely known
but we could not find a reference, that is why we need
Proposition 2.14 besides Proposition 2.13. Furthermore,
Proposition 2.14 provides us an additional fact which will be
useful later. Let $\cV = \cS \cap D$ and let $x_1$
be the center of the germ $\cS_1$, i.e. $\delta_1 (x_1)=x$
where $x$ is the center of $\cS$. 
Since $X$ is a local holomorphic UFD
with respect to $A \hookrightarrow A'$ there is a small
neighborhood $V$ of $x$ in $X$ such that we can suppose that
the defining ideal of $\cS$ is principal in the ring of
holomorphic functions on $V$. It follows from the
proof the defining ideal of $\cS_1$ is principal
in the ring of holomorphic functions on $V_1 = \delta_1^{-1}(V)$.
Note that $V_1$ is not already a small neighborhood of $x_1$.
It contains the set $\delta_1^{-1}(\cV )$.

\bigskip

\subsection{Preliminary Decomposition}
We shall fix first notation for this subsection.

\medskip

\conv (1) When we speak about the modification $A \hookrightarrow A'$
in this subsection we suppose that $f=g^n$ where
$g\in A$, the zero multiplicity of $g$ at generic points
of each irreducible component of $D$ is one,
$E$ is non-empty irreducible, and $I=I_f$.
Furthermore, we suppose that the defining ideal $\cI_{A'}(E)$
of $E$ in $A'$ is generated by $g$. Recall that 
by Proposition 2.4 the last condition
holds when both $A$ and $A'$ are UFDs and $g$ is irreducible.

(2) Furthermore, we shall consider affine domains $A_i =\C [X_i], \, i \geq
0$ in this subsection such that
$A\hookrightarrow A_i \hookrightarrow A'$. These embeddings generate
morphisms of algebraic varieties $\delta_i : X_i \to X$ and $\rho_i :
X' \to X_i$ such that $\sigma = \delta_i \circ \rho_i$.
By Remark 2.4 there exist ideals $I_i$ in $A$ and
$K_i$ in $A_i$ such that $A\hookrightarrow A_i$ is an affine 
modification with locus $(I_i,g^{n_i})$ for some $n_i >0$ and
$A_i \hookrightarrow A'$ is an affine modification with 
locus $(K_i,f)$. Hence the exceptional divisor $E_i$ of the first
modification coincides with the divisor $D_i$ of the second one.

(3) We suppose that $K_i$ is the $f$-largest ideal 
of the modification $A_i \hookrightarrow A'$ whence by Proposition 2.2
the closure of the geometrical center $C_i$ of
$\rho_i$ coincides with its reduced center $\bC_i$.
%(4) We shall deal with the situation when $A_i \hookrightarrow A_{i+1}$
%is an affine modification with locus $(J_i,g)$. The corresponding
%morphism of affine algebraic varieties will be denoted by 
%$\sigma_{i+1} : X_{i+1} \to X_i$.

\medskip

\lemma {\it 
Let $A_1 \hookrightarrow A'$ be an affine modification as in
Convention 2.2. Suppose that $A_1$ is normal and 
the closure of $C_1=\rho_1 (E)$ in $X_1$ is an irreducible
component $D^1_1$ of $D_1$.
Let $E_0$ be the Zariski open subset of $E$ that consists
of all points $x' \in E$ such that
there exists a neighborhood
of $x'$ in $E$ which contains no other points
from $\rho_1^{-1} (\rho_1 (x'))$ but $x'$. Put $D_0= \rho_1 (E_0)$
and let $D^2_1$ be the union of irreducible components of $D_1$
different from $D^1_1$. Then

(i) $D_0 = D^1_1\setminus D^2_1$ and $E_0=\rho_1^{-1}(D_0)$;

(ii) the restriction of $\rho_1$ to $(X'\setminus E) \cup E^0$
is an isomorphism between this variety and
$(X_1 \setminus D_1)\cup D_0$;

(iii) in particular, if $E=E_0$ (this is so, for instance,
when $D^2_1$ does not meet $D^1_1$)
then $\rho_1$ is an embedding, and if $D_1 =D_1^1$
then $\rho_1$ is an isomorphism.}

\proof
Since the restriction of $\rho_1$ to $X' \setminus E$
is an embedding, for every $x' \in E^0$
there exists a Zariski open neighborhood
$V_{x'}$ of $x'$ in $X'$ which  contains no other points
from $\rho_1^{-1} (\rho_1 (x'))$ but $x'$.   
Put $x_1=\rho_1 (x')$.
Since $X_1$ is normal $x_1$ cannot
be a fundamental point of the birational map 
$\rho_1^{-1}$ by the Zariski Main
Theorem [Ha, Ch. 5, Th. 5.2]. That is, $\rho_1^{-1}$ is a morphism
in a neighborhood of $x_1$ whence $\rho_1^{-1}$ is an embedding
in this neighborhood which proves (ii).

We denote by $X_1^0$ the algebraic subvariety
$(X_1\setminus D_1)\cup (D_1^1 \setminus D^2_1)$.
Note that the complement to $(X_1\setminus D_1)\cup D_0$ in $X_1^0$ is
a constructive subset of codimension at least 2. Since $X'$ is affine
and $X_1^0$ is normal we can extend morphism $\rho_1^{-1}$ to 
a holomorphic map from $X^0_1$ to $X'$
[Rem, Lemma 13.10] which is regular (e.g., see [Ka2]).
This implies that $\rho_1^{-1}|_{X^0_1} : X^0_1 \to 
X'$ is an embedding whence
$D_0 \supset D_1^1 \setminus D^2_1$. In order to show the
reverse inclusion assume that $x' \in E_0$ and $x_1 = \rho_1 (x')$
is a point from $D_1^1 \cap D_1^2$. 
Since $\rho_1^{-1}$ is an embedding in a neighborhood of $x_1$
we see that the exceptional divisor of $\rho_1$ must contain
a component different from $E$.
This contradiction yields (i). The last statement follows 
immediately from (i) and (ii).
\qed

\medskip

We shall need the following technical notion.

\medskip

\defin Let $A \hookrightarrow A'$ be an affine modification,
$A \hookrightarrow A_1$ be a basic modification
such that $A_1 \subset A'$,
and $S = \{ h^n \mid \, n \in \N \}$ be a multiplicative
system in $A$ where $h \in A$.  Suppose that
$(h \circ \sigma )^{-1} (0)$ does not contain $E$
%any irreducible
%component of $E$ 
and $S^{-1}A_1 = S^{-1}A'$.
Then we call
$A \hookrightarrow A'$ a pseudo-basic modification
(with respect to $A \hookrightarrow A_1$).
That is, this pseudo-basic modification becomes 
the basic modification $S^{-1}A \hookrightarrow S^{-1}A'$
after the localization.

Note that if the assumption of Lemma 2.3 holds and
$A \hookrightarrow A_1$ is basic then it follows from this Lemma
that $A \hookrightarrow A'$ is pseudo-basic.

\medskip

\lemma {\it Let the assumption of Lemma 2.3 and Convention 2.2 hold, and
let $A \hookrightarrow A_1$ be a basic modification
of rank $s \geq 1$ (i.e. $A \hookrightarrow A'$ is pseudo-basic). 
Suppose that the reduced center $C$
of $\sigma$ is a connected component of the reduced
center of this modification $A \hookrightarrow A_1$.
Then the modification $A \hookrightarrow A'$ is locally basic.
Furthermore, if the reduced center of $A \hookrightarrow A_1$
coincides with $C$ then $\rho_1 : X' \to X_1$ is an isomorphism,
i.e. $A \hookrightarrow A'$ is basic.}

\proof
Since $C$ is a connected component of the reduced center of
$A \hookrightarrow A_1$ and $\delta_1$ is basic
(and, therefore, cylindrical by Proposition 2.7)
the exceptional divisor $E_1=D_1$
of $\delta_1$ is of form $D_1^1 \cup D^2_1$ where 
$D^1_1 = \delta_1^{-1} (C) \simeq \C^s \times C$ and $D^2_1$ does not meet
$D^1_1$. Lemma 2.3 (iii) implies now the desired conclusion.
\qed

\rem (1) Note that if $C$ is a point and belongs to the reduced
center of $\delta_1$ then $C$ is automatically a connected
component of the reduced center of $\delta_1$.

(2) By Proposition 2.12 instead of normality of $A_1$
one can require that $A$ is normal Cohen-Macaulay.

\medskip

\prop {\it 
Let $A\hookrightarrow A'$ and $A_i$ be as in Convention 2.2.
Suppose also that
$C$ is not contained 
in the singularities of $X$, its codimension in $X$ is
at least 2, and
the zero multiplicity of $g$ at generic points of $C$
is 1.  Let
$$A=A_0 \hookrightarrow \dots
\hookrightarrow A_{k-1} \hookrightarrow A_k, \, k\geq 0$$
is a strictly increasing sequence of affine domains such that
$A_k \subset A'$, and for every $i\leq k$ 

(i) the embedding
$A_{i-1} \hookrightarrow A_i$ is a basic modification with
locus $(J_i,g)$ and of rank $s_{i-1}$ where $s_{i-1}+1$ is the codimension
of $C_{i-1}$ in $X_{i-1}$ (see Convention 2.2 (3) for the definition of
$C_i$).

(1) Then $k \leq n$ (recall that $f=g^n$) and this sequence
can be extended to a strictly increasing sequence of affine domains
$$A_0 \hookrightarrow \dots
\hookrightarrow A_{m-1} \hookrightarrow A_m, \, k\leq m\leq n$$ 
for which (i) holds for every $i\leq m$,
and $A_{m-1} \hookrightarrow A'$ is pseudo-basic with respect
to $A_{m-1} \hookrightarrow A_m$.

(2) Suppose that $\sigma_i : X_i \to X_{i-1}$ is the
morphism associated with  the affine modification
$A_{i-1}\hookrightarrow A_i$. 
Then $\sigma_i (C_i)=C_{i-1}$ for $i\leq m-1$, and
$\rho_{m-1} (E)=C_{m-1}$.  

(3) Furthermore,
suppose that $A$ is normal Cohen-Macaulay, $A'$ is normal, and
that the closure $E_m^1$ of $\rho_m (E)$ is a connected component
of $E_m$ (resp. $E_m$ is irreducible).
Then $\rho_{m-1}$ is a locally basic (resp. basic) modification.}

\proof Let us show (2) first.
By Convention 2.2 (2) and (3)
the exceptional divisor of $\rho_i$ is $E$ whence
$\rho_i(E)=C_i$.
In particular, $\rho_{m-1} (E)=C_{m-1}$ and
$\sigma (E)=C_0$ for $\sigma = \rho_0$. 
Since $\sigma = \rho_i \circ \delta_i$
we see that $\delta_i (C_i)=C_0$.
This implies that $\sigma_i (C_i)=
C_{i-1}$ since $\delta_i  = \delta_{i-1} \circ \sigma_i$.

If $s_k=0$ then we put $m=k$ and get (1) automatically. Otherwise,
let us show now that the assumptions of this Proposition hold if we replace
$A  \hookrightarrow A'$ by $A_k  \hookrightarrow A'$.
It is enough to check this for $k=1$.
By Remark 2.8 for a generic point $x \in \rho_0 (E)=\sigma (E)$
the points from $\sigma_{1}^{-1}(x)$ are smooth in $X_1$.
Since $\rho_0 = \sigma_1 \circ \rho_1$
we see that $\sigma^{-1}_{1} (x)$
contains generic points of $\rho_{1} (E)$. Hence
generic points of $C_1$ are not
contained in the singularities of $X_{1}$, and
the zero multiplicity of the function
$g\circ \delta_1$ at these generic points is 1 by
Remark 2.8.
By Lemma 2.1 and Proposition 2.3 we can choose
a basic modification $A_k \hookrightarrow A_{k+1}$ with locus
$(J_k,g)$ such that $A_{k+1} \subset A'$ and the rank of the
modification is $s_k$.
Thus we can extend our strictly increasing sequence
of affine domains and we can always suppose that $k\geq 1$ in (1).
There are two possibilities: either this sequence becomes eventually
infinite or there exists $m$ such that $s_m=0$
which implies (1).
We need to show that the first possibility does not hold
and that $m\leq n$ with help of induction by $n$.

Assume first that (1) holds for $n-1>0$ and show it for $n$.
Let $b_0=g,b_1, \dots , b_s$ be a representative system of
generators for $A \hookrightarrow A_1$. By assumption of this
Proposition and 
Definition 2.5 there exists $h\in A$ such that $h^{-1}(0)$ does not
contain $\sigma (E)$, 
$X\setminus h^{-1}(0)$ is smooth,
$C\setminus h^{-1}(0)$ is a complete intersection
in $X\setminus h^{-1}(0)$ given by $b_0= \dots =b_s=0$.
If $S$ is the multiplicative system $\{ h^j \, | j \in \N \}$
in $A$ then the affine modification 
$S^{-1}A \hookrightarrow S^{-1}A'$ satisfies the analogue
of assumption of this Proposition and, furthermore,
$S^{-1}J_1$ is the defining ideal of $C\setminus h^{-1}(0)$
in $S^{-1}A$.
Since the locus of the affine modification
$S^{-1}A \hookrightarrow S^{-1}A_1$ is $(S^{-1}J_1,g)$ 
by Proposition 2.1,
Corollary 2.1 implies that for the affine modification
$S^{-1}A_1 \hookrightarrow S^{-1}A'$ the locus can be chosen
in the form $(L_1,g^{n-1})$.
Thus after the localizations of our strictly increasing sequence
of affine domains with respect to $S$ we have by the induction assumption
that the codimension of the reduced center of
the modification $S^{-1}A_m \hookrightarrow S^{-1}A'$
is 1 for some $m\leq n$. This implies that the same is
true for the reduced center of $A_m \hookrightarrow A'$,
i.e. $s_m=0$ which concludes this step of induction.

The next step of induction is for $n=1$.
By Proposition 2.3 in this case $S^{-1}J_1$ coincides with the $g$-largest
ideal of the affine modification $S^{-1}A \hookrightarrow S^{-1}A'$.
Hence $S^{-1}A_1=S^{-1}A'$. Since $h$ is chosen so
that $h^{-1}(0)$ does not contain $C$ this implies (1)
which concludes induction.

Note that when $A$ is normal Cohen-Macaulay so is $A_k$ by 
Propositions 2.9, 2.10, and 2.12.
Claim (3) is now a consequence of Lemma 2.4.
\qed

\medskip

Let $C_{m-1}^*$ be the complement in $C_{m-1}$ to the set
of points where $C_{m-1}$ meets the other components of
the reduced center of $\sigma_m$. Then the exceptional
divisor of the basic modification $\sigma_m$ contains
$E_m^* \simeq \C^{s_{m-1}} \times C_{m-1}^*$. Furthermore,
under the assumption of Proposition 2.15 (3) the restriction
of $\rho_m^{-1}$ to $(X_m \setminus E_m) \cup E_m^*$
is an embedding by Lemma 2.3. Hence

\medskip

\cor (cf. [Miy2, Lemma 2.3]) {\it Under the assumption of
Proposition 2.15 (3) the exceptional divisor $E$
contains a Zariski open cylinder $E_m^* \simeq
\C^{s_{m-1}} \times C_{m-1}^*$ such that $\rho_{m-1}|_{E_m^*}$
is the projection to the second factor.}

\bigskip

\subsection{Decomposition}
In this subsection we shall strengthen Proposition 2.15
in the case when $\dim X =3$ and $X$ is a holomorphic UFD.
Our main aim is to make $A_m=A'$ in this Proposition.
By Lemma 2.4 it is enough to require that the reduced
center of $\sigma_m$ is irreducible, i.e. it coincides
with $\bC_{m-1}$. This is true when the defining
ideal of $\bC_{m-1}$ in $\C [D_{m-1}]$ is principal
provided $C_{m-1}$ is a curve,
or more generally the defining ideal of each $C_i$
in $\C [D_i]$ is principal provided $C_i$ is a curve.
The last claim will be proven
by induction and the first step of induction is crucial.
%It turns out that under a mild
%additional assumption 
%the defining ideal of $C_1$ in $\C [D_1]$ is principal
%when $C_1$ is a curve.
But the proof of this step is difficult
and we present it in the next section (Proposition 3.1).
Another non-trivial fact whose proof is postponed till
next section says that the number of irreducible components
of the germ of $C_{m-1}$ at each point $z \in C_{m-1}$
coincides with the number of connected components in
$\rho_{m-1}^{-1}(z)$ (Lemma 3.3).
This helps us to show that $C_i$'s are contractible in
some cases. Furthermore, if we want to check that $C_i$'s
are smooth we need the following

\medskip

\lemma {\it 
Let Convention 2.2 hold and let
$A_1 \hookrightarrow A'$ be a basic modification
with the divisor $D_1$ of modification isomorphic to
$\C \times C$ where $C$ is a curve.
Suppose that $C$ has an irreducible singular point $z$
and that $C_1$ meets the line
$\C \times z$ in $D_1$ at $z_1=0 \times z$
but $C_1$ is different from this line.
Then $z_1$ is a singular point of $C_1$.}

\proof 
Assume the contrary, i.e. $C_1$ is smooth at $z_1$.
Since the situation is local we can suppose that $C$
is a closed curve in $\C^n$.
Consider a normalization $\nu_0 : C^{\nu} \to C$. It generates
a morphism $\nu =({\rm id},\nu_0 ): 
\C \times C^{\nu}  \to \C \times C \subset \C^{n+1}$.
Suppose that $(y, \bx ) =(y, x_1, \dots , x_n)$ is a coordinate system
in $\C^{n+1}$.
Let $g,b_1$ be an almost complete intersection in $A_1$
which generates this basic modification
$A_1 \hookrightarrow A'$, i.e. $b_1$ generates the defining
ideal of $C_1$ in $\C [D_1]$. We
treat $b_1$ as a polynomial $b_1(y,\bx )$ on $\C^{n+1}$.
Let $\beta = b_1|_{D_1} \circ \nu$, $C^{\nu}_1$ be 
the proper transform of $C_1$ (i.e. $C^{\nu}_1 = \beta^{-1}(0)$),
and let $o = \nu^{-1} (z_1)$. Since $C_1$ is smooth (and, therefore,
normal) and since
$\nu |_{C^{\nu }_1} : C^{\nu}_1 \to C_1$ is a
homeomorphism (i.e., this morphism is proper and
finite), $C^{\nu}_1$ is biholomorphic to $C_1$
(e.g., see [Pe1, Cor. 1.5])
whence $C^{\nu}_1$ is smooth at $o$.
Since the modification is basic the condition
on the gradients implies that the gradient of $\beta$
does not vanish at generic points of $C^{\nu}_1$.
Hence since $C^{\nu}_1$ is smooth at $o$ the gradient of $\beta$
does not vanish at $o$. Let $(v,t )$ be a local
coordinate system at $o$ where $t$ is a coordinate
on the second factor and $v$ is a coordinate on
the first factor of $\C \times C^{\nu} $. 
In particular, locally $\nu (v,t) = (v, \bx (t))$
whence the Taylor series of $\beta (v,t)=b_1(v, \bx (t))$ at
$o$ does not have a nonzero linear term $ct$ with $c \in \C$
(recall that $z$ is a singular point of $C$ whence none
of $x_i(t)$'s contains a linear term $ct$).
Hence, the linear part of this
power series must be $v$ up to a nonzero constant
factor (the factor is nonzero since otherwise the
gradient of $\beta $ at $o$ is zero).
Thus the Taylor series of $b_1$ at $z_1$ has a nonzero
linear term $cv$.
The implicit function theorem implies
that the germ of $C_1$ at $z_1$ is isomorphic to the germ
of $C$ at $z$ whence $C_1$ is singular at $z_1$.
Contradiction.
\qed

\medskip

\thm {\it Let $A \hookrightarrow A'$ be an affine modification
such that  Convention 2.2 (1) holds
(recall that this is true when $E$ is non-empty,
$A$ and $A'$ are UFDs and $f=g^n$ where $g \in A$ is irreducible),
$A$ is Cohen-Macaulay, and
$\dim X =3$. 
Suppose that either

($\alpha$) $X'$ is smooth
and $H_3 (X')=H_2(X' \setminus E)=H_3 (X' \setminus E)=0$, or

($\beta$) $E$ is a UFD and its Euler characteristics is 1.
%($\gamma$) $X'$ is smooth, $H_3 (X' \setminus E)=0$, $E$ has at most isolated
%singularities and Euler characteristics 1.

Suppose also that

(i) $D$ is isomorphic to $\C^2$;
\footnote{One can replace (i) with the condition that $D$ is a UFD.
In this case the statement of the theorem remains the same with
one exception : when $C_0$ is a point
and we want the center of $\sigma_1$ to be the defining
ideal of $C_0$ in $A$ then we have to allow
$\sigma_1$ to be not necessarily basic but only locally basic.}

(ii) $A$ is a local holomorphic UFD (see Definition 2.10) with respect to
the modification $A \hookrightarrow A'$ (which is true when
$X$ is smooth);

(iii) $C$ is not contained in the singularities of $X$.

Let $m,A_i,C_i$, and $\bC_i$ be the same as in Proposition 2.15
and Convention 2.2 (3).

(1) Then the algebras $A_i$'s can
be chosen so that $A_m=A'$, $C_i = \bC_i$ for every $i$, and
if $C_i$ is a curve its
defining ideal in $\C [D_i]$ is principal.

(2) Furthermore, each $C_i$ is either a point
or an irreducible contractible curve, and in the case when
$E$ has at most isolated singularities 
(which is true under condition ($\beta$)) these contractible
curves are smooth.}

\proof We use induction by $m$.
Suppose first that $C_0$ is a point. The assumption (i)
on $D$ and the fact that $C_0$ is a smooth point 
of $X$ by (iii)
allow us to choose $b_1, b_2 \in A$ such that
$g,b_1,b_2$ generate the defining ideal $I_1$ of $C_0$ in $A$.
Hence the exceptional divisor $E_1=D_1$ of the basic
modification $A \hookrightarrow A_1= A[I_1/g]$ with
locus $(I_1,g)$ is isomorphic to $\C^2$ and, in particular,
irreducible. Note that
$A \subset A'$ by Proposition 2.3 and Convention 2.2 (1).
If $m=1$ in Proposition 2.15 then
Lemma 2.3 implies that $\rho_1$ is an isomorphism,
and, in particular, $A_1 =A'$.

Suppose now that $m\geq 2$.
Note that $E_1$ is isomorphic to $\C^2$ and it is the
divisor of the modification $A_1 \hookrightarrow A'$
described in Convention 2.2.
By Propositions 2.10 and 2.12 $A_1$ is normal Cohen-Macaulay.
By Proposition 2.13 $A_1$ is a UFD and furthermore
it is a local holomorphic UFD with respect to
the modification $A_1 \hookrightarrow A'$ by Proposition 2.14.
Thus the assumptions of this Theorem hold also for
the modification $A_1 \hookrightarrow A'$. The decomposition
of this last modification into basic
modifications contains $m-1$ factors
and induction implies the desired conclusion in this case.    

If $\bC_0$ is a curve then its defining ideal in $\C [D]$
is principal whence the defining ideal $I_1$ of $\bC_0$ in $A$
is generated by $g$ and $b\in A$. The exceptional divisor
$E_1$ of $A \hookrightarrow A_1= A[I_1/g]
\subset A'$ is again irreducible, and $A_1 \subset A'$ as before.
If $m=1$ then Lemma 2.3 implies that $\rho_1$ is an isomorphism,
i.e. $A_1 =A'$.
Furthermore, for every $z \in C_0$ the number of components in
$\sigma^{-1}(z)$ is one, since $\sigma^{-1}(z)
\simeq \sigma_1^{-1}(z) \simeq \C$ (recall that $\sigma_1$ is basic).
Hence Lemma 3.3 below implies that the number of irreducible
components of the germ of $C_0$ at $z$ is one,
i.e. $z$ is not a double point of $C_0$.
The same Lemma says that the normalization of $C_0$ is $\C$
whence $C_0$ is contractible. This means that $C_0$ is
closed in the ambient affine algebraic variety $D$, i.e.
$C_0 =\bC_0=C$. 
When $C_0$ has singularities  then
$E$ has singularities
in codimension 1 which yields the last statement.

Actually in the case when $C_0$ is a curve,
instead of (i) one can assume that
$D$ is isomorphic to $\C \times G$ where $G$ is a curve
and the natural projection $C\to G$ is
dominant (this is, of course, true for $D \simeq \C^2$).
This curve $C_0$ is contractible by Proposition 3.1 below
(which implies, in particular, that $C =C_0$)
whence the projection $C\to G$ is finite.
The defining ideal of $C_0$ in $\C [D]$ is generated
by a function $b \in A$ (see Proposition 3.1 
below) whence
the defining ideal $I_1$ of $C$ in $A$ is generated
by $g$ and $b$. Thus the exceptional
divisor $E_1$ of the basic modification
$A \hookrightarrow A_1 =A[I_1/g]$ is again irreducible
and we see as before that $A_1$ is a Cohen-Macaulay UFD
which is also a local holomorphic UFD with respect to
the modification $A_1 \hookrightarrow A'$.
Thus when $m \geq 2$ the assumption of this
Theorem holds for $A_1 \hookrightarrow A'$ and
induction implies statement (1).

The curves $C_i$'s are contractible for $i\geq 2$ by the induction assumption.
Lemma 2.5 implies that if $C_0$ is not smooth then
$C_1$ cannot be smooth which yields the last statement of (2).
\qed

\bigskip

\section{The Geometry of The Exceptional Divisor and The Reduced
Center}

\subsection{The Exceptional Divisor}
We shall finish the proof of Theorem 2.3 in
this section. First we 
describe $E$ in the three-dimensional case
more accurately then Corollary 2.2 does.

\medskip

\lemma {\it Let $A \hookrightarrow A'$ be an affine modification
such that $E$ is irreducible and the geometrical center
$C_0$ of this modification is a curve.
%Suppose that Convention 2.2 holds and 
Suppose that $X, X'$ are of dimension 3.
Let $A \hookrightarrow A_1$ be a semi-basic modification of rank 1
and with locus $(I_1,f)$
such that $A_1 \subset A'$.
Then for 
%every $x\in C_0$
%each irreducible component of $\sigma^{-1}(x)$ is a curve and 
every finite subset $R \subset C_0$ there
exists a semi-basic modification $A \hookrightarrow A_2$ of rank 1
with locus $(I_2,f)$ such that $A_2 \subset A'$,
and $\rho_2$ is not constant on each of irreducible components
of $\sigma^{-1}(R)$ where $\rho_2$ is as Convention 2.2 (2).}

\proof 
Let a semi-regular sequence $f,b_1 \in A$ generate $I_1$, i.e.
$I_1 \subset I$ and $A_1=A[I_1/f]$.
Consider several basic modifications of this type.
That is, for $j=1, \dots ,k$ the sequence
$b_0=f, b_j$ is semi-regular and it
generates an ideal $I_j \subset I$. Let 
$A \hookrightarrow A_j=A[I_j/f]$ be an affine modification with
locus $(I_j,f)$. Recall that
$\delta_j : X_j \to X$ is the corresponding morphism of
affine algebraic varieties and
$E_j$ is the exceptional
divisor of $\delta_j$. These morphisms $\{ \delta_j \}$
define an affine variety $Y=X_1 \times_{X} X_2 \times_{X} \cdots
\times_X X_k$
and its subvariety $Y^*=(X_1\setminus E_1)\times_{X} \cdots \times_X
(X_k\setminus E_k)$. Since we can perturb the elements $b_j$
of our semi-regular sequences by Lemma 2.1 (4), we
can suppose that $I_1+ \dots + I_k=I$.
By Proposition 2.6 $X'$ can be viewed as
the closure $\bY^*$ of $Y^*$ in $Y$,
$\sigma$ can be viewed as the restriction of the natural projection
$Y \to X$ to $\bY^*$, and $\rho_j$ (from Convention 2.2 (2))
will be nothing but
the restriction of the natural projection $Y \to X_j$ to 
$\bY^*$.

Note that $X_j$ can be viewed as a closed subvariety of
$\C \times X$ such that
a coordinate $y_j$ on the first factor
of $\C \times X$ is chosen so that its
restrictions to $X_j$ coincides with $b_j/f$, and
$\delta_j$ coincides with the
restriction of the natural projection
$\C \times X \to X$.
Thus $X'$ can be viewed as a closed subvariety of
$\C^{k} \times X$ such that $(y_1, \ldots ,  y_k)$
is a coordinate system of the first factor,
the restriction of $y_j$ to $X'$ is $b_j/f$, and
$\sigma$ coincides with the restriction of the natural
projection $\C^{k} \times X \to X$ to $X' \, (=\bY^*)$.

Let $x \in  C_0$. 
Since $\dim \sigma^{-1}(x_1) =1$ for a generic point $x_1 \in C_0$
the dimension of every irreducible
component of $\sigma^{-1}(x)$ is at least 1 by the semi-continuity
of the dimension of the fibers of an algebraic morphism.
But this dimension cannot be equal to 2 since otherwise
$E$ contains at least two irreducible components in
contradiction with the assumption of this Proposition.
Thus $\sigma^{-1}(x)$ is a curve.
Since $\sigma (\sigma^{-1}(x)) =x$
the image of $\sigma^{-1}(x)$ in $\C^k$ under the
natural projection is a curve, i.e. the natural
projection of $\sigma^{-1}(x)$ to a generic
affine line in $\C^{k}$ is dominant.
By Lemma 2.1 after a small perturbation of
$b_2$ in $I$ the sequence
$b_0, b_2$ remains a semi-regular sequence,
i.e. we can suppose that the $y_2$-axis is a generic line.
Thus if $\tau : X' \to \C$ is the natural
projection to the $y_2$-axis then we can suppose that
$\tau$ is not constant on every irreducible component 
of $\sigma^{-1}(R)$. 
Note that $\tau = \theta_2 \circ \rho_2$ where
$\theta_2 : X_2  \to \C$ is the
the natural projection to the $y_2$-axis.
Hence the restriction of $\rho_2$ to every component of
$\sigma^{-1}(R)$ is not constant.
\qed

\lemma {\it Let Convention 2.2 hold (in particular,
$C_0$ is the geometrical center of $\sigma$) and
$X_1$ and $X'$ be normal varieties of dimension 3.
Let $A \hookrightarrow A'$ be a pseudo-basic modification
with respect to a basic modification $A \hookrightarrow A_1$.

(1) Then for every $z \in C_0$
the curve $\sigma^{-1}(z)$ is a disjoint union
of irreducible contractible curves.

(2) If $E$ has no double points then it is a topological
manifold.

(3) If in addition to (2) for every $z \in C_0$ the curve $\sigma^{-1}(z)$
is connected then $E$
is naturally homeomorphic to the product of $\C$ and a curve.}

\proof By Corollary 2.2 there exists $C^* \subset C$
for which $E$ contains a Zariski open cylinder
$E^* \simeq \C \times C^*$ such that $\sigma |_{E^*}$ is
the projection to the second factor. Therefore,
we need to consider only $z$ from the finite set
$R= C_0 \setminus C^*$. 
Let $(g, d)$ be an almost complete
intersection which generates
the basic modification
$A \hookrightarrow A_1$. Then $f=g^n, (b_1)^n$ is a semi-regular
sequence which is contained in $I$. Hence we
can choose a semi-regular sequence $f,b_2$ in $I$ which generates
a semi-basic modification
$A \hookrightarrow A_2$ as in Lemma 3.1.

Consider $Y=X_1 \times_{X} X_2$ and
and its subvariety $Y^*=(X_1\setminus E_1)\times_X
(X_2\setminus E_2)$.
Let $X_0$ be the closure $\bY^*$ of $Y^*$ in $Y$.
Recall that  $X_1$ (resp. $X_2$) can be viewed as a closed subvariety of
$\C \times X$ such that
a coordinate $y_1$ (resp. $y_2$) on the first factor
of $\C \times X$ is chosen so that its
restrictions to $X_1$ (resp. $X_2$) 
coincides with $d/g$ (resp. $b_2/f$), and
$\delta_j$ coincides with the
restriction of the natural projection
$\C \times X \to X$.
Thus $X_0$ can be viewed as a closed subvariety of
$\C^2 \times X$ such that $(y_1, y_2)$
is a coordinate system of the first factor.

Let $\delta_0 : X_0 \to X$, $\tau_i : X_0 \to X_i$
be the natural projections (note that
$\delta_0$ coincides with the restriction of the natural
projection $\C^{2} \times X \to X$ to $X_0$).       
Put $E_i^* = \delta_i^{-1} (C^*)$.
Let $F$ be the closure of $E^*_0$ in $\C^{2} \times X$.
Note that $\delta_0= \delta_i \circ \tau_i$ and $\sigma
= \delta_i \circ \rho_i$. Hence $\rho_1 = \tau_1 \circ \rho_0$.
By Corollary 2.2 $\rho_1 |_{E^*} : E^* \to E_1^*$
is an isomorphism whence
$\rho_0 |_{E^*} : E^* \to E_0^*$  and
$\tau_1|_{E_0^*} : E_0^* \to E_1^*$ are isomorphisms.

Let us show that $\rho_0 (\sigma^{-1}(z))$ is contained in
a disjoint union of contractible lines.
Note that $F \subset X_0$ and $\rho_0(E) \subset F$.
Hence $\rho_0 (\sigma^{-1}(z))$ is contained in
$F\setminus E^*_0$.
Since $\delta_1$ is basic
$E_1^* \simeq \C \times C^*$ where $y_1$ is a coordinate
on the first factor. The surface $F$ is contained
in $\C^2 \times C$ where $(y_1,y_2)$ is a coordinate system
on the first factor. Since 
$\tau_1|_{E_0^*} : E_0^* \to E_1^*$ is an isomorphism
and coincides with the restriction of the natural
projection $\C^2 \times C \to \C \times C, \, ((y_1,y_2),x) 
\to (y_1,x)$ for every $x\in C^*$ the equation
of $\tau^{-1}_1 (x)$ in $\C^2_{y_1,y_2}$ is of form
$\ba(x)y_2 + \sum_{i=0}^{k}\ba_i(x)y_1^i=0$ where
$\ba, \ba_i$ are regular functions on $C^*$ and $\ba$
is invertible. Consider $z \in C\setminus C^*$ and
an irreducible branch $\cC$ of of the germ of $C$ at $z$.
Let $\ba'$ be one of the restrictions of
$\ba_i$'s or $\ba$ to
$\cC \setminus z$ which has the largest
pole at $z$. Dividing by $\ba'$ we see that
the curves $\ba(x)y_2 + \sum_{i=0}^{k}\ba_i(x)y_1^i=0, x \in \cC$
approach to a curve 
$$a(z)y_2 + \sum_{i=0}^{k}a_i(z)(y_1)^i=0 \eqno (1)$$
where $a(x)$ and $a_i(x)$'s are rational continuous functions on $\cC$
which are regular on $\cC \setminus z$.
Not all coefficients before $y_i$'s in
equation (1) are zeros and we have three cases.

Case 1. If $a_0(z) \ne 0$ and the rest of coefficients
are zeros then equation (1) defines an empty set.

Case 2. If $a(z)=0$ and some $a_i(z) \ne 0$ for $i\geq 1$
then equation (1) defines a set of line parallel to the
$y_2$-axis.

Case 3. If $a(z) \ne 0$ then equation (1) defines a
contractible irreducible
curve $T$ such that $y_1$ is a coordinate on $T$. 
If $T$ is the closure of the image of an irreducible
component of $\sigma^{-1}(z)$ under $\rho_0$ then
the restriction of $\rho_1$ to this component gives
a dominant morphism to $\delta_1^{-1}(z)$
(this last curve is isomorphic to $\C$ since $\delta_1$
is basic). By Lemma 2.3 (and in its notation) 
$\delta_1^{-1}(z)$ is contained in $D_0$
whence $\rho_1^{-1}$ is an embedding in a neighborhood
of $\delta_1^{-1}(z)$. Since $\rho_1
=\rho_0 \circ \tau_1$ we see that $\rho_0^{-1}$ is
an embedding in a neighborhood of
$\delta_0^{-1}(z)$ and $\tau_1^{-1}$
is an embedding in a neighborhood
of $\delta_1^{-1}(z)$.
That is, in this case $z$ can
be treated as a point of $C^*$.
Thus  $\rho_0 (\sigma^{-1}(z))$ is contained in
a disjoint union of contractible lines.

Note that every irreducible component $T'$ of $\sigma^{-1}(z)$
is a limit of curves $\sigma^{-1}(x), \, x \in C^*$ which
are isomorphic to $\C$. Thus $T'$ admits a non-constant
morphism from $\C$, i.e. $T'$ is a once punctured curve.
This implies that $\rho_0$ maps $T'$ surjectively
on an irreducible component $T$ of $F\setminus E_0^*$
(recall that $\rho_0$ is non-constant by Lemma 3.1).
Hence it remains to show that $\rho_0|_{E}: E \to F$ is an injection.
Assume the contrary, i.e. there exist
different points $x_1',x_2'$ in $E\setminus E^*$
such that $\rho_0(x_1')=\rho_0(x_2')=x\in T$.
Let $V_i$ be a neighborhood of $x_i'$ in $E$.
Since $\rho_0$ is non-constant on every component of
$\sigma^{-1}(z)$ we see that $\rho_0 (V_i)$ contains
the germ $\cT$ of $T$ at $x$ and for a generic point
$x_0 \in \cT$ a neighborhood of $x_0$ is contained in $\rho_0 (V_i)$.
Thus $(\rho_0 (V_1) \cap \rho_0 (V_2) ) \setminus T$
is not empty.
Since $\rho_0|_{E^*}: E^* \to E_0^*$ is an isomorphism
we see that $V_1$ and $V_2$
meet each other. Hence $E$ is not separable which is
not true since $E$ is affine. Thus
$E^* \to F$ is an injection which implies (1).

Note that if $E$ has no double points then,
since $E^* \to F$ is an injection,
the equation $a(x)y_2 + q(x,y_1)=0$ with
$q(x,y_1) =\sum_{i=0}^{k}a_i(x)(y_1)^i$ and $x$
running over $\cC$
defines a homeomorphic image
of a germ of $E$.
Thus it suffices to check
the statements (2) and (3) for the variety given
this equation in Case 2 (two other cases are obvious).

Replacing $\cC$ by its normalization
one can suppose that $\cC$ is smooth
and, therefore, may be viewed as the germ of $\C$
at 0. Thus $x$ can be treated as a coordinate
on $\C$ now. If $k$ is the multiplicity of zero of $a(x)$
we can replace this equation $a(x)y_2 + q(x,y_1)=0$ with
$x^ky_2 + q(x,y_1)=0$.
The last equation
defines a variety $Z$. 
Put $q_1(x,y_1) = {\frac {\de } {\de y_1}} q(x,y_1)$.
Let $\{ c_i \}$ be the roots of $q(0,y_1)$ and let
$l$ be the zero multiplicity of $q_1(x,c_1)$.
Let $Z_1$ be the variety obtained from $Z$ by deleting
all lines $x=y-c_i=0$ where $i \geq 2$.
Consider the following locally nilpotent derivation
on the algebra of regular functions on $Z_1$:
$\de (x)=0, \de (y_1)=-x^{k-l}, \de (y_2) =q_1(x,y_1)/x^l$.
It defines a regular $\C_+$-action on $Z_1$ which
maps the line $x=y-c_1=0$ onto itself and acts
transitively on this line. Thus the preimage of
this line in a normalization $Z_1^{\nu}$ of $Z_1$
does not contain singular points of $Z_1^{\nu}$. Indeed, otherwise
each point of this preimage is singular (by the
transitivity of the $\C_+$-action) and the normal
variety $Z_1^{\nu}$ has singularities in codimension 1
which cannot be true. Hence $Z_1^{\nu}$ is smooth.
Since in the absence of double points any normalization is
a homeomorphism $Z_1$ is a topological manifold
which yields (2). Furthermore, it is easy to
see that $Z_1^{\nu}$ is isomorphic to $\cC \times \C$
which implies (3).
\qed

\medskip

\lemma {\it Let $X'$ be an affine threefold with $H_3(X')=0$ and $E$ be a
closed irreducible surface in $X'$ which admits a surjective
morphism $\tau : E \to C_{m-1}$ into a curve $C_{m-1}$ 
\footnote{We denote this curve by $C_{m-1}$ since it will play later
the role of the geometrical center of the modification
$\rho_{m-1}$ from Proposition 2.15 and $\tau$ will be $\rho_{m-1}|_{E}$.}
such that for a Zariski open
subset $C^*_{m-1}\subset C_{m-1}$ and $E^* = 
\tau^{-1}(C^*_{m-1})$ the morphism
$\tau |_{E^*}: E^* \to C^*_{m-1}$ is a $\C$-cylinder and
$L:=E\setminus E^*$ is a disjoint union of irreducible contractible
curves. Let $H_{2}(X'\setminus E) =H_{3}(X'\setminus E)=0$.
Suppose that $z \in C_{m-1} \setminus C^*_{m-1}$ and $\cC_z$ is
the germ of $C_{m-1}$ at $z$ (we treat $\cC_z$ as a bouquet of discs).
Put $\cC^z =\cC_z \setminus z$, $\cE^z = \tau^{-1}(\cC^z)$,
$\cE_z = \tau^{-1}(\cC_z)$,
and $L^z=L \cap \cE_z \, (=\tau^{-1}(z))$. 
Then there exists an isomorphism $H_0 (L) \simeq H_1 (E^*)$
such that for every germ $\cC_z$ as above the restriction
of this isomorphism generates an isomorphism
$H_0 (L^z) \simeq H_1 (\cC^z)$ (in particular,
the number of connected components of $L^z$ is the same
as the number of irreducible components of $\cC_z$). Furthermore,
the normalization of $C_{m-1}$ is $\C$.}

\proof Let $L:= E\setminus E^*$.
Consider the following
exact homology sequences of pairs:
%{\begin{picture}(500,90)
%\put(2,65){$\longrightarrow H_{j+1}(X')
%\rightarrow H_{j+1}(X',X'\setminus L)  \rightarrow 
%H_{j}(X'\setminus L)
%\rightarrow H_{j}(X') \rightarrow 
%H_{j}(X',X'\setminus L)\rightarrow $}
%\put(215,55){$\line(0,-1){20}$}
%\put(220,55){$\line(0,-1){20}$}
%\put(0,12){$ \ldots 
%\longrightarrow H_{j}(X'\setminus E)
%\rightarrow H_{j}(X'\setminus L)
%\rightarrow H_{j}(X'\setminus L, X'\setminus E)
%\rightarrow H_{j-1}(X'\setminus E) \rightarrow .$}
%\end{picture}
$$\longrightarrow H_{j+1}(X')
\rightarrow H_{j+1}(X',X'\setminus L)  \rightarrow
H_{j}(X'\setminus L)
\rightarrow H_{j}(X') \rightarrow
H_{j}(X',X'\setminus L)\rightarrow $$
and
$$ \ldots
\longrightarrow H_{j}(X'\setminus E)
\rightarrow H_{j}(X'\setminus L)
\rightarrow H_{j}(X'\setminus L, X'\setminus E)  
\rightarrow H_{j-1}(X'\setminus E) \rightarrow .$$
Note that $H_4(X')=0$ since $X'$ is an affine algebraic
variety [Mil, th. 7.1]. Taking into consideration the other assumptions on
$H_{j}(X'\setminus E)$ and $H_{j}(X')$ 
and Thom's isomorphisms $H_0 (L) \simeq
H_{4}(X',X'\setminus L)$ and
$H_1 (E^*) \simeq H_3(X'\setminus L,X'\setminus E)$
(e.g. see [Do, Ch. 8, 11.21]) we have
$$H_0(L) \simeq H_{4}(X',X'\setminus L) \simeq H_{3}(X'\setminus L)
\simeq H_{3}(X'\setminus L, X'\setminus E)\simeq H_1(E^*).$$
Note also that $H_1 (C^*_{m-1})\simeq H_1 (E^*)$ whence we have
an isomorphism between $H_0(L)$ and $H_1 (C^*_{m-1})$.

Let us have a closer look at this isomorphism.
Suppose that $L_i$ is an irreducible component of $L$. Consider
the germ $\cS_i'$ of a smooth complex surface
which is transversal to both $E$ and 
$L_i$ at a smooth $z'$ point of $L_i$. Recall that
the Thom class is the element $u_i \in H^4 (X', X'\setminus L_i)$
uniquely defined by the condition $u_i (\cS_i')=1$.
The Thom isomorphism $H_{4}(X',X'\setminus L) \to H_0 (L)$
is defined by the cap-product with the Thom class
$u= \sum_i u_i$ : $\eta \to u \cap \eta$
(in particular, $\cS_i'$ as an element of $H_{4}(X',X'\setminus L)$
is mapped under this isomorphism to the positive 
generator of $H_0 (L_i)$). We can suppose that 
$\cS_i'$ is diffeomorphic to a ball and its boundary $\de \cS_i'$
in $X'$ is diffeomorphic to a three-sphere. The isomorphism
$H_{4}(X',X'\setminus L) \simeq H_{3}(X'\setminus L)$ sends
$\cS_i'$ to $\de \cS_i'$ (which is viewed as
an element of $H_{3}(X'\setminus L)$). 
Then we can suppose that $\cS_i'$ meets $E^*$ transversally along
a disjoint union of the germs of curves $\G_i = \bigcup_j \G_i^j$,
that the boundary of each germ $\G_i^j$ in $E^*$ is
a smooth closed real curve $\gamma_i^j$, and that $\de \cS_i'$
meets $E^*$ transversally  along 
$\gamma_i = \bigcup_j \gamma_i^j$. Since $\cS_i'$ is a germ
and $\tau (z')=z$ we see that $\gamma_i \subset \cE^z$.
Let $T_i$ be a small neighborhood of $\gamma_i$ in $\de \cS_i'$.
By the excision theorem $T_i$ generates
the same element of $H_3(X'\setminus L, X'\setminus E)$ as
$\de \cS_i'$ does. Thus the isomorphism
$H_0 (L) \simeq H_3(X'\setminus L, X'\setminus E)$
sends the generator of $H_0 (L_i)$ to $T_i$. 
The advantage of $T_i$ is that it can be chosen as a fibration
over $\gamma_i$ where for each $j$ the fiber of $T_i$ over
the curve $\gamma_i^j$ can be viewed as a small complex
disc $\Delta_i^j$ which meets $E^*$ transversally
at a point of $\gamma_i^j$. Furthermore,
the connected component of $T_i$
over $\gamma_i^j$ is naturally diffeomorphic to
$\Delta_i^j \times \gamma_i^j$.
Recall that the Thom class $v_i^j$ of $\gamma_i^j$ in $X'\setminus L$
is an element of $H_2(X'\setminus L, X'\setminus E)$ uniquely
defined by the condition $v_i^j(\Delta_i^j)=1$.
The Thom isomorphism $H_3(X'\setminus L, X'\setminus E)
\to H_1 (E^*)$ is defined by the cap-product with
the Thom class $v= \sum_{i,j} v_i^j \, :
\theta \to v \cap \theta$. The properties of the
cap-product [Do, Ch. 7, 12.5 and 12.6] 
imply that the image of $T_i$ is $\gamma_i$
which is viewed as an element of $H_1 (E^*)\simeq H_1 (C^*_{m-1})$.
Note that $H_1 (C^*_{m-1}) = \oplus_{z \in C_{m-1} \setminus C^*_{m-1}}
H_1 (\cC^z) \oplus N$ where the group $N$ is not trivial
provided that either $C_{m-1}$ is of positive genus or
$C_{m-1}$ has more than one puncture.
Since $\gamma_i \subset \cE^z$ we see that the image
of the generator of $H_0(L_i)$ under the isomorphism
$H_0(L)\simeq H_1 (C^*_{m-1})$ is contained in $H_1(\cC^z)$.
Hence the image of $H_0(L)$ is contained in
$\oplus_{z \in C_{m-1} \setminus C^*_{m-1}} H_1 (\cC^z)$.
Thus $N$ is trivial and $H_0(L^z)\simeq H_1 (\cC^z)$.
This is the desired conclusion.
\qed

\medskip

The proof of Lemma 3.3 implies more.
Let $\cC_j, \, j=1, \dots , k$ be the irreducible
components of $\cC_z$. Then $\cC_j$ corresponds to
a generator $\alpha_j$ of $H_1 (\cC^z )$. If $L_i$
is an irreducible component of $\tau^{-1} (z)$ then
it corresponds to a generator $\beta_i$ of $H_0 (L^z)$.
By Lemma 3.3 the image of $\beta_i$ under the isomorphism
$H_0(L^z)\simeq H_1 (\cC^z)$ is $\sum_j m_i^j \alpha_j$.
One can extract from Lemma 3.3 the way to compute
these coefficients $m_i^j$.

\medskip

\lemma {\it Let the notation above hold and let
$\cS_i'$ be a germ of a holomorphic smooth surface which is
diffeomorphic to a complex two-dimensional ball
and transversal to $L_i$ and $E'$ at a generic point
$x'$ of $L_i$. Suppose $\cE_j$ is the closure of
$\tau^{-1}(\cC_j \setminus z)$ and $\cS_i'$ meets
$\cE_j$ along a curve $\G_i^j$. Then
the mapping $\tau |_{\G_i^j} : \G_i^j \to \cC_j$
is $m_i^j$-sheeted where $m_i^j$ is as before this Lemma.}

\proof
One can suppose that the boundary $\de \cS_i'$ of $\cS_i'$
meets $\cE_j$ transversally along a closed real curve $\gamma_i^j$.
It was shown in
the proof of Lemma 3.3 that the image of $\beta_i$
under the isomorphism $H_0(L^z)\simeq H_1 (\cE^z)$
is $\sum_j \gamma_i^j$ where $\gamma_i^j$ is viewed
as an element of $H_1 (\cE^z)$. Then the image of
$\gamma_i^j$ under the isomorphism $H_1(\cE^z)\simeq H_1 (\cC^z)$
coincides with $m_i^j \alpha_j$ where $m_i^j$ is the 
winding number of $\tau (\gamma_i^j)$ in $\cC_j$
around $z$. On the other hand $\gamma_i^j$ is the boundary
of $\G_i^j$. This implies that
$\tau |_{\G_i^j} : \G_i^j \to \cC_j$ is $m_i^j$-sheeted.
\qed

\medskip

\rem (1) If $E$ is a UFD then there is no need to assume in
Lemma 3.3 that $X'$ is smooth and $H_3(X')=H_2(X'\setminus E) =
H_3(X'\setminus E)=0$. Indeed, in this case
$\tau : E \to C_{m-1}$ generates a morphism
$\tau_{\nu} : E \to C_{m-1}^{\nu}$ where $C_{m-1}^{\nu}$
is a normalization of $C_{m-1}$. Each fiber $\tau_{\nu}^{-1}(z)$
consists of one connected component (since otherwise 
one can easily check that the defining
ideal of any of these connected components in $\tau_{\nu}^{-1}(z)$
is not principal, i.e. $E$ is not a UFD).
By Lemma 3.2 (3) $E$ is naturally homeomorphic to $\C \times C_{m-1}^{\nu}$
and the first claim of Lemma 3.3 holds automatically.
In order to have the second claim it is enough to require
that the Euler characteristics of $E$ is 1.
Then $C_{m-1}^{\nu}$ has Euler characteristics 1
whence it is isomorphic to $\C$.

(2) Furthermore, 
one may require  instead of that condition 
that $X'$ is smooth, $H_3 (X'\setminus E)=0$, 
and $E$ has at most isolated singularities
(and Euler characteristics 1).
Then that $E$ is a topological
manifold. The second homology group of this manifold is nontrivial
in the case when $\tau^{-1}(z)$ has more then one connected
component for some $z \in C_{m-1}$. On the other hand
the Thom isomorphism and the exact homology sequence of
the pair $(X', X' \setminus E)$ imply that this second homology groups
is isomorphic to $H_4(X')$ which is trivial.
Hence $\tau^{-1}(z)$ consists one connected
component for each $z \in C_{m-1}$ which yields again Lemma 3.3.

\bigskip

\subsection{The Reduced Center}
We shall describe some condition under which
the reduced center of an affine modification is
a complete intersection.

\medskip

\lemma {\it Let $C$ be a closed reduced irreducible curve in $\C^n$
and let $\bx = (x_1, \dots , x_n)$ be a coordinate system on $\C^n$.
Suppose that $D_1 = \C \times C $, that $v$ is a coordinate
on the first factor of $D_1$, and $\theta : D_1 \to C$ is
the natural projection. 

(1) Let $o$ be a singular point of $C$,
$\cV$ be the germ of $C$ at $o$,
and $\cH = \theta^{-1}(\cV )$.
Suppose that a function $h$ is holomorphic everywhere
in $\cH$ except for a finite number of points,
and that $h$ is a polynomial in $v$ over the ring
of functions on $\cV$.
Then $h$ is holomorphic in $\cH$.

(2) Let $h$ be a holomorphic function on $D_1$
whose zero set does not contain fibers of $\theta$
and let the zero multiplicity of $h$ at generic points
of this zero set be $n$. 
Suppose that $h$ is a polynomial in $v$ over the ring
of functions on $\cV$. Then $h^{1/n}$ is a holomorphic
function in $D_1$.

(3) Let $C_1$ be a closed reduced irreducible algebraic
curve in $D_1$ such that 

(i) the projection $\theta |_{C_1} : C_1 \to C$ is finite;

(ii) for every singular point $o \in C$
there exists $\cH$ as in (1)
such that in the ring of polynomials in $v$,
whose coefficients are holomorphic functions on $\cV$,
the defining ideal of $C_1\cap \cH$ is principal.

Then the defining ideal of
$C_1$ in $\C [D_1]$ is the principal ideal generated by
an irreducible regular function $b$ on $D_1$
and, if we treat $b$
as the restriction of a polynomial from $\C^{n+1}$ then
$$b (\bx ,v) =v^m+p_{m-1}(\bx ) v^{m-1} + \dots + p_0 (\bx )$$
where each $p_i \in \C [\bx ]$.}

\proof The argument is of local analytic nature and, therefore, it is
enough to consider the case when the normalization of $C$
is $\C$ (actually, we are only interested in the case when
$C$ is contractible).
Let $\nu_0 : \C \simeq C^{\nu } \to C$ be a normalization and
let $t$ be a coordinate on $C^{\nu }$. Then
$\nu = (\nu_0 , {\rm id }) : \C^2 \simeq \C \times C^{\nu} 
\to D_1$ is a normalization of $D_1$ and we can suppose that
$(v,t)$ is a coordinate system on this sample of $\C^2$.
For (1) put $\gamma = h \circ \nu$. This function is holomorphic
in a neighborhood of the $v$-axis in $\C^2$ by the
Riemann theorem about  deleting singularities, and it is of form
$$\gamma (v,t) = r_k(t)v^k+r_{k-1}(t) v^{k-1} + \dots + r_0 (t ).$$
The fact that $h$ is holomorphic everywhere on $\cH$
except for a finite number of points
implies that for every fixed $v=v_0$ except for
a finite number of values the function $\gamma (v_0,t)$
is contained in the ring generated by the coordinate
functions $x_1(t), \dots , x_n(t)$ of $\nu_0$.
This implies that each $r_i(t)$ belongs to
this ring whence $h$ is holomorphic in $\cH$ which is (1).

For (2) note that the function $h^{1/n}$ is holomorphic
everywhere in $\cH$ except for possibly points from
the finite set $h^{-1}(0) \cap \theta^{-1}(o)$.
Hence (1) implies (2).

Put $C_1^{\nu} = \nu^{-1} (C_1)$. It is the zero fiber of
an irreducible polynomial $\beta (v,t)$ on $\C^2$. Note that
the projection of $C_1^{\nu}$ to the $t$-axis
is finite since $\theta |_{C_1}$
is finite. Hence we can suppose that 
$\beta (v,t) =v^m+q_{m-1}(t ) v^{m-1} + \dots + q_0 (t )$,
i.e. $\beta$ is monic in $v$.
The function $b = \beta \circ \nu^{-1}$ is rational on $D_1$
and we are going to show that it is, in fact, regular.
It suffices to show that $b$ is holomorphic at each point
of $D_1$ (e.g., see [Ka2]).
Let $o$ be a singular
point of $C$. Since $\nu^{-1}$ is regular outside lines
of form $\theta^{-1}(o) \subset D_1$ it is enough to check that
$b$ is holomorphic at the points of $\theta^{-1}(o)$.

Let $O$ be the ring of germs of
analytic functions at the origin of $\C^n$ (whose coordinate system
is $\bx$).
Suppose that $h$ be the generator of the defining ideal of
$C_1 \cap \cH$ in the ring of holomorphic functions on $\cH$
that are polynomials in $v$.
By Cartan's theorems (e.g., see [GuRo, Ch. 8A, th. 18])
we can extend each coefficient of
$h$ (as a polynomial in $v$)
to a holomorphic function in a Stein neighborhood
of the origin in $\C^n$
%of $\cH$ in $\C^{n+1}$
%with coordinates $(\bx , v)$ 
whence we can treat $h$ as an
element of $O[v]$. 
Suppose that $o$ is the origin of $\C^n$.
Let $o_1, \dots , o_k$ be the set of all points from $C_1$ 
such that $\theta (o_i)=o$ for every $i$.
Let $c_i$ be the $v$-coordinate of $o_i$.
By the Weierstrass Preparation Theorem [Rem, Ch. 1, Th. 1.4]
there exists a unique Weierstrass polynomial
$\omega_1 \in O [v]$ such that
$h=\omega_1 (\bx , v-c_1)e_1$ where $e_1 \in O[v]$ does not
vanish at $o_1$. Applying this theorem again we see that
there exists a unique Weierstrass polynomial
$\omega_2 \in O [v]$ such that
$e_1=\omega_2 (\bx , v-c_2)e_2$ where $e_2 \in O[v]$ does not
vanish at $o_1$ and $o_2$. Hence
$h=\omega_1 (\bx , v-c_1)\omega_2 (\bx , v-c_2)e_2$. 
Repeating this process we get by induction that
$h=\omega e$ where $\omega \in O[v]$ is a monic polynomial,
whose zeros on the $v$-axis are $o_1, \dots , o_k$,
and $e\in O[v]$
does not vanish at $o_i$ for each $i$ (which implies that $e|_{\cH}$
is invertible).
Thus  $\omega |_{\cH}$ generates the same principal
ideal as $h$.
Therefore, we can suppose from the beginning that $h$
is a monic polynomial in $v$.

Hence $\gamma =h \circ \nu$ is monic as a polynomial
in $v$ over the ring of germs of
analytic functions at the finite set $\nu_0^{-1}(o) \in \C$.
Note that $\gamma = \beta \alpha$ where
$\alpha$ does not vanish since the zero
multiplicity of $\gamma$ and $\beta$ at generic points
$C_1^{\nu}$ is 1. Hence since $\alpha$ is a rational function
in $v$ it is constant on each line parallel to the $v$-axis.
Furthermore, this constant is 1 since both $\gamma$ and $\beta$
are monic (look at the quotient $\gamma /\beta$ as $v$
approaches $\infty$ along any of these lines).
Thus $\beta =\gamma$ whence $b$ coincides with $h$ in $\cH$
and, therefore, $b$ is holomorphic.
\qed

\medskip

\prop {\it Let the assumptions of Convention 2.2 and
Proposition 2.15 hold. Suppose that $\dim X =3$, $m \geq 2$ where
$m$ is from Proposition 2.15, and either

($\alpha$)
$X'$ is smooth, and $H_3 (X')=H_2(X' \setminus E)
=H_3(X'\setminus E)=0$, or

($\beta$) $E$ is a UFD and its Euler characteristics is 1.
%($\gamma$) $X'$ is smooth, $H_3(X'\setminus E)=0$, $E$ has at most isolated
%singularities and Euler characteristics 1.

Suppose also that
$D_1=E_1$ is isomorphic to $\C \times C$
(i.e. $C$ is a curve and it is
the reduced center not only of $\sigma$ but 
of $\sigma_1=\delta_1$ as well)
and $\theta : C_1 \to C$ is finite 
where $\theta =\sigma_1 |_{C_1}$.
Let $X$ be a local holomorphic UFD with respect to
$A\hookrightarrow A'$ (see Definition 2.10).

Then $C_1$ is closed in $D_1$ (i.e. $C_1$ is also the reduced center
of $\rho_1$), its defining ideal in $\C [D_1]$
is principal, and $C$ is contractible.}

\proof Put $\tau = \rho_{m-1} |_E$ and
suppose that 
$L^z$ be as in Lemma 3.3 and Remark 3.1
where $z \in C_{m-1} $. Let $L_i$ be one the
components of $L^z$. 
Consider a germ $\cS_i'$ of a smooth analytic surface 
transversal to both $E$ and $L_i$ at a generic point $z'$
of $L_i$,
and consider $\cC_j , \cE_j, \G_i^j$ as in Lemma 3.3.
That is, the germ $\cS_i' \cap E$ coincides with
$\bigcup_j \G_i^j$.
Let $\cS_i$ be the germ $\sigma (\cS_i')$
of a surface in $X_{m-1}$ at $z=\sigma (z')$ (indeed,
this is the germ of a surface by Lemma 2.2). 
Then $\cS_i \cap D_{m-1}$ consists of components $\cC_j$
as above. 

Let $m_i^j$ be as before Lemma 3.4.
Then the preimage of a generic point $x\in \cC_j$
under $\sigma |_{\cS_i'}$ consists
of $m_i^j$ points. Let $x'$ be one of these points.
Recall that $\cS_i'$ is smooth at $x'$
and transversal at $x'$ to the exceptional
divisor $E$. Furthermore, we can suppose that
$\cS_i'$ is transversal to the fiber of $\sigma$ 
through $x'$. Since $\sigma$ can be viewed
as a usual monoidal transformation in a neighborhood of 
a generic point $x\in C_{m-1}$ (by Remark 2.8)
the last fact implies that the germ of $\cS_i$ at $x$ consists
of $m_i^j$ irreducible smooth components each of which meets $D_{m-1}$
transversally. 

Let $\gamma = \sigma_{m-1} \circ \cdots \circ \sigma_2$.
Put $\cS_i^1 =\gamma (\cS_i)$ and let $\cG_l$ be 
the image of $\cC_j$ under $\gamma$
(that is, $\cG_l$ is an irreducible component of the
germ of $C_1$ at $z_1=\gamma (z)=\rho_1 (z')$). Suppose that
the mapping $\gamma |_{\cC_j} : \cC_j \to \cG_l$ is $n_j$-sheeted.
Let $R_l$ be the set of natural $j$'s such that
$\gamma (\cC_j) = \cG_l$.  Since $\gamma$ is a composition
a basic modifications which are nothing but usual monoidal
transformations over generic point $x_1=\gamma (x)$ of
$\cG_l$ one can see that the germ of $\cS_i^1$ at $x_1$
consists of $\sum_{j \in R_l} n_jm_i^j$ smooth 
irreducible components which meet $D_1$ transversally at
$x_1$. By the assumption of this Proposition
$X$ is a local holomorphic UFD with respect to
$A\hookrightarrow A'$ whence Proposition 2.14 implies
the defining ideal of $\cS_i^1$
is generated by the germ $h_i$ of a holomorphic function.
Since the zero multiplicity of $h_i$ at smooth points of $\cS_i^1$
is 1 we see that the zero multiplicity of
$h_i|_{D_1}$ at $x_1$ is $\sum_{j \in R_l} n_jm_i^j$.

Let $L^z$ consists of $k$ components $L_1, \dots , L_k$.
Since we have an isomorphism $H_0 (L^z) \simeq H_1 (\cC^z )$
by Lemma 3.3, $j$ changes from 1 to $k$ and the matrix
$(m_i^j)$ is invertible. Hence there exist integers
$s_1, \dots , s_k$ such that the germ of the function
$h= (h_1|_{D^0})^{s_1} \cdots (h_k|_{D^0})^{s_k}$ has
the zero multiplicity $n_1$ at generic points of
$\cG_1$ and the zero multiplicity 0 at generic points of
the other irreducible component of the
germ of $C_1$ at $z_1$ (i.e. $h$ is different from zero
or $\infty$ at generic points of these other components).
Let $\cD_1$ be the germ of $D_1$ at $z_1$ and let $\cD^1,
\cD^2 , \ldots $ be its irreducible
components. Since the restriction of $h$ to $\cD^i$
is well-defined everywhere on $\cD^i$ except may be for $z_1$
we see that $h|_{\cD^i}$ becomes homomorphic function
after normalization in virtue of the Riemann theorem
about deleting singularities. In particular, it is continuous at $z_1$.
Furthermore, if $\cG_1 \subset \cD^1$ then
$h|_{\cD^1}$ vanishes at $z_1$. Since $\cD^1 \cap \cD^2$
is a curve which contains $z_1$ we see that
$h|_{\cD^2}$ vanishes also at $z_1$ whence
the set of zeros of $h|_{\cD^2}$ is the germ of a curve,
and this germ can be only $\cG_1$. But it cannot be true
since then $\cG_1$ is not contained
in the set of double points of $D_1$ (recall that
the morphism $\theta : C_1 \to C$ is not constant).
Thus $\cD_1$ is irreducible
whence the geometrical center $C_0$ of $\sigma$
has no double points. Since $C_0$ admits
a non-constant morphism from $C_{m-1}$ and by Lemma 3.3 the normalization of
$C_{m-1}$ is $\C$ this implies that $C_0$ is contractible
and coincides with the reduced center $C$ of $\sigma$.

Since $C_1$ admits a non-constant morphism from $C_{m-1}$
it is a once punctured curve and, therefore, it is
closed in $D_1$. Note that $h$ is holomorphic everywhere
on $\cD_1$ except, may be, $z_1$. By Lemma 3.5 $h$ is
holomorphic on $\cD_1$. Consider the function $e$ on $\cD_1$
such that $e^{n_1} =h$. 
%It is easy to see again that
%$e$ is holomorphic everywhere on $\cD_1$ except, may be, $z_1$.
By Lemma 3.5 it is holomorphic. 
Recall that the domain of $h_1$ includes the strip $\sigma_1^{-1}(\cV )
\simeq \C \times \cV$ where $\cV$ is the germ of $C$
at $\sigma (z')$ (see Remark 2.9).
Hence $e$ is defined
in this strip. 
Let $v$ be a coordinate on the first factor of the strip.
For every holomorphic
function $e_1$ on this strip, which is polynomial in $v$
and which vanishes on $\cG_1$,
the quotient $e_1/e$ is again holomorphic by Lemma 3.5.
Therefore, $e$ is the generator of the defining ideal
of $\cG_1$. Hence the defining ideal of the germ of $C_1$
at $z_1$ is principal.
Note that $\theta : C_1 \to C$ is a finite morphism
since $C_1$ is a once puncture curve. By Lemma 3.5
the defining ideal of $C_1$ in the ring of regular
functions on $D_1$ is principal.
\qed

\medskip

This concludes the proof of Theorem 2.3.

\bigskip

\section{Applications of the Decomposition}

\subsection{The proof of Miyanishi's theorem}

We shall reduce first the Miyanishi theorem to a problem
about affine modifications.

\medskip

\lemma {\it
Let $X'$ be an affine algebraic variety of dimension 3
such that
$X'$ is a UFD and there exists
a Zariski open subset $Z$ of $X'$
which is a cylinder over a smooth curve $U$. 

(i) Then $U$ is a
Zariski open subset of $\C$ and the natural projection
$p_0 : Z \to U$ can be extended to a regular function $p :X' \to \C$
whose general fibers are still isomorphic to $\C^2$.

(ii) Furthermore, let $x,y,z$ be coordinates on $X = \C^3$.
Then there exists an affine modification 
$\sigma : X'\to X$ such that its coordinate form is
$\sigma = (p,p_1,p_2)$ and the divisor
of this modification coincides with the zeros of some
polynomial $f(x)$ on $\C^3$.}

\proof
Let $\bF_c$ be the closure of the fiber $F_c =\{ p_0 =c \} \subset Z$
in $X'$ (where $c \in U$).
Assume that $\bF_c \cap \bF_{c'} \ne \emptyset$ for some
$c \ne c' \in U$.  Since $X'$ is a UFD there exists
a regular function $g$ on $X'$ whose zero set
coincides with the divisor $\bF_{c'}$.
Thus the zero locus of $g|_{\bF_c }$ is $\bF_c \cap \bF_{c'} $. 
On the other hand $g|_{\bF_c }$ is nowhere zero on
$\bF_c \setminus (\bF_c \cap \bF_{c'} ) \supset F_c \simeq \C^2$
whence this function must be a nonzero constant on $F_c$ and, therefore,
$\bF_c$. Contradiction.
Thus $\bF_c \cap \bF_{c'} = \emptyset$ for every
$c' \ne c \in U$.

Assume $F_c \simeq \C^2$ if different from $\bF_c$.
Assume that one of the irreducible components of $\bF_c \setminus F_c$
is a point. Then a normalization $G$ of $\bF_c$ contains a sample
of $\C^2$ and one of the irreducible components of $G\setminus \C^2$
is also a point $o$. By the theorem about deleting singularities
of holomorphic functions 
in codimension 2 for normal complex spaces [Rem, Ch. 13] every
holomorphic function on $\C^2$ can be extended to this point $o$
whence $\C^2$ is not Stein. Contradiction.
Thus $\bF_c \setminus F_c$ is a curve.
Since $\bF_c \cap \bF_{c'} = \emptyset$ we see that the closure of
$\bigcup_{c \in U} (\bF_c \setminus F_c)$ is a divisor in $X'$.
Since $X'$ is a UFD there exists a regular function $h$ on $X'$
whose zero set coincides with this divisor.
Thus the zero locus of $h|_{\bF_c }$ is $\bF_c \setminus \bF_{c} $
and we get a contradiction in the same way we did for function $g$.
Hence $F_c = \bF_c$.

This implies that $p_0 : Z \to U$ can be extended to continuous
map $p$ from $X'$ to the completion $\bU$ of $U$,
and $p^{-1}(U)=Z$. In particular, general fibers of $p$
are isomorphic to $\C^2$. Since $X'$ is a UFD $p$
must be holomorphic [Rem, Ch. 13] and, therefore, regular
(e.g., see [Ka2]).

Since $X'$ is a UFD (i.e. every effective
divisor is the zero divisor of some regular function on $X'$)
we see that the Zariski open subset $Z$ of $ X'$ is also a UFD whence
$U$ is a
UFD. This implies that $U$ is rational, i.e. 
its completion is $\bU = \pr^1$.

Show that $p : X' \to \pr^1$ is not surjective.
Assume the contrary. Let $X_0=p^{-1}(\C )$ and
$q=p |_{X_0}$.
We can suppose that $Z \subset X_0$, i.e. $U \subset \C$.
Extend the isomorphism $Z \simeq U\times \C^2 \subset \C
\times \C^2$ to a
rational map from $X_0$ to $\C^3$ (with 
coordinate $x,y,z$) and then multiply
the two last coordinates by polynomials in $q$ to make
this mapping regular. We obtain a birational morphism
$\sigma : X_0 \to \C^3$. This is an affine modification by Theorem 2.1.
It is clear that $q =x \circ \sigma$ and the divisor
of this modification is given by the zeros of some polynomial
$f(x)$ in $x$. Since $q: X_0 \to \C$ is surjective 
Proposition 2.5 implies that every invertible function
on $X_0$ is of form $h\circ \sigma$ where $h$ is invertible
on $\C^3$. Therefore, each invertible function on $X_0$ is constant.
On the other hand the divisor $p^{-1}(\infty )$ is the
zero divisor of some function $g$ on $X'$ since $X'$ is a UFD.
Hence $g|_{X_0}$ is invertible and non-constant.
Contradiction. Thus one can suppose that
$p =q$ and $X'=X_0$ which concludes the proof.
\qed

\medskip

Now we shall consider the case when the polynomial $f(x)$
constructed in the previous Lemma coincides with $x^n$.

\medskip

\lemma {\it
Let $A$ be the polynomial ring $\C [x,y,z]$ in three variables
$x,y$, and $z$, let $f=x^n$, and let $A\hookrightarrow A'$
be an affine modification. Suppose that $q(x)$ is
a polynomial in $x$ such that $q(0) \ne 0$. 
Let $B=A[1/q], J=I[1/q]$, and let $B\hookrightarrow B'$ be
an affine modification with locus $(J,f)$. Suppose that
$B'$ is a UFD and one of the following conditions hold

($\alpha$) $E$ is non-empty, $X'$ is smooth, and $H_3 (X')=0$;

($\beta$) $E$ is a UFD and
and its Euler characteristics is $e(E)=1$.

Let $E$ have at most isolated singularities.
Then $A'$
is also a polynomial ring $\C [x,u,v]$ in variables $x,u$, and $v$.}

\proof
Let $\delta : Y' \to Y= \C^3 \setminus \{ q^{-1}(0) \}$ be
the affine modification of the varieties
which corresponds to the modification $B \hookrightarrow B'$
with locus $(J,f)$.
Since affine modifications commute with localizations
in the sense of Proposition 2.1 we see that $E$
is also the exceptional divisor of this modification $\delta$.
If $B' \ne B$ we can suppose that the reduced center of
this modification is at least of codimension 2 in $Y$.
(Indeed, otherwise each element $h \in J$ must vanish on
the plane $\{ x=0 \} \subset \C^3$ whence $h=xh_1$.
We can replace the locus $(J,x^n)$ by the locus $(J/x, x^{n-1})$.
After several replacement like this one we shall obtain
an element of $J$ which does not vanish on the plane.)

Show that Theorem 2.3 is applicable which is clear if
condition ($\beta$) holds.
Let $T=\C^3 \setminus \{ xq(x)=0 \}$
and $Z= \C^3 \setminus \{ x=0 \}$. Note that 
$H_i (Z)=H_i(T)=0$ for $i \geq 2$ and
$T$ is isomorphic to $Y'\setminus E$.
Hence $H_2(Y'\setminus E) =H_3(Y'\setminus E)=0$.
We can glue $Y'$ and $Z$
along $T \simeq Y'\setminus E$. Then we have
$X' = Y' \bigcup_T Z$.
The Mayer-Viertoris theorem implies
that $H_3(X')=H_3(Y')=0$. Hence condition
($\alpha$) in this Lemma implies condition ($\alpha$)
in Theorem 2.3.

By Theorem 2.3 we have now a sequence of basic
modifications
$$B =B_0 \hookrightarrow B_1 \hookrightarrow \ldots
\hookrightarrow B_m =B' $$
which corresponds to the sequence of morphisms
$$Y'= Y_m \stackrel{\delta_m}{\rightarrow}
 Y_{m-1} \to  \ldots \to Y_1 \stackrel{\delta_1}{\rightarrow} Y= \C^3
\setminus \{ q^{-1}(0) \}.$$
The natural embeddings 
$B_i \hookrightarrow B'$ are affine modifications which generate morphisms
$\theta_i : Y' \to Y_{i}$.
Let $C_{i}=\theta_i (E)$.
By Theorem 2.3 each $C_{i}$ is either a point or 
a smooth contractible irreducible curve
which is a connected component of the geometrical
center of $\delta_i$.

Our first aim is to show that
$B'$ is a localization of the
polynomial ring $\C [x,u,v]$ with respect to 
to the multiplicative system $\{q^n(x) | \, n \in \N \}$
where $u,v \in A'$.
Suppose that $C_0=C$ is a point (say, the origin 
$o= \{ x=y=z=0 \}$). Let $M$ be the maximal ideal in $B$ that
vanishes at $o$. By Theorem 2.3 $B_1 = B[M/x]$.
Hence $B_1 = A_1[1/q]$ where $A_1$ 
is the polynomial ring $\C [x, y/x,z/x]$ in three variables.
Suppose that $j$ is the first number for which
$C_j$ and, therefore, every $C_k$ with $k>j$
are curves (recall that $\delta_i |_{C_i} : C_i \to C_{i-1}$
must be surjective by Proposition 2.15).
By induction we can suppose that $B_j =A_j[1/q]$ where 
$A_j$ is a polynomial ring $\C [x, \xi , \zeta ]$
(in particular, the divisor $D_j$ 
of $\theta_j$ and $\delta_{j+1}$ is the $\xi \zeta$-plane).

By Theorem 2.3 $C_j$ is isomorphic to $\C$ and
by the Abhyankar-Moh-Suzuki theorem one can assume that
it is given by $x=\xi =0$. 
Let $I_{j}$ be the ideal generated by $x$ and $\xi $.
By Theorem 2.3 
$B_{j+1}=B_j[I_{j}/x]$. Then by Proposition 2.1
$B_{j+1} = A_{j+1}[1/q]$ where $A_{j+1}$
is the polynomial ring $\C [x, \xi /x,\zeta ]$ in three variables.
Hence our first aim can be achieved by induction.

By Proposition 2.11
the projection $Y' \to Y$ generates isomorphisms
of the homology groups and the fundamental groups.
As we mentioned in the beginning of the proof
$X' = Y' \bigcup_T Z$ and, similarly, $\C^3 = Y\bigcup_T Z$.
The Mayer-Viertoris and van Kampen-Seifert theorems imply
that the homology groups and the fundamental group
of $X'$ are trivial since this is true for $\C^3$.
Hence $X'$ is contractible.
Note also that $X'$ is smooth since $Y'$ is smooth.
This implies that $X'$ is a UFD (e.g., see [Ka1]).
Thus if we put $q(x)$ equal to 1 then the assumptions
of this Lemma are true (since we know now that $A'$
is a UFD). 
Therefore, $B'=A'[1/q]=A' =\C [x,u,v]$.
\qed

\bigskip

\lemma {\it
Let $A =\C [x,y,z]$, 
let $f$ is a polynomial in one variable $x$,
and let $A \hookrightarrow A'$ be an affine modification.
Suppose that $A'$ is a UFD and that for every root $c$ of
$f$ the polynomial $f-c$ is not a unit in $A'$
(or, equivalently the mapping $f \circ \sigma : X' \to \C$
is surjective). Let one of the following conditions hold

($\alpha$) $X'$ is smooth and $H_3 (X')=0$;

($\beta$) $e(X')=1$ and every irreducible component
of $E$ is a UFD.

%(i) Then the Euler characteristic of every component
%of the exceptional divisor of 
%this modification $\sigma : X' \to X \simeq \C^3$ is 1.
If every irreducible component of $E$
has at most isolated singularities (which is
automatically true under condition ($\beta$)) then
$X'$ is isomorphic to $\C^3$ and $x \circ \sigma$ is a variable
on this sample of $\C^3$.} 

\proof
Let $f(x)=x^nq(x) $ where $q(0) \ne 0$,
$J=I[1/q]$ and $B=A[1/q]$.
By Proposition 2.1 $B'=B[J/x^n]$ coincides with $A'[1/q]$.
Hence $B$ and $B'$ are UFDs and 
the exceptional divisor $E^0$
of the modification
$B \hookrightarrow B'$ is not empty by the assumption on $f\circ \sigma$.
It is irreducible by Proposition 2.4.
This makes Lemma 4.2 applicable to this modification
under condition ($\alpha$).
Show that the same is true under ($\beta$). By 
Proposition 2.15 we can present $B \hookrightarrow B'$
as a composition of basic modifications.
Hence $E^0$ is homeomorphic to $C^0 \times \C^s$
where $C^0$ is a point (see Lemma 2.4 and Remark 2.10)
or an irreducible curve (see Lemma 3.2 and Remark 3.1)
and $s=2$ or
1 respectively.
That is, $e(E_0) \leq 1$ in any case. Let
$D_0$ be the coordinate plane $x=0$. 
Then by the additivity of Euler characteristics
[Du] $e(X')$ differs from $e(X)=e(\C^3)=1$ by the sum
of terms of form $e(E_0)-e(D_0)$ (these terms should
be considered for each root of $f$).
Since $e(X')=1$ we see that $e(E_0)=e(D_0)=1$ which 
makes Lemma 4.2 applicable.

Suppose that $L$ is the ideal in $A$ generated by $I$ and $x^n$,
i.e. $I[1/q]=L[1/q]=J$.
By Lemma 4.2 $B'=A^1 [1/q]$ 
where $A^1=A[L/x^n]= \C [x,u,v]$ is a polynomial
ring in three variables. 
Let $K$ be the ideal in $A^1$ generated by $I/x^n$.
By Theorem 2.1 (3) $A'=A^1[K/q]$.
Now the induction by
the degree of $f$ implies implies the desired conclusion.
\qed

\bigskip

\rem In fact the assumption that $A'$ is a UFD
can be replaced by a weaker one. Namely, one can
assume only that for every root $c$ of $f$ there
exists a polynomial $r(x)$ with $r(c) \ne 0$ such that
the localization $A[1/r]$ is a UFD.

\medskip

Lemmas 4.1 and 4.3 imply

{\bf  Miyanishi's Theorem (Lemma I).} {\it
Let $X'$ be an affine algebraic variety of dimension 3
such that $X'$ is a UFD, all invertible functions
on $X'$ are constants, and

(1) the Euler characteristics of $X'$ is $e(X')=1$;

($2'$) there exists a Zariski open subset $Z$ of $X'$
which is a $\C^2$-cylinder over a curve $U$ (i.e. $Z$ is isomorphic
to the $\C^2 \times U$);

($3$) each irreducible component of $X' \setminus Z$ is a UFD.

Then $X'$ is isomorphic to $\C^3$.

(4) Furthermore, the curve $U$ is a Zariski open subset of $\C$,
the natural projection from $Z$ to $U$ can be extended to
a regular function on $X'$, and this function is a variable.

The statement of this Lemma remains true if conditions (1) and (3)
are replaced by

($1'$) $X'$ is smooth and $H_3(X')=0$;

($3'$) each irreducible component of $X' \setminus Z$
has at most isolated singularities.}

\bigskip

\subsection{How to present $X'$ as a closed algebraic
subvariety of $\C^N$}

In this section we shall study $X'$ which satisfies assumptions
($1'$) and ($2'$) of Lemma I, but we do not
require ($3'$). We want to present $X'$ explicitly
as a closed affine algebraic subvariety of some Euclidean
space $\C^N$. 

\medskip

\lemma {\it
Let $A=\C [x,y,z]$ and $f(x)=x^n$. Let
$A \hookrightarrow A'$ be an affine modification.
Suppose that $X'$ is a smooth UFD, the only invertible functions
on $X'$ are constants, and $H_3(X')=0$.
Then
 
(1) either $X'$ is isomorphic to $\C^3$ or at least

(2) $X'$ can be viewed as the subvariety of $\C^{3+m}$
given by polynomial equations
\[ \begin{array}{l}
xv_1-q_0(y,z )=0\\
xv_2-v_1^{n_1} + q_1(y,z,v_1 )=0\\
\, \, \, \, \, \, \, \, \, \, \, \, \, \, \, \, \, \, \ldots \\
xv_m-v_{m-1}^{n_{m-1}} + q_{m-1}(y,z,v_1, \dots , v_{m-1} )=0\\
\end{array} \]
where the usual degree of $q_j$ with respect to $v_i$
is less than $n_i$ for every $i=1, \dots , j$.
Furthermore, one can suppose that
$q_0(y,z)=y^k-z^l$ where $(k,l)=1, k>l\geq 2$, and
$m \geq 2$.}

\proof
Since $X'\setminus E \simeq \C^3 \setminus \{ x=0 \}$
we have $H_2(X'\setminus E)=H_3(X'\setminus E)=0$.
By Theorem 2.3 the modification $\sigma : X' \to X$
is a composition of basic modifications
$$X'= X_m \stackrel{\sigma_m}{\rightarrow}
X_{m-1} \to  \ldots \to X_1 \stackrel{\sigma_1}{\rightarrow} X.$$
Let $A_j =\C [X_j]$ and $C_j$ be as in Convention 2.2 (3).
%By Remark 2.4 $A' =A_j[K_j/x^n]$
%where $K_j$ is the $x^n$-largest ideal of the modification
%$A_j \hookrightarrow A'$.
%Suppose that $C_j =\cV_{X_j} (K_j)$.
By Theorem 2.3 each $C_j$ is either a point
or an irreducible contractible curve. If $C_0$ is
a point then, as it was shown in the proof of Lemma 4.2,
$A_1 \simeq \C [x, y/x,z/x]$ whence $X_1 \simeq \C^3$.
Therefore, we can suppose that $C_0$ is a curve
whence each $C_i$ is a curve since
$\sigma |_{C_i} : C_i \to C_{i-1}$ is surjective by Proposition 2.15.
If $C_0$ is a smooth curve, i.e. $C_0 \simeq \C$ then
it was shown in the proof of Lemma 4.2 that we can suppose that
$A_1 \simeq \C [x, y/x, z]$, i.e. $X_1 \simeq \C^3$.
Therefore, we consider the case when $C_0$ is not smooth.
By the Lin-Zaidenberg theorem [LiZa] one can assume
that the equations of this curve in $\C^3$ are
$x=y^k-z^l=0$ where $(k,l)=1$ and $k>l\geq 2$.
Let $I_1$ be the ideal in $A$ generated by $x$
and $y^k -z^l$. By Theorem 2.3
$A_1=A[I_1/x]$.
By Theorem 2.2 this implies
that $A_1=\C [x,y,z, (y^k-z^l)/x]$ and $X_1$
is the hypersurface in $\C^4$ with coordinates
$(x,y,z,v_1)$ given by
$$xv_1 = q_0(y,z):= y^k-z^l.$$
The exceptional divisor $E_1$ is $\sigma_1$ is the
intersection of this hypersurface with the hyperplane $x=0$.
By Theorem 2.3 and by Lemma 3.5
$C_1$ is the zero fiber of a regular function on $E_1$
which is  of form $v_1^{n_1} + q_1 (y,z,v_1)$ where the usual
degree of $q_1$ with respect to $v_1$ is at most $n_1 -1$.
Let $I_2$ be the ideal in $A_1$ generated by $x$ and 
$v_1^{n_1} + q_1 (y,z,v_1)$.
By Theorem 2.3 $A_2=A_1[I_2/x]$.
Therefore, by Theorem 2.2 $X_2$ may be viewed as the subvariety
of $\C^5$ (with coordinates $(x,y,z,v_1,v_2)$)
given by the equations
\[ \begin{array}{l}
xv_1-q_0(y,z )=0\\
xv_2-v_1^{n_1} + q_1(y,z,v_1 )=0.
\end{array} \]
Repeating the above argument we see that
$X'$ can be viewed as the subvariety of $\C^{3+m}$
given by the equations
\[ \begin{array}{l}
xv_1-q_0(y,z )=0\\
xv_2-v_1^{n_1} + q_1(y,z,v_1 )=0\\
\, \, \, \, \, \, \, \, \, \, \, \, \, \, \, \, \, \, \ldots \\
xv_m-v_{m-1}^{n_{m-1}} + q_{m-1}(y,z,v_1, \dots , v_{m-1} )=0\\
\end{array} \]
where the usual degree of $q_j$ with respect to $v_i$
is less than $n_i$ for every $i=1, \dots , j$.

In order to check that $m>1$ when $X'$ is smooth it is enough to note
that $X_1$ is singular at the origin.
\qed

\bigskip

\prop {\it
Let $A=\C [x,y,z]$ and $f$ be a polynomial in $x$.
Suppose that $A \hookrightarrow A'$ is an affine modification.
Let $c_0, c_1, \dots$ be the roots of $f$.
Suppose that $X'$ is a smooth UFD, the only invertible functions on $X'$
are constants, and $H_3(X')=0$.
Then
 
(1) either $X'$ is isomorphic to $\C^3$ or at least
there exists a root of $f$ (say $c_0$ and assume that $c_0=0$)
such that

(2) $X'$ can be viewed as the subvariety of $\C^{N}$
given by a system of polynomial equations
\[ \begin{array}{l}
xv_1-q_0(y,z )=0\\
xv_2-v_1^{n_1} + q_1(y,z,v_1 )=0\\
\, \, \, \, \, \, \, \, \, \, \, \, \, \, \, \, \, \, \ldots \\
xv_m-v_{m-1}^{n_{m-1}} + q_{m-1}(y,z,v_1, \dots , v_{m-1} )=0\\
(x-c_1)u_{1,1}-r_{1,0}(y,z )=0\\
(x-c_1)u_{1,2}-u_{1,1}^{n_{1,1}} + r_{1,1}(y,z,u_{1,1} )=0\\
\, \, \, \, \, \, \, \, \, \, \, \, \, \, \, \, \, \, \ldots \\
(x-c_1)u_{1,m_1}-u_{1,m_1-1}^{n_{1,m_1-1}} + 
r_{1,m_1-1}(y,z,u_{1,1}, \dots , u_{1,m_1-1} )=0\\
(x-c_2)u_{2,1}-r_{2,0}(y,z )=0\\
\, \, \, \, \, \, \, \, \, \, \, \, \, \, \, \, \, \, \ldots \\
\end{array} \]
where $q_0(y,z)=y^k-z^l, \, (k,l)=1, k>l\geq 2$ and $m>1$. Furthermore,

-the usual degree of $q_j$ with respect to $v_i$
is less than $n_i$ for every $i=1, \dots , j$ and;

- $r_{s,j}$ are polynomials such that
the usual degree of $r_{s,j}$ with respect to $u_{s,i}$
is less than $n_{s,i}$ for every $i=1, \dots , j$.

The variety $X'$ is irreducible, it is a complete intersection in $\C^N$,
i.e. the ideal $I'$ of polynomials that vanish on $X'$
is generated by the left-hand sides of this system.}

\proof
We consider the case when $f$ has two roots 0 and 1
(i.e. $f(x)=x^n(x-1)^s$) since
the general case differs only by more complicated notation.
%Let $X_1 =\C^3 \setminus \{ x=0 \} , X_1' = X' \setminus \{ x=0 \}$,
%$X_2 =\C^3 \setminus \{ x=1 \} , X_2' = X' \setminus \{ x=1 \}$,
%and $X_0 =\C^3 \setminus \{ f(x)=0 \} \simeq X' \setminus \{ f(x)=0 \}$.
%Note that $X_0$ can be treated as a subvariety of $X_1,X_2, X_1$,
%or $X_2'$. Thus we can glue some of these varieties along $X_0$.
%The results are $X^1 =X_1 \bigcup_{X_0}X_2', X^2=X_1' \bigcup_{X_0}X_2$.
%By Proposition 2.3 $\C [X^1]=A[I_1/x^n]$ where $I_1$
Consider $X^1={\rm spec} \, A[I_1/x^n]$ where $I_1$
is generated by $I$ and $x^n$. By Lemma 4.4 $X^1$
coincides with the zero set of the system
\[ \begin{array}{l}
xv_1-q_0(y,z )=0\\
xv_2-v_1^{n_1} + q_1(y,z,v_1 )=0\\
\, \, \, \, \, \, \, \, \, \, \, \, \, \, \, \, \, \, \ldots \\
xv_m-v_{m-1}^{n_{m-1}} + q_{m-1}(y,z,v_1, \dots , v_{m-1} )=0\\
\end{array} \]
where the usual degree of $q_j$ with respect to $v_i$
is less than $n_i$ for every $i=1, \dots , j$ and $m>1$.
Let $X^2:= {\rm spec} \, A[I_2/(x-1)^s]$ where $I_2$
is generated by $I$ and $(x-1)^s$.
By Lemma 4.4 $X^2$ can be given by the zeros of the system
\[ \begin{array}{l}
(x-1)u_{1}-r_{0}(y,z )=0\\
(x-1)u_{2}-u_{1}^{1,n_{1}} + r_{1}(y,z,u_{1} )=0\\
\, \, \, \, \, \, \, \, \, \, \, \, \, \, \, \, \, \, \ldots \\
(x-1)u_{m_1}-u_{m_1-1}^{1,n_{m_1-1}} + 
r_{m_1-1}(y,z,u_{1}, \dots , u_{m_1-1} )=0\\
\end{array} \]
where $r_{j}$ are polynomials such that
the usual degree of $r_{j}$ with respect to $u_{i}$
is less than $n_{i}$ for every $i=1, \dots , j$.
By Remark 2.6 we see that $X'$ is isomorphic to
the common zeros of the system
\[ \begin{array}{l}
xv_1-q_0(y,z )=0\\
xv_2-v_1^{n_1} + q_1(y,z,v_1 )=0\\
\, \, \, \, \, \, \, \, \, \, \, \, \, \, \, \, \, \, \ldots \\
xv_m-v_{m-1}^{n_{m-1}} + q_{m-1}(y,z,v_1, \dots , v_{m-1} )=0\\
(x-1)u_{1}-r_{0}(y,z )=0\\
(x-1)u_{2}-u_{1}^{1,n_{1}} + r_{1}(y,z,u_{1} )=0\\
\, \, \, \, \, \, \, \, \, \, \, \, \, \, \, \, \, \, \ldots \\
(x-1)u_{m_1}-u_{m_1-1}^{1,n_{m_1-1}} + 
r_{m_1-1}(y,z,u_{1}, \dots , u_{m_1-1} )=0.\\
\end{array} \]

In order to prove the last statement
consider $X_1$ given in $\C^4$ by the equation
$$xv_1-q_0(y,z )=0.$$
Clearly, $X_1$ is a complete intersection and it is irreducible.
Consider $X_2$ given in $\C^5$ by the equations
\[ \begin{array}{l}
xv_1-q_0(y,z )=0\\
xv_2-v_1^{n_1}+q_1(y,z,v_1)=0.
\end{array} \]
Let $A_i=\C [X_i]$. Note that $A_2=A_1[v_1^{n_1}/x]$,
is a Davis modification. In particular, $X_2$ is irreducible.
By Theorem 2.2
the ideal of polynomials that vanishes
on $X_2$ is generated by the left-hand sides of these
two equations. Now the induction implies the similar
conclusion about $X'$.
\qed

\medskip

Suppose that $X'$ satisfies the assumptions ($1'$) and ($2'$) of Lemma I.
By Lemma 4.1 this $X'$ satisfies the assumption of Proposition 4.1
whence we have

\medskip

\cor {\it Suppose that $X'$ satisfies the assumptions ($1'$) and ($2'$) of
Lemma I. Then $X'$ satisfies also either (1) or (2) from
Proposition 4.1.}

\bigskip

\section{The Makar-Limanov Invariant}

\subsection{General Facts about Locally Nilpotent Derivations}
Recall the following

\defin 
A derivation $\de$ on $A$
is a linear endomorphism which satisfies the
Leibniz rule, i.e. $\de (ab)=a\de (b) + b \de (a)$.
Two derivations are called
{\it equivalent} if they  have the same kernel. 
A derivation $\de$
is called locally nilpotent if for
each $a \in A$ there exists an $k = k(a)$ such that $\de^k(a) = 0$.

\medskip

Every locally nilpotent derivation defines
a degree function $\dg_{\de }$ on the domain $A$ 
with natural values (e.g., see [FLN])
given by the formula $\dg_{\de } (a) = \max \{ k \, | \,
\de^k (a) \neq 0 \}$ for every nonzero $a \in A$. 
The first five statements of the following theorem
can be found in [M-L1], [M-L2], [KaM-L1], and [KaM-L2].

\bigskip

\thm {\it 
Let $\de$ be a nonzero locally nilpotent 
derivation of $A$.

(1) $A$ has transcendence degree one over $\Ker {\de}$. 
The field ${\rm Frac}(A)$ of fractions of $A$ 
is a purely transcendental extension of ${\rm Frac}(\Ker {\de })$,
and $\Ker {\de}$ is algebraically closed in $A$.

(2) Let $b \in A$ and $\dg_{\de }(b)=1$.
Then for every $a \in A$ such that $\dg_{\de }(a)=k$
there exist $a',a_0, a_1, \ldots , a_k \in \Ker {\de }$
for which $a', a_k \neq 0$ and $a'a = \sum^k_{j=0} a_j b^j$.

(3) Every two locally nilpotent derivations
$\de$ and $\delta$ on $A$ are equivalent iff they
generate the same degree function.
Furthermore, there exist
$\alpha , \beta \in A^{\de }$ such that $\alpha \de = \beta \delta$.   

(4) Suppose $a_1, a_2 \in A$. Then
$a_1a_2 \in \Ker {\de }\setminus \{ 0 \}$ implies
$a_1, a_2 \in \Ker {\de }$.
In particular, every unit $u \in A$ belongs to $\Ker {\de}$.

(5) Suppose that $a_1^k+a_2^l\in \Ker \de \setminus \{ 0 \}$ where $k,l\geq 2$
are relatively prime. Then $a_1, a_2 \in \Ker \de$.

(6) {\rm (cf. [Za], proof of Lemma 9.3)}
Let $A = \C [X]$ and let $F=(f_1, \dots , f_s): X \to Y\subset \C^s$
and $G: Y \to Z\subset \C^j$ be dominant morphisms of reduced affine
algebraic varieties. Put $H= G \circ F=(h_1, \dots , h_j) : X \to Z$.
Suppose that for generic point $\xi \in Z$ there exists
a (Zariski) dense subset $T_{\xi}$ of $G^{-1}(\xi )$ such that
the image of any non-constant morphism from $\C$ to
$G^{-1}(\xi )$ does not meet $T_{\xi }$.
Then the fact that $h_1, \dots , h_j \in A^{\de }$ implies
that $f_1, \dots , f_s \in A^{\de }$.}

{\bf Proof of (6).}
Consider the $\C_+$-action on $X$ generated by $\de$. Choose
a generic point $\xi \in Z$ and any point $\zeta \in H^{-1}(\xi )$.
Consider the orbit $O_{\zeta}$ of this point $\zeta$ under the action.
Since $h_1, \dots , h_j \in A^{\de }$ the subvariety $H^{-1}(\xi )$
is invariant under the action which means that $O_{\zeta }$
is contained in $H^{-1}(\xi )$. If $F(O_{\zeta })$ is not a point
it cannot meet $T_{\xi }$. But $F^{-1}(T_{\xi })$ is (Zariski)
dense in $H^{-1}(\xi )$. Thus there is a point $\theta$ from
$F^{-1}(T_{\xi })$ in every (Zariski) neighborhood of $\zeta$.
Since $O_{\theta}$ is a point $O_{\zeta}$ is also a point
whence each orbit of the action is contained in a fiber
of $F=(f_1, \dots , f_m)$. This implies the desired conclusion. 
\qed

\medskip

\rem In fact, (6) implies (4) and (5). Indeed,
let $a_1,a_2 \in A$ and either $a_1a_2 \in \Ker \de \setminus \{ 0 \}$ or
$a_1^k+a_2^l \in \Ker \de \setminus \{ 0 \}$ where $k>l\geq 2$ and $(k,l)=1$.
Then $a_1,a_2 \in A^{\de }$. In order to see this we apply (6) in the case
when $F=(a_1,a_2): X \to Y$ where $Y$ is the closure
of $F(X)$ in $\C^2$, and $G: Y \to Z:=G(Y) \subset \C$ is given
either by $(x,y) \to xy$ or $(x,y) \to x^k +y^l$.

More generally, we have the following result which will not
be used in this paper.

\medskip

\cor {\it Let $\de$ be a nonzero locally nilpotent derivation
of $A$, and let $a_1,a_2 \in A$. Suppose that $p \in \C [x,y]$
is a non-constant polynomial which is not equivalent to
a linear one. Let $p(a,b) \in \Ker \de$ and let one of the
following conditions be true

(i) $p(a,b)\notin \C$, or

(ii) $p(a,b)= c \in \C$ and the curve $p^{-1}(c)
\subset \C^2$ does not contain an irreducible contractible
component.

Then $a$ and $b$ belong to $\Ker \de$.}

\medskip

%\defin A semi-degree function $\dg:A\longrightarrow \R\bigcup
%\{-\infty\}$
%on an algebra $A$ is a map which satisfies the following
%axioms:
%(d1) $\dg 0=-\infty,$ and $\dg a\in \R$ for all $a\ne 0;$ $\dg 1=0.$
%(d2) $\dg (f+g)\le \max\{\dg f,\dg g\}$ for all $f,g\in A$.
%(d3) $\dg fg\leq \dg f+\dg g$ for all $f,g\in A.$
%If we have equality in (d3) for all $f,g \in A$ then $\dg$ is
%a degree function.
%A semi-degree function determines
Consider an ascending filtration
$\cF =\{\fr{t}\}$ on
$A,$ where $t \in \R$ and $\fr{t} \subset  \fr{s}$ for $t<s$.
Put $\fro{t} = \bigcup_{s<t} \fr{s}$.

\bigskip

\defin
Consider the linear space $\Gr A=\oplus_{t\in \R} \Gr^t A$
where $\Gr^t A=\fr{t}/\fro{t},$ and introduce the
following multiplication on $\Gr A$. Suppose that
$f_1 \in \fr{t_1}/\fro{t_1}$ and
$f_2 \in \fr{t_2}/\fro{t_2}$. Put
$(f_1 + \fro{t_1}) 
(f_2 + \fro{t_2})$ equal to
$f_1 f_2+ \fro{t_1+t_2}$ if
$f_1 f_2 \in \fr{t_1+t_2} \setminus \fro{t_1+t_2}$ 
and 0 otherwise (of course, the last possibility does not
hold in the case when the filtration is generated
by a degree function). Extend this multiplication
using the distributive law.
Then we call $ \hA = \Gr A$ the {\it associated graded algebra} 
of the filtered algebra $(A,\,\cF )$.

\medskip

\defin Define the mapping $\gr:A\to \Gr A$ by
$\gr f=\hf = f+ \fro{t}$ when
$f \in \fr{t} \setminus \fro{t}$.
If the filtration is generated by a degree function then
this mapping $\gr$ is a multiplicative homomorphism.

\bigskip

\lemma [KaM-L2] {\it Suppose that $\cF$ is a weight filtration on $A$
(see the definition of a weight filtration on an algebra
of regular functions after Lemma 5.2 below).
Then for every derivation $\de$ on $A$ there exists the smallest
$t_0 \in \R$ (which is called the degree of $\de$) such that
$\de (\fr{t} ) \subset \fr{t+t_0}$ for every $t \in \R$.
Furthermore, there exists $a \in \fr{t}$ for some $t$ such that
$a \in \fr{t+t_0} \setminus \fro{t+t_0}$.}

\bigskip

\defin
Consider  the function $\df : A \setminus 0 \to \R \cup \infty$ given by
$$\df(a) = d_A (a) - d_A (\de (a)).$$
Every nonzero $\de \in \LND (A)$
defines a nonzero $\hde \in \LND(A)$ as follows:
$\hde (\ha )= {\widehat {\de (a)}}$ if $\df (a)$
coincides with the negative degree of $\de$, and
$\hde (\ha )=0$ otherwise. We call $\widehat {\de}$
the associate locally nilpotent derivation for $\de$.

\bigskip

\defin The Makar-Limanov invariant of $A$ is
$$\ML (A)=\bigcap_{\de \in \LND(A)} \Ker \de$$
where $\LND(A)$ is the set of all locally nilpotent derivations
on $A$.
%{\footnote{In some other papers (including
%joint papers of the author and Makar-Limanov)
%this invariant is denoted by $\AK (A)$.}
When $\ML (A)= \C$ we call the invariant trivial
(this is so when $A$ is
the ring of polynomials).

\medskip

\rem For a locally nilpotent derivation $\de$ and every
$t \in \C$ the mapping $\exp ({t\de }) : A \to A$ is an automorphism
whence it generates a $\C_+$-action on $X$ [Ren]. When $\de$
is nonzero this action is nontrivial. 
Hence $\ML (A)$ coincides with the subset of $A$
which consists of those regular functions on $X$ that
are invariant under any regular $\C_+$-action.

\medskip

The method of computation of this
invariant which we are going to exploit,
is based on two ideas.

First, there is no need to consider
all the set $\LND (A)$ in the definition of $\ML (A)$.
It is enough to consider its subset $S$ such that it
contains at least one representative from every
equivalence class. Then $\ML (A)=\bigcap_{\de \in S} \Ker \de$.

Second, one can study $\LND (\hA )$ which may be easier
and then one can use the knowledge of $\LND (\hA )$ in order
to find all $\de \in S$. This second step requires
a more convenient description of $\hA$ for some specific
filtrations.

\bigskip

\subsection{The associate algebra $\hA'$}

Let $A'=\C^{[N]}/I'$ where $I'$ is a prime ideal
in the ring of polynomials $\C^{[N]}$ in $N$ variables
$x_1, \dots , x_N$.
For every $a \in A'$ we denote by $[a]$ the set of polynomials
$p \in \C^{[N]}$ such that $p |_{X'} =a$. Each nonzero polynomial $p$
is the sum of nonzero monomials, and the set of this monomials
will be denoted by $M(p)$.

\medskip

\defin 
A {\it weight degree function} on the 
polynomial algebra $\C^{[N]}$ is a degree function $d$
such that $d(p) = \max \{ d(\mu )\, | \, \mu \in M(p) \},$ 
where $p \in \C^{[N]}$ is
a non-zero polynomial. Clearly,
$d$ is uniquely determined by the {\it weights}
$d(x_i) \in \R ,\,\,i=1,\dots,N.$ 
A weight degree function $d$ defines a grading
$\C^{[N]} = \oplus_{t \in \R} \C^{[N]}_{d,\,t},$ where
$\C^{[N]}_{d,\,t} \setminus \{0\}$ consists of all the $d-$homogeneous
polynomials of $d-$degree $t$. Accordingly, 
for any $p \in \C^{[N]} \setminus
\{0\}$ we have a decomposition
$p = h_{t_1}+ \ldots + h_{t_k}$ into a sum of $d-$homogeneous
components $h_{t_i}$ of degree $t_i$ where
$t_1<t_2 < \ldots t_k=d(p)$. We call
$\bp :=h_{d(p)}$ the {\it principal component } of $p$.

\bigskip

\defin Let $d$ be a weight degree function on $\C^{[N]}$.
For $a \in A' \setminus \{0\}$ set 
$$d_{A'}(a) = \inf_{p \in [a]} d(p).$$

\bigskip

Let $\hI_d'$ be the (graded) ideal of $\C^{[N]}$ generated
by the principal components of the elements of $I'$.

\medskip

\lemma [KaM-L2] {\sl For every nonzero $a \in A'$ we have

(1) there exists a polynomial $p \in [a]$ such that $\bp \notin \hI_d'$;

(2) $d_{A'}(a) = d(q)$ for a polynomial
$q \in [a]$ iff
$\bq \notin {\hI_d'}.$ In particular,
$d_{A'}(a) = \min_{q \in [a]}
\{d(q)\} $;

%(3) $d_{A'}$ is a semi-degree function and $d(ab) < d(a) + d(b)$ for $a,b
%\in A'$ only in the case when there exist $p\in [a]$ and $q \in [b]$
%such that $\bp , \bq \notin \hI_d'$ and $\bp \bq \in \hI_d'$.
%In particular, if
(3) if $\hI_d'$ is prime then $d_{A'}$ is a degree function on $A'$.}

\bigskip

Note that $\fr{t}':=\{a\in A'\;\vert\; d_{A'}(a)\le t\}$ where $t \in \R$
gives a filtration on $A'$. We shall call such a filtration
on the algebra $A'$ of regular functions 
of an affine algebraic variety $X'$ a {\it weight filtration}.
This filtration
generates the associate graded algebra $\hA_d'$
and the mapping $\gr_d : A' \to \hA_d'$ which is a
multiplicative homomorphism in the case when
$\hI_d'$ is prime.

\bigskip

\prop [KaM-L2] {\it
The associated graded algebra is isomorphic to
$${\hA_d'} \simeq \C^{[N]}/{\hI_d'} = \C[{\hX_d'}]\,,$$
where ${\hX_d'}$ is the affine variety in $ \C^N$ defined
by the ideal
${\hI_d'}.$
Furthermore, for every nonzero $a \in A'$ we have
$\gr_d (a) = \bp |_{\hX_d' }$
where $p \in [a]$ and $d(p)=d_{A'}(a)$.}

\bigskip

\conv Consider
the coordinate system $(x,y,z,v_1, \dots , v_m,u_{1,1}, \dots
, u_{j,i}$, $\dots )$ in the space $\C^N$
which appeared in Proposition 4.1.
Let $q_0(y,z)=y^k - z^l, m_i, n_{j,i}$ be as in Proposition 4.1.
Put $d_x=d(x), d_y=d(y),d_z=d(z), d_i =d(v_i)$
and $d_{j,i}=d(u_{j,i})$ where $d$ is a weight degree function.
From now on we are going to study only those weight
degree functions on $\C^{[N]}$ that satisfy the following

(1) $kd_y=ld_z$ (in particular, $\bq_0 =q_0 =y^k-z^l$);

(2) $d_1+d_x=kd_y$, and $d_1, d_x$ are $\Q$-independent
(this implies that the elements in the following pairs $(d_x, d_y),
(d_x,d_z),(d_1,d_y),(d_1,d_z)$ are $\Q$-independent);

(3) $d_x <0$ and $d_1 >> d_y >0$;

(4) $d_x +d_{i+1} =m_id_i$ for $i \geq 1$;

(5) $d_x+d_{j,i+1}=n_{j,i}d_{j,i}$ for every $j,i \geq 1$.

\bigskip

\prop {\it
Let $X'$ be the zero set of the system of polynomial
equations from Proposition 4.1 and $A'=\C [X']$.
Then under Convention 5.1 the associate graded
algebra $\hA_d'=\C [\hX_d' ]$ where $\hX_d'$ is isomorphic
to the zero set of the following system
\[ \begin{array}{l}
xv_1-q_0(y,z )=0\\
xv_2-v_1^{n_1}=0\\
\, \, \, \, \, \, \, \, \, \, \, \, \, \, \, \, \, \, \ldots \\
xv_m-v_{m-1}^{n_{m-1}}=0\\
-c_1u_{1,1}=0\\
-c_1u_{1,2}-u_{1,1}^{n_{1,1}} =0\\
\, \, \, \, \, \, \, \, \, \, \, \, \, \, \, \, \, \, \ldots \\
-c_1u_{1,m_1}-u_{1,m_1-1}^{n_{1,m_1-1}} =0\\
-c_2u_{2,1}=0\\
\, \, \, \, \, \, \, \, \, \, \, \, \, \, \, \, \, \, \ldots \\
\end{array} \]
Furthermore, the defining ideal $\hI_d'$ of $\hX_d'$ is prime
(i.e., $\hA_d'$ is a domain)
and it is generated by the left-hand sides of the equations above.}

\proof
By Proposition 4.1 the ideal $I'$ is generated by the
left-hand sides of the equations from that Proposition.
The principal components
of these left-hand sides coincides with the left-hand
side of the system above. For the first equation
it follows from Convention 5.1 (1) and (2).
For the second equation it follows from Convention 5.1 (3)
and the fact that the usual degree of $q_1$ from
Proposition 4.1 with respect to $v_1$ is at most $n_1-1$.
Similarly, taking into consideration the assumption
on the usual degrees of $q_j$ with respect to $v_i$
and $r_{s,j}$ with respect to $u_{s,i}$ we obtain
the claim about principal components.

The polynomial system from this Proposition defines
an algebraic variety $Y$ which is irreducible and a complete
intersection, i.e. its defining ideal $K$ is generated
by the left-hand sides of the equations above.
(Indeed, one can apply, for instance, the same argument
with the Davis theorem which we used while proving that
$X'$ is irreducible and a complete intersection in Proposition 4.1.)
Clearly, $K \subset \hI_d'$. In order to show the reverse
inclusion we need to prove that every element of $\hI_d'$
vanishes on $Y$.

Consider $\C^{N+1}$ which contains $\C^N$ from Proposition 4.1
as a coordinate $N$-plane. Suppose that the coordinates are
$x,y,z,v_i,u_{i,j}$ (as in Proposition 4.1), and $\xi$.
Consider the subvariety $Z \subset \C^{N+1}$ given
by the same system as in Proposition 4.1 and the additional
equation $x\xi -1=0$ (one can see that $Z$ is the localization
of $X'$ with respect to the multiplicative system $\{ x^n | \, n\geq 0
\}$). In this system  of equations
we can replace $xv_1 - q_0(y,z)=0$ by $v_1 - \xi q_0(y,z)=0$,
$xv_2 -v^{n_1}+q_1(y,z,v_1)=0$ by $v_2 -\xi (v^{n_1}+q_1(y,z,v_1))=0$,
etc.. Clearly, the defining ideal $J$ of $Z$ is generated
by $x\xi -1=0$ and these replacements. More precisely,
every element of $J$ is of form
$$\alpha_0 (x\xi -1)+\alpha_1 (v_1 - \xi q_0(y,z)) +
\alpha_2 (v_2 -\xi (v^{n_1}+q_1(y,z,v_1))) + \dots .$$
Extend the weight degree function $d$ to $\C^{[N+1]}$ by putting
$d(\xi ) = -d_x$. Then the above form of elements of $J$
and Convention 5.1 imply that $\hJ_d$ is generated by
$x\xi -1, xv_1 - q_0(y,z), 
v_2 -\xi v^{n_1}$, etc., where $\hJ_d$ is the ideal generated
by the principal components of elements of $J$. 
Hence $\hJ_d =K[\xi ] = K[1/x ]$. This means that $\hJ_d$
defines a variety $\hZ \subset \C^{N+1}$ such that the image
of its natural projection to $\C^N$ is $Y\setminus \{ x=0 \}$.
Since $I' \subset J$ we see that $\hI_d' \subset \hJ_d$ whence
every element of $\hI_d'$ vanishes on $\hZ$ and, therefore, on $Y$.
\qed

\bigskip

\rem The variety $\hX_d'$ is independent on the
choice of $d$ satisfying Convention 5.1 and it is isomorphic
to the zero set of the following polynomial equations
in the space $\C^{3+m}$ with coordinates
$(x,y,z,v_1, \dots ,v_m)$
\[ \begin{array}{l}
P_1(x,y,z,v_1)=xv_1-q_0(y,z )=0\\
P_2(x,v_1,v_2)= xv_2-v_1^{n_1}=0\\
\, \, \, \, \, \, \, \, \, \, \, \, \, \, \, \, \, \, \ldots \\
P_m(x,v_{m-1},v_m)= xv_m-v_{m-1}^{n_{m-1}}=0.\\
\end{array} \]
Therefore, we shall write further $\hI' ,\hA'$, and $\hX'$ instead of
$\hI_d' , \hA_d'$, and $\hX_d'$ provided it does not cause
misunderstanding.

\bigskip

\subsection{ Locally nilpotent derivation of Jacobian type}

We shall study locally nilpotent derivations on the
associate graded algebra $\hA' = \C [\hX' ]$ 
which was introduced in the previous subsection.
We shall use the presentation of $\hX'$ given in Remark 5.3.
In particular, $P_1, \dots , P_m$ have the same meaning as
in that Remark.

For polynomial $q_1, \dots , q_{m+3}$ on $\C^{m+3}$ we denote
by $J(q_1, \dots , q_{m+3})$ their Jacobian with respect
to $(x,y,z,v_1, \dots ,v_m)$. 
In this subsection and the next one we denote $q|_{\hX'}$
by $\tq$ for every polynomial $q$ on $\C^{m+3}$.
Note that $(\tx , \ty , \tz )$ is a local (holomorphic) coordinate
system at each point of $\hX_0 = \hX' \setminus \{ x=0 \}$.
For $a_1,a_2,a_3 \in \hA'$ we denote by $J_0(a_1,a_2,a_3)$
the Jacobian of these regular functions on $\hX_0$
with respect to $\tx , \ty$, and $\tz $.
This is a rational function on $\hX'$ but
$\tx^m J_0(a_1,a_2,a_3)$ is already regular on $\hX'$
since $x^m$ is the determinant of the matrix
$\{ \de P_i/\de v_j \, | \, i,j=1, \dots , m \}$.
Furthermore, if $a_i =\tq_i$ then
$J(P_1, \dots , P_m, q_1,q_2, q_3)|_{\hX '} =
\tx^m J_0(a_1,a_2,a_3)$ [KaM-L2] (it is useful to keep this
equality in mind when we shall make direct computations in Proposition
5.3).

\bigskip

\defin Fix $a_1, a_2 \in \hA'$ and let $a\in \hA'$.
Then $\de (a) = \tx^mJ_0(a_1,a_2,a)$ is called
a derivation of Jacobian type on $\hA'$.

\bigskip

\lemma [KaM-L2] {\it Let $\delta$ be a nontrivial locally nilpotent
derivation on $\hA'$ and let $a_1, a_2 \in \Ker \delta$ be
algebraically independent. Then $\de (a) = \tx^mJ_0(a_1,a_2,a)$
is a locally nilpotent derivation which is equivalent to $\delta$.}

\bigskip

We say that $a \in \hA'$ is $d$-homogeneous if $a$ is
the restriction to $\hX'$ of a $d$-homogeneous polynomial.
We are going to study locally nilpotent derivations of
Jacobian type on $\hA'$ such that $a_1$ and $a_2$ in the above
definition are irreducible $d$-homogeneous
where $d$ is a weight degree function satisfying
Convention 5.1. 

\bigskip

\lemma {\it
Let $a \in \hA'$ be an irreducible $d$-homogeneous element.
Then up to a constant factor
$a$ is of one of the following elements
$\tv_i, \tx, \ty,\tz ,$ or $\ty^k +c \tz^l$ where $c \in \C^*$ 
and $k, l$ are the
same as in Proposition 4.1.}

\proof
Let $q$ be $d$-homogeneous and $a=\tq$ (in particular,
$q$ is irreducible).
It follows from the explicit form of the polynomial system
in Remark 5.3 and Convention 5.1
that we can suppose that each monomial from
$M(q)$ is non-divisible by $xv_i$ for every $i=1, \dots ,m$.
The restriction of every function $v_i$ to $\hX'$ is
of form $q_0^s/x^j$ where $s,j>0$. Note that if
we extend $d$ naturally to the field of rational functions
then $d(v_i)=d(q_0^s/x^j)$ in virtue of Convention 5.1. 

Assume that $\mu_1 , \mu_2 \in M(q)$
are such that $\mu_1$ is divisible by $x$ but $\mu_1$ is not.
Then $\mu_1$ and $\mu_2$ coincides with the restriction to $\hX'$
of the functions $x^{j_1}y^{\alpha_1}z^{\beta_1}$ and
and $y^{\alpha_2}z^{\beta_2}q_0^s/x^{j_2}$ where $j_1 >0, j_2 \geq 0$.
Since $d(\mu_1 )=d(\mu_2 )$ we have
$d(x^{j_1}y^{\alpha_1}z^{\beta_1})=
d(y^{\alpha_2}z^{\beta_2}q_0^s/x^{j_2})$ whence
$(j_1+j_2)d_x =d(y^{\alpha_2-\alpha_1}
z^{\beta_2-\beta_1}q_0^s)$. Since 
$d_y= (l/k)d_z$ and $d(q_0)=kd_y$ we see that $d_x$ and $d_y$ are
$\Q$-dependent which contradicts Convention 5.1. 

Thus the assumption is wrong and if one monomial from
$M(q)$ is divisible by $x$ then every monomial is.
Therefore, we suppose that none of the monomials
from $M(q)$ is divisible by $x$ since we are interested
in the case when $a$ is irreducible. Let 
$\mu_1 , \mu_2 \in M(q)$ and $\mu_i =y^{\alpha_i}z^{\beta_i}\nu_i$
where $\nu_i$ is a monomial which depends on $v_1, \dots , v_m$ only.
The restriction of $\mu_i$ to $\hX'$ coincides
with $y^{\alpha_i}z^{\beta_i}q_0^{s_i}/x^{j_i}$. Thus
$$d(y^{\alpha_1}z^{\beta_1}q_0^{s_1}/x^{j_1})=
d(y^{\alpha_2}z^{\beta_2}q_0^{s_2}/x^{j_2}).$$
The same argument as above shows that $j_1=j_2$
since otherwise $d_x$ and $d_y$ are $\Q$-dependent.
Hence
$$d(y^{\alpha_1}z^{\beta_1}q_0^{s_1})=
d(y^{\alpha_2}z^{\beta_2}q_0^{s_2}).$$
Since $d(q_0)=kd_y =ld_z$ and $(k,l)=1$ we see that
$\alpha_i =\alpha_0 +t_ik$ and
$\beta_i = \beta_0 + \tau_il$ where
$\alpha_0$ is one of the numbers $0,1, \dots , k-1$,
$\beta_0$ is one of the numbers $0,1, \dots , l-1$ and
$t_1-t_2 +\tau_1 -\tau_2 =s_2-s_1$.

Therefore, the restriction of $q$ to $\hX'$
coincides with 
$$y^{\alpha_0}z^{\beta_0}\varphi (y^k,z^l) q_0^s/x^j$$
where $\varphi (y^k,z^l)$ is $d$-homogeneous and
the restriction of $q_0^s/x^j$ to $\hX'$
coincides with the restriction of some monomial $\nu$
which depends on $v_1, \dots , v_m$ only.
Now the statement of Lemma follows from the fact
that that $\varphi (y^k,z^l)$ is
the product of factors of type $c_1y^k +c_2z^l$ where $c_1,c_2 \in C$.
\qed

\bigskip

\cor {\it
Let $a=\tq$ where $q\notin \C [y,z]$ is a $d$-homogeneous
polynomial which does not depend on $x$. Then $q$ is divisible
by some $v_i$.}

\bigskip

\prop {\it Let $X'$ be as in Proposition 4.1.
Suppose also that $X'$ is smooth.
Let $\de (a) = \tx^m J_0(a_1, a_2,a)$ be a nontrivial
locally nilpotent derivation of Jacobian type on $\hA'$ and
let $a_1$ and $a_2$ be irreducible $d$-homogeneous. Suppose that
$m\geq 2$. Then 

(1) $(a_1,a_2)$ coincides (up to the order) with
one of the pairs $(\tx , \ty)$ or $(\tx , \tz)$,

(2) $\tx \in \Ker \de $ and $\dg_{\de } (v_i) \geq 2$ for 
every $i=1, \dots , m$.}

\proof
If $(a_1, a_2)$ is one of the pairs in (1)
it is easy to check that $\de$ is nontrivial and locally nilpotent,
and (2) holds also. We need to show that if we
use other possible irreducible $d$-homogeneous elements as
$a_1,a_2$ (recall that such elements were described in Lemma 5.4)
then $\de$ is not a nontrivial locally nilpotent derivation.
Since we want a nontrivial derivation we need to consider only
$a_1$ and $a_2$ which are algebraically independent in $\hA'$.

Case 1. Let $(a_1,a_2)=(\ty ,\tz )$. The direct computation shows that
$\de (\tx )= \tx^m$ whence $\de$ cannot be locally nilpotent.
Indeed, one can see that $\dg_{\de } (\de (\tx ))=
\dg_{\de } (\tx ) -1$. On the other hand
$\dg_{\de } (\tx^m )= m \dg_{\de } (\tx )$ which yields
a contradiction.

Case 2. Either $a_1$ or $a_2$ is of form $\ty^k +c\tz^l$ where $c\in \C^*$
and $k$ and $l$ are as in Proposition 4.1.
By Theorem 5.1 $\ty ,\tz \in \Ker \de$ since
$\ty^k +c\tz^l \in \Ker \de$.
By Lemma 5.3 $\de$ is equivalent to the derivation
$\tx^m J_0(\ty , \tz ,a)$ whence this case does not hold.

Case 3. Let $(a_1,a_2)=(\tv_{i_1},\tv_{i_2})$ where $i_1 <i_2$.
Consider the identical morphism $F: \hX' \to \hX' \subset \C^{m+3}$ and
morphism $G : \hX' \to \C^2$ given by
$(\tx , \ty , \tz , \tv_1 , \dots , \tv_m) \to (\tv_{i_1}, \tv_{i_2})$. 
Recall that $\tv_{i_k} =\tq_0^{s_k}/\tx^{j_k}$.
It is easy to check that $\tv_{i_1}$ and $\tv_{i_2}$ are
algebraically independent in $\hA'$ which means that the pairs
$(s_1,j_1)$ and $(s_2, j_2)$ are not proportional (in fact,
it can be shown by induction that $j_2s_1 -j_1s_2 >0$).
Consider a generic point $\xi \in \C^2$ and the fiber
$G^{-1} (\xi )$. Since 
$(s_1,j_1)$ and $(s_2, j_2)$ are not proportional
one can see that each component of this fiber is a curve
in $\C^{m+3}$ given by equations $v_i=c_i,x=c'$, and
$q_0(y,z) =y^k-z^l=c$ where $c_i,c' \in \C$ and $c \in \C^*$.
This curve is hyperbolic and thus it does not admit non-constant
morphisms from $\C$. By Theorem 5.1 if
$\de$ is locally nilpotent then
$\tx , \ty , \tz , \tv_i \in \Ker \de$ whence $\de$
is trivial.
Therefore, this case does not hold.

Case 4. Let $(a_1, a_2) = (\tx , \tv_i )$. The same argument
is in Case 3 works.
 
Case 5. Let $a_1, a_2)=(\ty , \tv_i)$.
Consider again the identical morphism $F: \hX' \to \hX'$ and
morphism $G : \hX' \to \C^2$ given by
$(\tx , \ty , \tz , \tv_1 , \dots , \tv_m) \to (\ty , \tv_{i})$. 
Recall that $\tv_i$ is of form
$\tv_{i} =\tq_0^{s}/\tx^{j}$ where $j\geq 2$ if $i>1$.
This implies that the curve $G^{-1}(\xi )$ where 
$\xi =(c_1,c_2) \in \C^2$
is isomorphic to the curve $(c_1^k-z^l)^s -c_2x^j =0$.
When $j\geq 2$ and $s$ is not divisible by $j$
the last curve does not have contractible 
components for generic $\xi$.
Theorem 5.1 (6) implies that if in this case $\de$ is locally nilpotent
then $\de$ must be trivial. If $j \geq 2$ and $s$ is divisible
by $j$ then each irreducible component of $G^{-1}(\xi )$
is contractible and contains double points of 
$G^{-1}(\xi )$. Since $G^{-1}(\xi ) \subset \hX'$ is invariant under
the $\C_+$-action generated by the locally nilpotent
derivation $\de$ the singular points of
this curve must be fixed under this action whence the
action itself is trivial on $G^{-1}(\xi )$. Therefore,
it is trivial on $\hX'$ whence $\de$ is again trivial.

It remains to consider the case when $j=1$,
i.e. $(a_1,a_2)=(\ty , \tv_1) = (\ty ,\tq_0/\tx )$. The direct computation
shows that $\de (\tx )$ coincides (up to a constant multiple) with
$\tx^{m-1}\tz^{l-1}$. Since $m\geq 2$ we see that $\de$
cannot be nontrivial locally nilpotent (indeed, compare
$\dg_{\de } (\tx )$ and $\dg_{\de } (\de (\tx ))$)
and we have to disregard this case.

Case 6. When $(a_1,a_2)=(\tz ,\tv_i )$ the same argument as in
Case 5 works.
\qed

\bigskip

\subsection{ The computation of $\ML (A')$}

By Definition 5.4 
each nontrivial $\de \in \LND (A')$ generates a nontrivial
$\hde_d \in \LND (\hA' )$ which depends on the choice of the
weight degree function $d$. Similarly the mapping
$\gr_d : A' \to \hA'$ depends on $d$. The following
relation between $\de$ and $\hde_d$ is simple but essential for us
$${\rm if} \, \, \dg_{\de } =s \, \, {\rm then} \, \,
\dg_{\hde_d}(\gr_d (a) ) (a) \leq s.$$
\defin  A locally nilpotent derivation $\de$ on $A'$
is called perfect if its
associate derivation $\hde_d$ is of form
$\hde_d (a) =\tx^m J_0(a_1, a_2,a)$ where $a_1,a_2\in \hA_d'$ are
irreducible $d$-homogeneous and algebraically independent.
The set of all perfect derivations will
be denoted by $\de \in \Per (A')$.

\bigskip

\prop [KaM-L2] {\it If the associate graded algebra
$\hA_d'$ is a domain (i.e., $\hI_d'$ is prime) then for
every nontrivial locally nilpotent derivation on $\hA'$
there exists an equivalent perfect derivation. In particular,
$$\ML (A' )= \bigcap_{\de \in \Per (A' )}\Ker \de .$$}

\bigskip

\prop {\it Let $A'$ be as in Proposition 4.1 and let
$d$ satisfy Convention 5.1. Suppose that $X'$ is smooth.
For every $\de \in \LND (A')$ we have $x \in \Ker \de$ whence
$\ML (A' ) \ne \C$.}

\proof Since $A'$ is a domain by Proposition 5.2,
it is enough to consider $\de \in \Per (A')$ by Proposition 5.4.
Consider $a \in A'$ with $\dg_{\de }(a) \leq 1$. 
Show that there exist a polynomial $q$ with $q|_{X'} =a$
such that none of monomial $\mu \in M(q)$ is divisible
by $v_i$ or $u_{s,j}$ for all $i,s,j$
(i.e. $q \in \C [x,y,z]$). It follows from
the explicit from of the polynomial system in Proposition 4.1
that we can suppose that none of $\mu \in M(q)$ is divisible
by $xv_i$ or $xu_{s,j}$. Thus $M(q)$ is the disjoint
union $M_1(q) \cup M_2(q)$ where $M_1(q)$ consists of monomials
which depends on $x,y,z$ only and $M_2(q)$ consists of monomials
which do not depend on $x$ and do not belong to $\C [y,z]$.
Show that $M_2(q)$ is empty. Let $\mu \in M_2(q)$. 
Under Convention 5.1 one can keep $d_y,d_z$ fixed, decrease $d_x$,
and increase $d_i,d_{s,j}$ so that $d(\mu ) > d(\nu )$
for every $\nu \in M_1(q)$. Hence if $\bq_d$ is 
the principal component of $q$ then $M(\bq_d ) \subset M_2(q)$.
The relation between $\de$ and $\hde_{d}$ shows that
$\dg_{\hde_d}(\gr_d (a) ) \leq 1$. The element
$\gr_d (a) =\bq_d|_{\hX'} $ is the product of irreducible
$d$-homogeneous elements of $\hA'$. By Corollary 5.2
one of these elements is $\tv_i$ whence
$\dg_{\hde_d} (\tv_i) \leq \dg_{\hde_d}(\gr_d (a) ) \leq 1$
which contradicts Proposition 5.3. Thus $M_2(q)$ is empty.
Let $b \in A'$ with $\dg_{\de } (b)=1$. Recall that
by Theorem 5.1 there exist $a', a_0, \dots , a_s \in \Ker \de$
such that $a' \tv_1 = \sum_{j=0}^s a_jb^j$ where $s=\dg_{\de }\tv_i$.
Hence $\tv_1 = (q(x,y,z)/r(x,y,z))|_{X'}$ where $a'=r(x,y,z)|_{X'}$.
On the other hand $A' \subset \C [x,y,z, 1/f(x)]$ where $f$ is as
in Proposition 4.1.
Since $v_1 \notin \C [x,y,z]$ we see that $r(x,y,z)$ must be
divisible by some $x-c$ where $c$ is a root of $f$.
Therefore, $x-c \in \Ker \de$ as a divisor of an element
from $\Ker \de$ whence $x \in \Ker \de$.
\qed

\medskip

This implies

\medskip

{\bf Lemma II.}
{\it Let $X'$ be an affine algebraic variety of dimension 3
such that $X'$ is a UFD, all invertible functions
on $X'$ are constants, and
let the assumption ($1'$) and ($2'$) of Lemma I hold,
but (3) does not.
Then $X'$ is an exotic algebraic structure on $\C^3$
with a nontrivial Makar-Limanov invariant.}      

\medskip

\rem
In [KaZa1] we described some conditions under which one can extend
locally nilpotent derivations from $A$ to $A'$. Using this
technique it is not difficult to show that $\ML [A']= \C [x]|_{X'}$.

\medskip

Lemmas I and II imply

\medskip

\thm {\it
Let $X'$ be an affine algebraic variety of dimension 3
such that $X'$ is a UFD ,
all invertible functions
on $X'$ are constants, and

(1) the Euler characteristic of $X'$ is
$e(X')=1$; 

($2'$) there exists a Zariski open subset $Z$ of $X'$
which is a $\C^2$-cylinder over a curve $U$; 

(3) each irreducible component of $X' \setminus Z$ is a UFD.

Then $U$ is isomorphic to a Zariski open subset
of $\C$ and $p$ can be extended to a regular function on $X'$. Furthermore, 
$X'$ is isomorphic to $\C^3$ and $p$ is a variable.

The same conclusion remains true if we replace (1) and (3) by

($1'$) $X'$ is smooth and $H_3(X')=0$;

($3'$) each irreducible component of $X' \setminus Z$
has at most isolated singularities. 

\medskip

In the case when conditions ($1'$) and ($2'$) hold but (3) does not,
$X'$ is an exotic algebraic structure on $\C^3$ (that is,
$X'$ is diffeomorphic to $\R^6$ as a real manifold
but not isomorphic to $\C^3$) with a nontrivial
Makar-Limanov invariant.}

In order to finish the proof of the Main Theorem we need now to show
that condition ($2'$) above is equivalent to condition
(2) in the Main Theorem. This will be done in [KaZa2].

\bigskip

\address
{\noindent Department of Mathematics and Computer Science \\
University of Miami \\
Coral Gables, FL  33124, U.S.A.\\
\email kaliman@@cs.miami.edu}

\end{document}